\documentclass{IEEE-Con-Sys-mag}
\makeatletter
\@ifundefined{mv@sans}{\DeclareMathVersion{sans}}{}
\makeatother
\usepackage{xcolor}
\definecolor{iceblue}{RGB}{230,242,247}
\definecolor{ivoryA}{HTML}{FCFAF2}  
\usepackage{tcolorbox}
\usepackage[compress]{cite}
\usepackage{amsmath}

\usepackage{amssymb}
\usepackage{amsfonts}
\newcommand*{\QEDB}{\hfill\ensuremath{\square}}
\newcommand{\tcb}{\textcolor{black}}
\newcommand{\tcs}{\textcolor{black}}
\newcommand{\tcr}{\textcolor{black}}
\newcommand{\tcrr}{\textcolor{red}}
\begin{document}

\title{Hybrid Set-Seeking Systems\stitle{Model-Free Feedback Optimization via Hybrid Inclusions}}

\author{\vspace{-0.5cm}{J}ORGE I. POVEDA, ANDREW R. TEEL}
\vspace{0.2cm}

\affil{\vspace{0.5cm}J. I. Poveda (\tcrr{poveda@ucsd.edu}) is with the Deptartment of Electrical and Computer Engineering, University of California, San Diego, USA\\
A. R. Teel (\tcrr{teel@ucsb.edu}) is with Department of Electrical and Computer Engineering, University of California, Santa Barbara, USA}

\maketitle

\dois{}{}

\chapterinitial{H}\tcs{ybrid dynamical systems, which combine continuous-time and discrete-time dynamics, are pervasive in modern engineering. They naturally arise in cyber-physical systems (CPS) that integrate physical processes with digital components, spanning domains such as transportation \cite{farid2016hybrid}, power grids \cite{ochoa2023control}, robotics \cite{ames2017hybrid}, and power electronics \cite{senesky2003hybrid}. Beyond cyber-physical systems, hybrid dynamics also appear in a wide range of applications, including mechanical systems with impacts~\cite{mechanicalImpact}, chemical process control~\cite{lennartson2002hybrid}, manufacturing systems~\cite{pepyne2002optimal}, and synchronization mechanisms in biological systems~\cite[Ex. 1.2]{bookHDS}. In the context of control and decision-making, hybrid systems can model algorithms that involve logic modes, timers, resets, impulsive actions, sampled-data structures, switching dynamics, event-triggered mechanisms, and general ``if–then'' rules, \cite{bookHDS,RSanfeliceBook,lygeros2008hybrid}. Hybrid systems are also closely related to frameworks widely studied in computer science and engineering, such as the well-known hybrid automata \cite{henzinger1996theory,lygeros2003dynamical,tavernini1987differential}}.

Despite the prevalence of hybrid systems in engineering and science, the development of extremum-seeking (ES) methods for systems with hybrid dynamics remained relatively unexplored until recently. ES systems are feedback control mechanisms designed to achieve real-time steady-state optimization in dynamical systems where both the plant dynamics and the cost function are unknown. This "model-free" characteristic makes ES an adaptable and powerful feedback control scheme for various complex engineering applications where precise mathematical models are unavailable. As a result, ES has become widely adopted across different industries and areas, including engine optimization \cite{TeelPopovic2001}, wind turbine control \cite{ESCwind}, underwater vehicles \cite{PID_ES,Fish_ESC}, ABS brake systems \cite{ESC_ABS}, autotuning systems \cite{benosman_ijrnc16}, photolithography in semiconductor manufacturing \cite{ren2012laser}, fuel cells \cite{ESCfuelcells}, control of mobile robots \cite{PatentBenosmanPoveda},  HVAC systems \cite{PatentJhonson}, model predictive control \cite{SB16}, biological systems \cite{abdelgalil2022sea}, traffic control \cite{ochoa2022high,yu2021extremum}, mechanical systems \cite{elgohary2025extremum,suttner2022extremum,suttner2020extremum}, sampled-data systems \cite{SampledDataRevisited,zhu2022sampled}, etc. For surveys on ES, we refer the reader to \cite{scheinker2024100,surveyESCTan,sternby1979review}, as well as to \cite[Sec. 13.3]{astrom1994adaptive}.

The idea of using feedback control mechanisms for real-time black-box optimization in dynamical systems dates back to the early 1920s with LeBlanc's work \cite{Leblanc}, and was further advanced in the 1950s \cite{Morosanov} and 1960s \cite{Meerkov,Meerkov1,krasovskii1963dynamics,chinaev1969self} in the adaptive controls' literature under the name of ``extremum control'' \cite[Sec. 13.3]{Astrom:Book}. However, most of the stability results for ES nonlinear systems emerged in the 1990s and 2000s, utilizing multi-time scale methods for ODEs and identifying the appropriate time scales associated with the components of the algorithms \cite{Krstic2000}. Since then, different analytical tools for ES have been developed to solve model-free optimization and decision-making problems in single-agent systems \cite{KrsticBookESC,BookESCZhang,BookES_Stabilization,DurrLieBracket,MoaseHamerstein,Benosman:book}, multi-agent systems \cite{Frihauf12a,Dougherty:16,PovedaQuijanoACC13}, PDEs \cite{ThiagoPDE18}, stochastic systems \cite{LiuStochastic}, and more.

As the foundational theory behind ES relied on perturbation tools for locally Lipschitz ordinary differential equations (ODEs), particularly averaging \cite{Verhulst_book,KhalilBook} and singular perturbations \cite{Kokotovic_SP_Book,TeelMoreauNesic2003,KhalilBook}, the majority of the ES systems and algorithms were designed for applications with closed-loop smooth dynamics. In particular, the study of ES systems that incorporate discrete-time behaviors---such as ``events'' triggered by specific conditions, abrupt state ``jumps'' or ``resets'' aimed at enhancing performance, or algorithms that ``switch'' between different dynamics during the seeking process---has received comparatively less attention. Such systems, which are naturally hybrid, also arise across different applications where either the plant is a  hybrid system or the controller must be hybrid to address fundamental limitations of traditional smooth control techniques. These limitations can be related to achieving desired transient or steady-state performance or handling complexities such as dynamic communication networks, resource-constrained actuators, occluded operational spaces, or dynamic adversarial environments. ES systems with discrete states, such as logic modes, also arise when general ``if-then'' rules are incorporated into the decision-making algorithms, which is now a standard practice in the design and implementation of controllers and optimization algorithms for autonomous and \cite{Future_control} and cyber-physical systems \cite{SanfeliceCPS}.   
Despite the limited available theoretical tools for their analysis and design, ES algorithms with logic states and switching rules have been explored for particular applications since the early 1970s \cite{luxat1971stability}, when a switching hysteresis-based mechanism was designed to stabilize a peak-holding adaptive controller in electrical machines with neglected dynamics. Modern approaches of ES that incorporate discrete-time behaviors via switches or resets have also been studied  for the efficient control of photovoltaic power systems \cite{moura2013lyapunov} and HVAC systems in industry \cite{PatentJhonson}. However, a comprehensive theory for the study of hybrid ES systems has lagged behind their practical applications, despite its potential to enable the design of novel feedback architectures and techniques capable of addressing a broader class of problems than those handled by traditional ES methods. 

As discussed later in the paper, ES systems with hybrid dynamics
typically exhibit stability properties with respect to general \emph{sets} rather
than isolated equilibria. We therefore refer to these controllers as \emph{hybrid set-seeking systems}, which are generally characterized by non-smooth dynamics modeled through differential and difference inclusions rather than standard differential or difference equations. Building on this perspective, the primary aim of this paper is to provide an accessible introduction to non-smooth and hybrid set-seeking systems formulated as hybrid inclusions. We begin with well-established smooth ES algorithms extensively studied in the literature \cite{Krstic2000,NesicChina,DerivativesESC,TanAndNesic2006Local,mimmo2024uniform}, which can be viewed as hybrid set-seeking schemes with negligible discrete-time dynamics. We then extend this perspective by incorporating hybrid dynamics through averaging tools, enabling the solution of complex model-free optimization and decision-making problems in static maps. Next, we show how singular perturbation methods for hybrid systems can address seeking problems in plants governed by dynamical systems rather than static maps. Taken together, these developments yield perturbation-based tools for multi-time-scale hybrid inclusions, providing a unified framework for analyzing and designing hybrid seeking algorithms that overcome the limitations of traditional smooth ES methods. The effectiveness of this framework is demonstrated through diverse applications and numerical examples, showcasing the flexibility and robustness of the proposed tools.
\subsection{Perturbation Theory and Seeking Systems}
The key mathematical property that enables the design and analysis of most ES
algorithms modeled as oscillatory, locally Lipschitz ODEs
of the form
\begin{equation}\label{oscillatingODE}
\dot{x}=f(x,\omega t),\quad x\in\mathbb{R}^n,\quad t\in\mathbb{R}_{\geq 0},
\end{equation}
is that, for a sufficiently large oscillation frequency $\omega$, the solutions of
\eqref{oscillatingODE} can be approximated on finite intervals by those of the so-called
\emph{average dynamics}:
\begin{equation}\label{integralaverage}
\dot{\bar{x}}=f_{\text{ave}}(\bar{x})
:=\frac{1}{T_{\omega}}\int_{0}^{T_{\omega}} f(\bar{x},\omega t)\,dt,
\quad \bar{x}\in\mathbb{R}^n,
\end{equation}
where $T_{\omega}=2\pi/\omega$ is the oscillation period of $f(x,\cdot)$. For example, to minimize an unknown smooth scalar static map $J:\mathbb{R}\to\mathbb{R}$ representing the plant, the simplest ES continuous-time algorithm based on periodic probing, studied in \cite{Krstic2000,DerivativesESC,TanAndNesic2006Local} can be modeled by the following equations:
\begin{subequations}\label{ESCODEVanilla}
\begin{align}
u&=\hat{u}+\varepsilon_a\sin(\omega t),\\
\dot{\hat{u}}&=-\frac{2}{\varepsilon_a} J(u)\sin(\omega t),
\end{align}
\end{subequations}
where $\varepsilon_a>0$ and $\omega>0$ are tunable parameters. A Taylor expansion of $J$ around $\hat{u}$ for small values of $\varepsilon_a$, and a direct computation of \eqref{integralaverage} using the identities $\int_{0}^{T_{\omega}}\sin(\omega t)dt=0$ and $\frac{1}{T_{\omega}}\int_{0}^{T_{\omega}}\sin(\omega t)^2dt=\frac{1}{2}$ show that the average dynamics of \eqref{ESCODEVanilla} corresponds to the following $\varepsilon_a$-perturbed system, with state $\bar{u}\in\mathbb{R}$:
\begin{equation}\label{perturbedgradientflow01}
\dot{\bar{u}}=f_{\text{ave}}(\bar{u})=-\frac{\partial J(\bar{u})}{\partial \bar{u}}+O(\varepsilon_a),
\end{equation}
where $O(\varepsilon_a)$ indicates bounded perturbations (on compact sets) of order $\varepsilon_a$. In turn, as $\varepsilon_a\to0^+$, the behavior of system \eqref{perturbedgradientflow01} can be approximated by a standard gradient flow of $J$, with state $\tilde{u}\in\mathbb{R}$ and dynamics
\begin{equation}\label{gradientflow01}
\dot{\tilde{u}}=-\frac{\partial J(\tilde{u})}{\partial \tilde{u}}.
\end{equation}
Under suitable regularity (e.g., locally Lipschitz) properties in system \eqref{gradientflow01}, we can use standard robustness results for differential equations (see ``Closeness of solutions for Perturbed Systems'') to assert that the solutions of \eqref{perturbedgradientflow01} and \eqref{gradientflow01} are $(T,\epsilon)$-close when $\varepsilon_a$ is sufficiently small \cite{KhalilBook,TeelPraly}. Namely, for each compact set of initial conditions, each $T>0$, and each $\epsilon>0$ there exists $\varepsilon_a^*>0$ such that for all $\varepsilon_a\in(0,\varepsilon^*_a)$ the solutions of \eqref{perturbedgradientflow01} and \eqref{gradientflow01} satisfy the bound
\begin{equation}\label{closeness1}
\sup_{t\in[0,T]}|\tilde{u}(t)-\bar{u}(t)|\leq \frac{\epsilon}{2}.
\end{equation}
Similarly, by leveraging results on averaging theory for differential equations \cite{TeelMoreauNesic2003} (see "Averaging Theory: From ODEs to Hybrid Systems"), it can be shown that, for a fixed and sufficiently small $\varepsilon_a>0$, the solutions of \eqref{ESCODEVanilla} and \eqref{perturbedgradientflow01} are also $(T,\epsilon)$-close when $\omega$ is sufficiently large, i.e., for each compact set of initial conditions, each $T>0$, and each $\epsilon>0$ there exists a sufficiently large frequency $\omega^*>0$ such that for all $\omega>\omega^*$ the solutions of \eqref{ESCODEVanilla} and \eqref{perturbedgradientflow01} satisfy the bound
\begin{equation}\label{closeness2}
\sup_{t\in[0,T]}|\hat{u}(t)-\bar{u}(t)|\leq \frac{\epsilon}{2}.
\end{equation}
By combining the bounds \eqref{closeness1} and \eqref{closeness2}, we can conclude that the trajectories of the ES dynamics \eqref{ESCODEVanilla} and those of the gradient flow \eqref{gradientflow01} are also $(T,\epsilon)$-close in compact sets. This property immediately allows us to obtain qualitative attributes for the behavior of the solutions of the model-free ES dynamics \eqref{ESCODEVanilla} based on the behavior of the solutions of the model-based gradient flow \eqref{gradientflow01}, provided the parameters $(\omega,\varepsilon_a)$ are appropriately tuned. Among such attributes, we can infer stability properties for system \eqref{ESCODEVanilla} based on the stability properties of \eqref{gradientflow01}, provided that such stability properties are ``robust'' to small perturbations, an attribute that, under the Lipschitz continuity of \eqref{gradientflow01}, can be readily established using, for example, Lyapunov functions. This perturbation-based methodology for the study of smooth ES lies at the heart of the design and stability analysis of most algorithms, including those whose average dynamics \eqref{integralaverage} are not simple perturbed gradient flows, but instead follow alternative well-posed optimization dynamics, such as primal-dual dynamics \cite{PovedaNaliAuto20}, projected gradient flows \cite{chen2025continuous}, Riemannian gradient flows \cite{PoQu12}, Newton flows \cite{MISONewton}, pseudo-gradient flows for games \cite{Frihauf12a}, decentralized optimization dynamics \cite{ye2023distributed}, etc, see \cite{DerivativesESC} for a unifying framework.

Fortunately, the perturbation methodology described above extends to systems with hybrid dynamics under appropriate, mild regularity conditions and with an appropriate measure of closeness. For the latter, we  note that, unlike smooth ODEs, solutions to hybrid systems are typically discontinuous, making the use of the uniform distance, as in \eqref{closeness1}-\eqref{closeness2}, to measure the closeness between solutions of nominal and perturbed systems problematic. For instance, if \( \hat{u} \) is a solution of a hybrid system that has a discontinuity at time \( t' \), and \( \bar{u} \) is a solution of a \(\varepsilon_a\)-perturbed version of the same system that has a discontinuity at time \( t' + \varepsilon_a \), then no matter how small \(\varepsilon_a\) is, \( \hat{u} \) and \( \bar{u} \) will never be \((T, \epsilon)\)-close at time \( t' \), regardless of the choice of \(\varepsilon_a\). Thus, a different measure of closeness is required to link the qualitative behaviors of the perturbed and unperturbed systems. Although this issue is well-known in the hybrid systems literature \cite[Example 3]{CSM_hybrid}, we will illustrate this phenomenon in a stylized hybrid set-seeking system that exhibits "jumps" in the states, resulting in Zeno-like behavior reminiscent of a bouncing ball without air resistance. %By using tools from hybrid inclusions and a generalized notion of solution \cite{bookHDS}, we will also show how the perturbation-based methodology can be naturally extended to ES systems that incorporate discrete-time behaviors during the seeking process.

\subsection{Generalized Solutions for Hybrid Seeking Systems}
To extend the perturbation-based methodology to seeking systems with hybrid dynamics, we must generalize both the mathematical models we consider and their corresponding notions of solutions. Specifically, in studying hybrid seeking systems, we will parameterize solutions using a continuous-time index \( t \in \mathbb{R}_{\geq 0} \) and a discrete-time index \( j \in \mathbb{Z}_{\geq 0} \). The primary benefit of this parameterization is that it allows solutions to be defined in \emph{hybrid time domains}. For such solutions, the concept of \((T, \epsilon)\)-closeness can be naturally extended by measuring the distance between the \emph{graphs} of the nominal and perturbed solutions, that is, by using set-convergence (between graphs) as opposed to the standard uniform convergence. This can be achieved, provided the hybrid dynamics exhibit certain mild regularity properties. We will illustrate these regularity properties through various examples of hybrid seeking systems with different types of hybrid time domains. By defining solutions in such domains and employing graphical convergence to analyze the limiting behavior of perturbed solutions, we can seamlessly extend the perturbation-based methodology traditionally used in extremum-seeking control \cite{KrsticBookESC,DerivativesESC,TanAndNesic2006Local} to general hybrid set-seeking systems.
\subsection{Seeking Systems modeled as Hybrid Inclusions}
Most ES dynamics studied via averaging are modeled by ODEs of the form \eqref{oscillatingODE}. However, it is also useful to generalize these models for the study of hybrid seeking systems. For instance, hybrid systems often exhibit discontinuities, or "jumps," in the state, governed by a set-valued update rule of the form \( x^+ \in G(x) \), where \( x^+ \) denotes the state value immediately after the jump. This type of update rule, known as a "difference inclusion", allows \( x^+ \) to take any value within the set \( G(x) \), which may contain infinitely many points for a given \( x \). These models are powerful because they can capture more general and flexible behaviors and decision-making rules compared to standard difference equations of the form $x^+=g(x)$. Similarly, hybrid ES algorithms may also involve continuous-time dynamics described by differential inclusions of the form \( \dot{x} \in F(x) \), which can arise from discontinuous vector fields requiring appropriate regularizations (e.g., Krasovskii, Filippov, etc) for analysis, or from optimization and decision-making algorithms that are inherently set-valued \cite{goebel2024set}.

By incorporating set-valued updates in the continuous-time and/or discrete-time dynamics, hybrid inclusions can also provide a robust mathematical framework for modeling and analyzing seeking algorithms that are influenced by external discontinuous coordinating signals that govern system behavior. Such signals can arise in a variety of contexts, including switching dynamics \cite{Liberzon_Book}, time-varying optimization problems \cite{PovedaAutoOnline}, event-triggered systems \cite{RSanfeliceBook}, and systems with resetting timers \cite{HDS_haddad}, to name just a few examples. As such, most of the models in our study of hybrid set-seeking algorithms will be represented using hybrid inclusions.

The remainder of this paper is divided into three main sections. In Part 1, we present the preliminaries on hybrid systems and introduce the mathematical language necessary to study hybrid set-seeking systems. Part 2 focuses on hybrid seeking algorithms for plants characterized by static maps, while Part 3 extends these results to systems with dynamic plants. The main developments are complemented by several supporting discussions. "Continuous-Time Set-Valued Dynamical Systems" and "Discrete-Time Set-Valued Dynamical Systems" offer brief introductions to dynamical systems modeled as differential and difference inclusions, respectively. "Averaging Theory: From ODEs to Hybrid Inclusions" provides a mathematical overview of averaging theory applicable to hybrid set-seeking systems, and "Singularly Perturbed Hybrid Dynamical Systems" introduces the mathematical tools needed for studying hybrid systems in singular perturbation form. The notation used throughout this article is defined in Table \ref{tab1}.

\begin{summary}
\summaryinitial{T}his article focuses on the analysis and synthesis of extremum-seeking  systems that incorporate \emph{hybrid dynamics} in the loop, encompassing both continuous-time and discrete-time behaviors. Such extremum-seeking systems arise from the need to use advanced control and optimization algorithms to overcome the limitations of standard smooth feedback techniques and to meet stringent transient and steady-state demands in high-performance applications. They also emerge in settings where traditional controllers must be implemented on plants that naturally exhibit hybrid behaviors, such as in cyber-physical and autonomous systems, which rely on digital devices for computation, actuation, and sensing.

To explore hybrid extremum-seeking dynamics using control-theoretic tools, we begin by reviewing the key technical concepts that enable the development of perturbation theory for hybrid inclusions. These concepts provide the foundations for both averaging and singular perturbation theory in these systems. We then demonstrate how these tools can be applied to the analysis and synthesis of various hybrid extremum-seeking algorithms for plants modeled as static maps. In particular, we present several examples of set-valued and switching extremum-seeking algorithms under various switching scenarios: arbitrarily fast, under dwell-time and average dwell-time constraints, and under average activation time constraints. Additionally, we explore state-based switching extremum seeking in the context of obstacle avoidance problems, as well as Gradient-Newton switching extremum-seeking schemes. Other topics include extremum seeking with momentum and resets, with intermittent updates, slowly-varying parameters, and using safety mechanisms to incorporate constraints. In all these applications, we show how the perturbation-based approach, traditionally used for studying smooth extremum seeking control, can naturally be extended to hybrid systems, provided they satisfy suitable mild regularity conditions and their solutions are modeled in a hybrid time domain. Finally, in part 3 we address dynamic plants and discuss how the proposed hybrid extremum-seeking controllers can be integrated with such plants to solve model-free feedback optimization and decision-making problems. Building on established results for extremum seeking and using various illustrative examples, this article aims to provide an introductory and tutorial-style presentation of hybrid extremum-seeking systems, making it accessible to readers from different backgrounds, including those with no prior experience in hybrid systems.
\end{summary}

\section{Part 1: Preliminaries in Hybrid Systems}
%
%\subsection{Preliminary Concepts and Examples}
%
Hybrid dynamical systems (HDS), or simply ``hybrid systems'', are dynamical systems that combine continuous-time and discrete-time phenomena. Therefore, a natural mathematical representation of such systems would involve both differential equations and discrete-time recursions. In this article, we follow the formalism presented in \cite{bookHDS} to study general hybrid systems; see also \cite{CSM_hybrid} for a tutorial dedicated exclusively to such systems. In its simplest form, a hybrid system with state $x\in\mathbb{R}^n$, can be represented by the following expressions:
\begin{subequations}\label{eq:HDSmodel01}
\begin{align}
&x\in C,~~~~~~~\dot{x}= f(x),\label{eq:flow_map0}\\
&x\in D,~~~~~x^+=g(x)\label{eq:jump_map0},
\end{align}    
\end{subequations}
where \eqref{eq:flow_map0} describes the continuous-time dynamics, and \eqref{eq:jump_map0} describes the discrete-time dynamics. The notation in \eqref{eq:flow_map0} indicates that the state $x$ is allowed to evolve according to the continuous-time dynamics $\dot{x}=f(x)$ whenever the condition $x\in C$ holds, where the set $C\subset\mathbb{R}^n$ is called the \emph{flow set}. 
\begin{figure}[h!]
\begin{tcolorbox}[colback=iceblue!60, colframe=iceblue!60]
\centering
\includegraphics[width=0.95\linewidth]{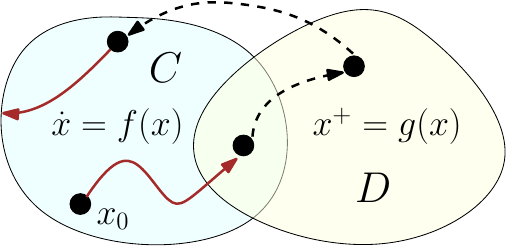}
\end{tcolorbox}
 \caption{Abstract representation of a solution of a hybrid system \eqref{eq:HDSmodel01} from the initial point $x_0$.\label{fig0}}
\end{figure} 
Similarly, the notation in \eqref{eq:jump_map0} indicates that $x$ is allowed to evolve according to the discrete-time dynamics $x^+=g(x)$ whenever the condition $x\in D$ holds, where the set $D\subset\mathbb{R}^n$ is called the \emph{jump set}, and where $x^+$ denotes the value of $x$ after an  ``instantaneous'' change --or jump-- in the state. For consistency, we call the functions $f:\mathbb{R}^n\to\mathbb{R}^n$ and $g:\mathbb{R}^n\to\mathbb{R}^n$ the \emph{flow map} and the \emph{jump map}, respectively. In this way, the hybrid system \eqref{eq:HDSmodel01} can be represented by the tuple $\mathcal{H}=\{f,C,g,D\}$, which we call the \emph{data} of the system. Note that when $D=\emptyset$, the system \eqref{eq:HDSmodel01} recovers a constrained continuous-time dynamical system. Similarly, when $C=\emptyset$, the system \eqref{eq:HDSmodel01} recovers a constrained discrete-time dynamical system. %In this way, system \eqref{eq:HDSmodel01} generalizes existing models of dynamical systems that are purely continuous and purely discrete.  
Figure \ref{fig0} shows an abstract representation of a possible "solution" to \eqref{eq:HDSmodel01}, starting from $x_0$, and experiencing an initial interval of flow, followed by two consecutive jumps and a final interval of flow that eventually hits the boundary of $C$, forcing the solution to stop. We will rigorously formalize this notion of "solution" and illustrate its application to study hybrid set-seeking systems.
\begin{table}
\caption{Notation used in the paper.\label{tab1}}
\begin{tabular*}{18.5pc}{|p{49pt}|p{148pt}|}
\hline
Symbol& ~ \\
\hline
$\mathbb{R}$ & Set of real numbers \\
$\mathbb{R}^n$ & The n-dimensional Euclidean space \\
$\mathbb{Z}_{\geq0}$&  Set of non-negative integers \\
$\dot{x}$ & Time derivative of the state of a system \\
$x^+$& State of a system after a jump  \\
$\emptyset$& The empty set \\
$\overline{\mathcal{A}}$& The closure of the set $\mathcal{A}$  \\
$\text{co}(\mathcal{A})$& The convex hull of the set $\mathcal{A}$ \\
$\text{bd}(A)$ & The boundary of the set $\mathcal{A}$\\
$x^\top$& The transpose of the vector $x$\\
$(x,y)$& Equivalent notation for  $[x^\top,y^\top]^\top$ \\
$|x|$& Euclidean norm of the vector $x$ \\
$\mathbb{B}$ & The closed unit ball, of appropriate dimension, in the Euclidean norm \\
$\mathbb{S}^1$&  The set $\{x=(x_1,x_2)\in\mathbb{R}^2:x_1^2+x_2^2=1\}$ \\
$\mathbb{S}^n$&  The nth-Cartesian product $\mathbb{S}^1\times\ldots\times\mathbb{S}^1$. \\
$f:\mathbb{R}^n\to\mathbb{R}^m$& A function from $\mathbb{R}^n$ to $\mathbb{R}^m$ \\
$F:\mathbb{R}^n\rightrightarrows\mathbb{R}^m$& A set-valued mapping from $\mathbb{R}^n$ to $\mathbb{R}^m$ \\
$\mathbf{R}_o$& The matrix $\mathbf{R}_o=\left[\begin{array}{cc}
0 & 1\\
-1 & 0
\end{array}\right]$ \\
$\mathcal{K}_{\infty}$ & Class of functions $\alpha:\mathbb{R}_{\geq0}\to\mathbb{R}_{\geq0}$ that are zero at zero, continuous, strictly increasing, and unbounded.\\ 
$\mathcal{KL}$ & Class of functions $\beta:\mathbb{R}_{\geq0}\times\mathbb{R}_{\geq0}\to\mathbb{R}_{\geq0}$ that are nondecreasing in their first argument, nonincreasing in their second argument, satisfy $\lim_{r\to0^+}\beta(r,s)=0$ for each $s\in\mathbb{R}_{\geq0}$, and $\lim_{s\to\infty}\beta(r,s)=0$ for each $r\in\mathbb{R}_{\geq0}$. \\
\hline
\end{tabular*}
%\label{tab1}
\end{table}
%
%
%\vspace{-1cm}~
%
\begin{example}[Systems with Periodic Resets]\label{example1}
To illustrate how the data in \eqref{eq:HDSmodel01} can be used to model a common class of systems encountered in practice in control, such as systems with periodic resets, let $x=(x_1,x_2)\in\mathbb{R}^2$, and consider the following hybrid dynamics:
\begin{align*}
&x\in C:=[0,T]\times\mathbb{R},~~~\dot{x}=\left(\begin{array}{c}
\dot{x}_1\\
\dot{x}_2
\end{array}\right)=f(x):=\left(\begin{array}{c}
1\\
-\gamma x_2
\end{array}\right),\\
&x\in D:=\{T\}\times\mathbb{R},~~x^+=\left(\begin{array}{c}
x_1^+\\
x_2^+
\end{array}\right)=g(x):=\left(\begin{array}{c}
0\\
\frac{1}{2}x_2
\end{array}\right),
\end{align*}
where $T,\gamma>0$ are tunable parameters.  In this system, the state $x_1$ can be seen as a timer that grows linearly during flows and is reset to zero every time the condition $x_1=T$ is satisfied. Whenever such resets occur, the state $x_2$ is also reset to one half of its previous value. Since resets of $x_1$ force $x$ to satisfy conditions $x\in C$ and $x\notin D$, every jump in the system is followed by a period of flow of duration $T$, during which the continuous-time
\begin{figure}[h!]
 \begin{tcolorbox}[colback=iceblue, colframe=iceblue]
\centering
\includegraphics[width=\linewidth]{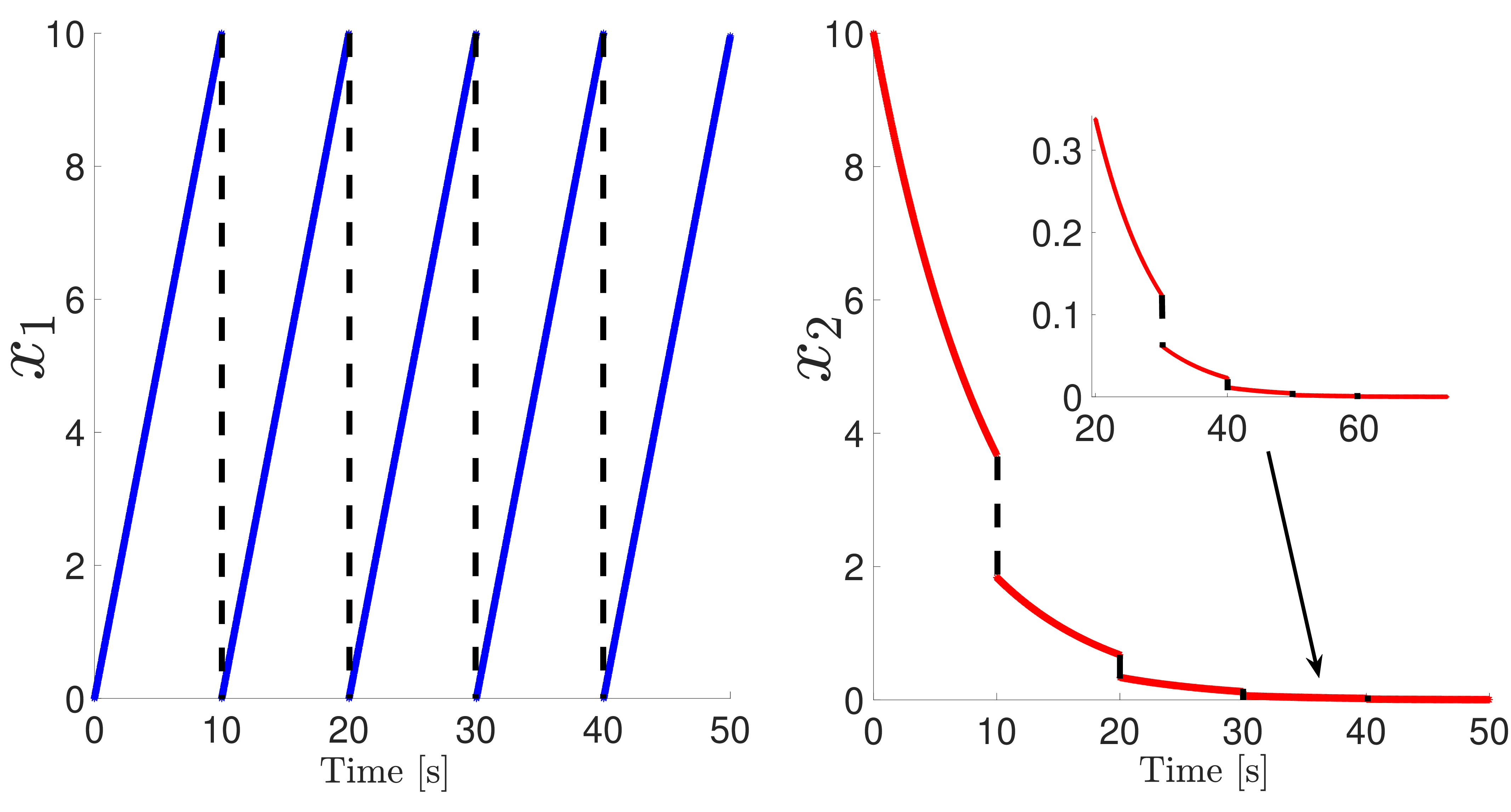}
\end{tcolorbox}
 \caption{Trajectory generated by the hybrid system of Example 1, with initial condition $x_0=(0,10)$. Left: Trajectory of $x_1$. Right: Trajectory of $x_2$.\label{fig1}}
\end{figure}
 \begin{figure*}[t!]
\begin{tcolorbox}[colback=iceblue, colframe=iceblue]
 \centering
\includegraphics[width=0.49\linewidth]{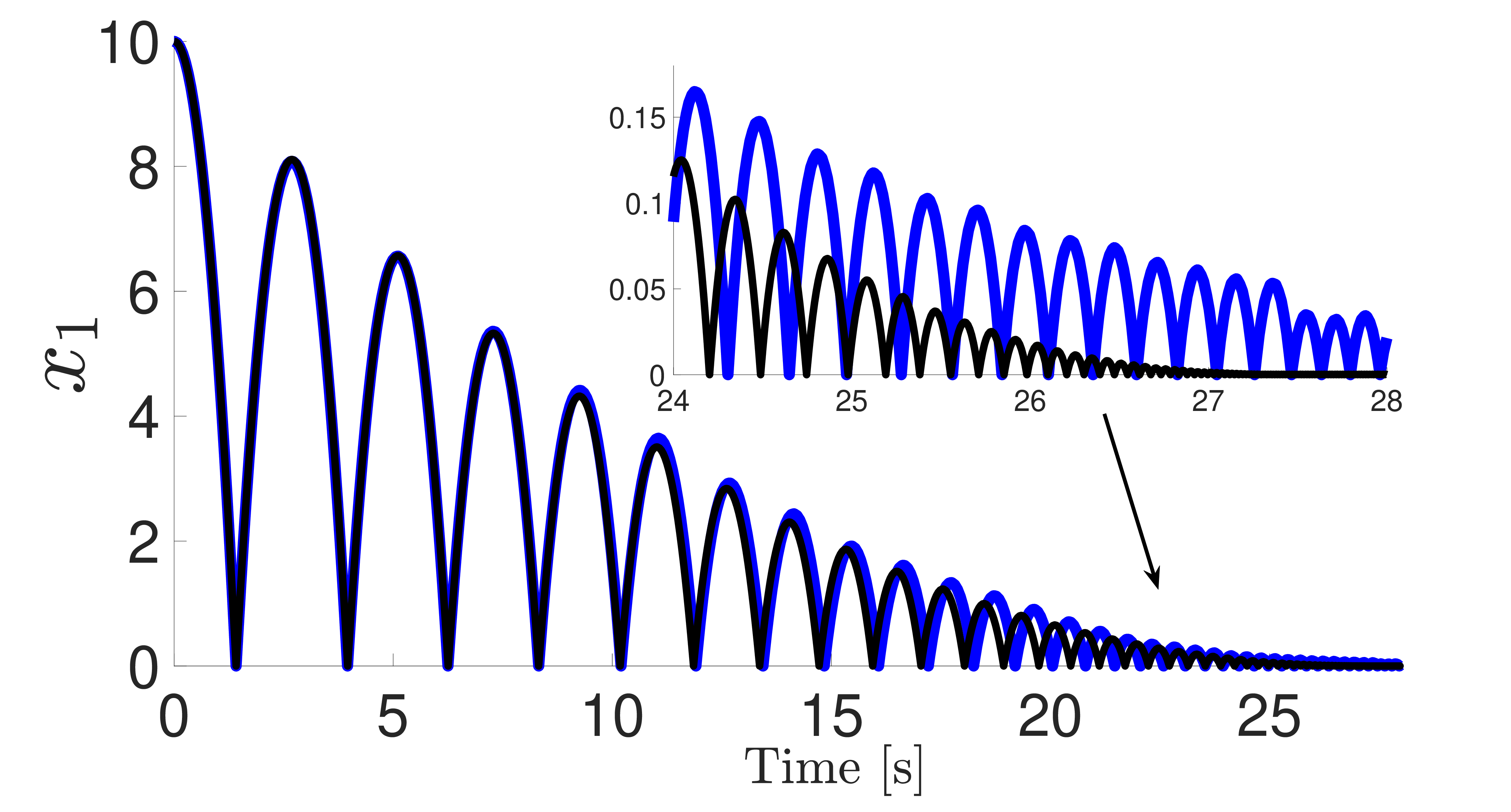}\includegraphics[width=0.49\linewidth]{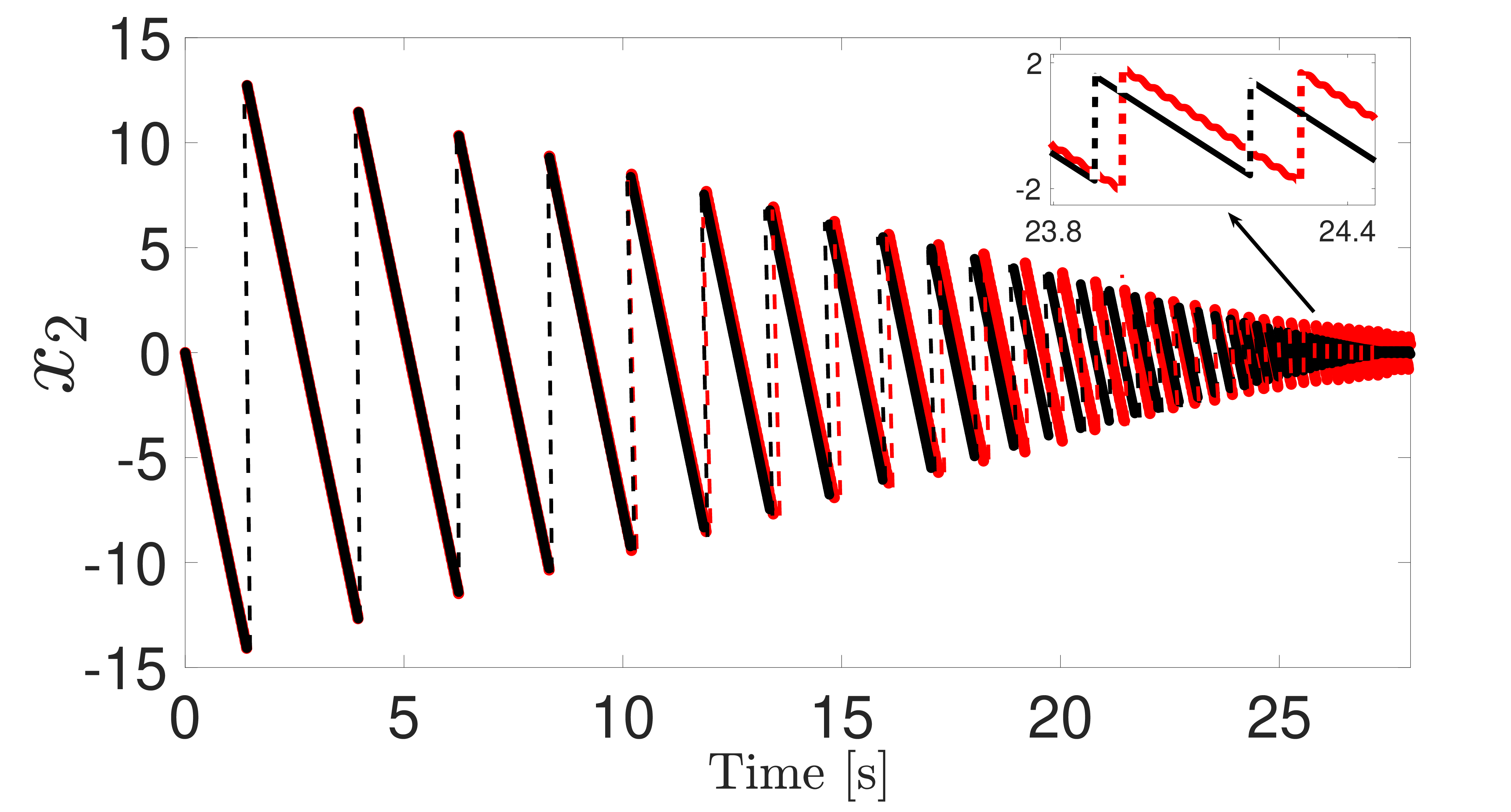}
\end{tcolorbox}
 \caption{Trajectories of the "bouncing" seeking system with time-varying function $\gamma(\tau)$ using $\omega=100$ (shown in color blue), and with constant  $\gamma=\bar{\gamma}=10$ (shown in color black). Left: Trajectories of the position $x_1$. Right: Trajectories of the velocity $x_2$, which experience jumps, illustrated with dashed lines.\label{fig2}}
\end{figure*}
 dynamics of $x_2$ make it decrease exponentially at a rate dictated by $\gamma$. Figure \ref{fig1} shows the trajectories of $x_1$ and $x_2$ generated by this simple hybrid system, using $T=10$, $\gamma=\frac{1}{10}$, and from the initial conditions $(0,10)$. The dashed lines indicate the jumps in the values of the states $x$ triggered by the discrete-time dynamics $x^+=g(x)$. Note that for this system, the intersection of the flow set and the jump set is not empty. In fact, we have $D\cap C=\{T\}$. However, for initial conditions satisfying $x_1=T$, the continuous-time dynamics $\dot{x}_1=1$ cannot flow and extend $x_1$ without leaving the set $C$, for example, no solutions can evolve via flows from this initialization. On the other hand, from this same initial condition, the state $x_1$ can be extended in time via jumps by immediately resetting $x_1$  and $x_2$. This type of behavior will be formalized later when we discuss the notion of solutions to hybrid systems of the form \eqref{eq:HDSmodel01}. In the meantime, we note that, while the state $x_2$ converges to zero in Figure \ref{fig1}, the state $x_1$ continues to oscillate restricted to the set $[0,T]$. Therefore, for this hybrid system it is reasonable to study the stability properties with respect to the compact \emph{set} $\mathcal{A}=[0,T]\times \{0\}$.  This simple example serves to illustrate the fact that, unlike differential equations of the form \eqref{ESCODEVanilla}, stability properties in hybrid seeking systems will be studied with respect to general sets, as opposed to equilibrium points. %Naturally, this feature will also apply to hybrid versions of the ES algorithm \eqref{ESCODEVanilla}.
 \QEDB 
\end{example}
Dynamic timers, such as those in Example \ref{example1}, are common in hybrid systems to model time-triggered jumps, which can arise in sampled-data systems, switching algorithms, event-triggered control, and reset-based controllers, to name just a few examples. Incorporating exogenous timers into hybrid systems is also useful for studying time-varying dynamics, which are particularly relevant in ES systems, such as \eqref{ESCODEVanilla}. In this case, the overall dynamics can be written as a hybrid system with  states \((x,\tau) \in \mathbb{R}^n \times \mathbb{R}\), and the dynamics given by
\begin{subequations}\label{eq:HDSmodel0}
\begin{align}
&(x,\tau)\in C\times\mathbb{R}_{\geq0},~~~~\dot{x}= f(x,\tau),~~~\dot{\tau}=\omega,\label{eq:flow_map001}\\
&(x,\tau)\in D\times\mathbb{R}_{\geq0},~~x^+=g(x),~~~~\tau^+=\tau\label{eq:jump_map001},
\end{align}    
\end{subequations}
where $\omega>0$ dictates the rate of evolution of $\tau$. Note that \eqref{eq:HDSmodel0} can be written as \eqref{eq:HDSmodel01} using the extended state $y=(x,\tau)$.

The following example illustrates a class of time-varying hybrid systems with oscillating flow maps that are of interest in this article, which also exhibit a well-known behavior that can emerge in hybrid systems: Zeno phenomena. 

\begin{sidebar}{Continuous-Time Set-Valued Dynamical Systems}
\sdbarinitial{S}et-valued continuous-time dynamics modeled as \emph{differential inclusions} typically arise in hybrid seeking systems due to discontinuous feedback-based optimization or decision-making algorithms, plants with uncertain parameters, or auxiliary states used to generate families of switching signals. Here, we review some basic information on differential inclusions \cite[Ch.4]{bookHDS}, \cite{TeelPraly,SAA3}.

\setcounter{sequation}{0}
\renewcommand{\thesequation}{S\arabic{sequation}}
\setcounter{stable}{0}
\renewcommand{\thestable}{S\arabic{stable}}
\setcounter{sfigure}{0}
\renewcommand{\thesfigure}{S\arabic{sfigure}}

\section{Differential Inclusions}
Consider a constrained differential inclusion of the form
\begin{equation}\label{DI0}
x\in C,~~~\dot{x}\in F(x),
\end{equation}
where $C\subset\mathbb{R}^n$ and $F:\mathbb{R}^n\rightrightarrows\mathbb{R}^n$. System \eqref{DI0} is said to satisfy the \emph{Basic Conditions} if $C$ is a closed set, and $F$ is outer-semicontinuous and locally bounded relative to $C$, with $F(x)$ being non-empty and convex for each $x\in C$. Given $T>0$, an absolutely continuous function $x:[0,T]\to\mathbb{R}^n$ is said to be a solution of \eqref{DI0} if: a) $x(0)\in C$, and b) for almost all $t\in[0,T]$:
\begin{equation}
    \frac{dx(t)}{dt}\in F(x(t))~~\text{and}~~x(t)\in C.
\end{equation}
When \eqref{DI0} satisfies the Basic Conditions, the existence of solutions can be studied under the so-called viability conditions \cite[Lemma 5.26]{bookHDS}. To state such conditions, recall that the \emph{tangent cone} to a set $S\subset\mathbb{R}^n$ at a point $x\in\mathbb{R}^n$, denoted $T_S(x)$, is given by all vectors $s\in\mathbb{R}^n$ for which there exist $x_i\in S$, $\tau_i>0$ with $x_i\to x$ and $\tau_i\to0^+$, such that $s=\lim_{i\to\infty}\frac{x_i-x}{\tau_i}$. In particular, under the Basic Conditions, we have that: 
\begin{enumerate}[(a)]
\item (Necessity) If $x:[0,T]\to\mathbb{R}^n$, for some $T>0$, is a solution of \eqref{DI0}, then
\begin{equation}
F(x(0))\cap T_C(x(0))\neq\emptyset.
\end{equation}
\item (Sufficiency) Given $x_0\in C$, if there exists a neighborhood $U$ of $x_0$ such that for all $x\in U\cap C$,
\begin{equation}
F(x)\cap T_C(x)\neq\emptyset,
\end{equation}
then there exists $T>0$ and a solution $x:[0,T]\to\mathbb{R}^n$ to \eqref{DI0} with $x(0)=x_0$.
\end{enumerate}
 Note that, in general, given an initial condition $x(0)=x_0\in\mathbb{R}^n$, solutions to \eqref{DI0} from $x_0$ might not be unique.

\section{Krasovskii Solutions of ODEs}
Differential inclusions \eqref{DI0} often arise when studying discontinuous differential equations for which solutions might not exist from some initial conditions. In particular, given a set $C\subset\mathbb{R}^n$, a function $f:\mathbb{R}^n\to\mathbb{R}^n$, and an absolutely continuous function $x:[0,T]\to\mathbb{R}^n$, the function $x$ is said to be a Krasovskii solution of the constrained system $x\in C$, $\dot{x}=f(x)$ if 
\begin{equation}\label{krasovskiicons}
x(t)\in\overline{C},~\text{and}~\dot{x}(t)\in F_K(x(t)):=\bigcap _{\epsilon>0}\overline{\text{co}}~f((x(t)+\epsilon\mathbb{B})\cap C),
\end{equation}
for almost all $t\in[0,T]$. When $f$ is continuous, the mappings $F_K$ and $f$ coincide. The construction \eqref{krasovskiicons} is also valid whenever $f$ is set-valued but does not necessarily satisfy the Basic Conditions. Indeed, by construction (see \cite[Example 6.6]{bookHDS}), if $f$ is locally bounded relative to $\overline{C}$, then the differential inclusion \eqref{krasovskiicons} satisfies the Basic Conditions.
\section{Robustness: Continuous-Time Hermes Solutions}
For differential equations with a continuous right-hand side, small additive perturbations to the states do not significantly affect the system's behavior. In contrast, small perturbations can drastically impact the behavior of systems with a discontinuous right-hand side. The effect of arbitrarily small perturbations on such systems can be analyzed using the concept of Hermes solutions \cite{SAA4}. The function $x:[0,T]\to\mathbb{R}^n$ is said to be a Hermes solution of the constrained system $x\in C$, $\dot{x}=f(x)$ if there exists a sequence of absolutely continuous functions $x_i:[0,T]\to\mathbb{R}^n$ and a sequence of measurable functions $e_i:[0,T]\to\mathbb{R}^n$ such that
\begin{align*}
&x_i(t)+e_i(t)\in C,~~\text{for all}~t\in(0,T),\\
&\dot{x}_i(t)=f(x_i(t)+e_i(t)),~~~\text{for almost all}~t\in[0,T],
\end{align*}
the sequence $\{x_i\}_{i=1}^{\infty}$ converges uniformly to $x$ on $[0,T]$, and the sequence $\{e_i\}_{i=1}^{\infty}$ converges uniformly to the zero function on $[0,T]$. In words, Hermes solutions capture the limiting effect of arbitrarily small additive state perturbations acting on the system $x\in C,~\dot{x}=f(x)$. The connection between Hermes solutions and Krasovskii solutions is established in the following lemma, corresponding to \cite[Thm. 4.3]{bookHDS}.

\vspace{-0.1cm}
\begin{lemma}[Continuous-Time Krasovskii Solutions]\label{hermesandkrasovskii}
Suppose that $f$ is locally bounded and let $x:[0,T]\to\mathbb{R}^n$ be an absolutely continuous function. Then $x$ is a Hermes solution of $\dot{x}=f(x)$ if and only if it is a Krasovskii solution of the same system.
\end{lemma}

\vspace{-0.2cm}\noindent 
Put simply, Lemma \ref{hermesandkrasovskii} suggests that to accurately characterize the behavior of the solutions of the system $\dot{x} = f(x)$ under arbitrarily small additive perturbations of the state, we should consider its Krasovskii solutions. A different concept of solutions explored in the literature of discontinuous differential equations is that of Filippov solutions \cite{SAA3}. An absolutely continuous function $x:[0,T]\to\mathbb{R}^n$ is said to be a Filippov solution of $\dot{x}=f(x)$ if it satisfies 
\begin{equation}\label{filippovcons}
\dot{x}(t)\in F_F(x(t)):=\bigcap _{\epsilon>0}\bigcap _{\mu(S)=0}\overline{\text{co}}~f((x(t)+\epsilon\mathbb{B})\backslash S),
\end{equation}
where $\mu$ is the Lebesgue measure. While the constructions \eqref{krasovskiicons} and \eqref{filippovcons} are similar, the construction in \eqref{filippovcons} ignores the behavior of $f$ on sets $S$ with zero measure. Hence, Filippov solutions are Krasovskii solutions, but the converse is generally not true. Thus, in light of Lemma \ref{hermesandkrasovskii}, Filippov solutions are problematic as they do not completely capture the behavior of the original dynamics under small additive state perturbations. For further discussions on the difference between Filippov and Krasovskii solutions and for a proof of Lemma 1 for the special case where $C=\mathbb{R}^{n}$, we refer the reader to \cite{SAA5}.

% \section{Set-valued seeking dynamics}

% Here is some random text showing how references are handled for a sidebar. Reference \cite{SAA1} gives one random reference. A basic equation is
% \begin{sequation}
% f = ma,\label{seq1}
% \end{sequation}
% and the impact of \eqref{seq1} is highlighted in Figure \ref{sfig1}. 

%Reference \cite{SAA2} gives another random reference, but \cite{AA1} and \cite{BB1} are references in the main body of the text that can be cited in the sidebar.

\end{sidebar}

\begin{example}["Bouncing" Seeking Systems]\label{seekingballexample1}
Bouncing balls are classic examples of hybrid systems. Their dynamics are governed by Newtonian physics while the ball is in the air and experience instantaneous jumps upon hitting the ground. These models are most accurate when the ball is rigid and the deformations caused by impacts are negligible. Inspired by classic bouncing ball systems \cite[Ex. 1.1]{bookHDS}, we can consider a hybrid system of the form \eqref{eq:HDSmodel0}, with state $x=(x_1,x_2)\in\mathbb{R}^2$, where $x_1$ models the vertical position of a point-mass ball bouncing vertically on a horizontal surface and $x_2$ models the vertical velocity of the ball. Following the notation in \eqref{eq:HDSmodel0}, the data of this system is given by 
\begin{subequations}\label{seekingballsystems}
\begin{align}
C&:=\{(x_1,x_2)\in\mathbb{R}^2: x_1\geq0\},\label{flowsetbouncing}\\
\left(\begin{array}{c}
\dot{x}_1\\
\dot{x}_2\\
\end{array}\right)&=
f(x,\tau):=\left(\begin{array}{c}
x_2\\
-\gamma(\tau),
\end{array}\right),~~\dot{\tau}=\omega,\label{flowsseekingball}\\
D&=\{(x_1,x_2)\in\mathbb{R}^2: x_1=0,~x_2\leq0\},\label{jumpsetball}\\
\left(\begin{array}{c}
x_1^+\\
x_2^+
\end{array}\right)&=g(x):=\left(\begin{array}{c}
0\\
- \lambda x_2
\end{array}\right),~~~\tau^+=\tau,\label{jumpmapball}
\end{align}
\end{subequations}
where the time-varying function $\gamma(\cdot)$ is given by
\begin{equation}\label{oscillatorygravity}
\gamma(\tau)=2\overline{\gamma} \sin(\tau)^2,
\end{equation}
and where \( \overline{\gamma} \) denotes the nominal gravitational acceleration constant and $\lambda\in(0,1)$ is a restitution coefficient. According to \eqref{flowsetbouncing}, whenever the altitude is non-negative, the ball flows through the air via \eqref{flowsseekingball} under oscillatory acceleration $-\gamma(\tau)$. On the other hand, according to \eqref{jumpsetball}-\eqref{jumpmapball}, when the ball touches the ground (i.e., $x_1=0$) and its velocity is non-positive, the mapping $g$ resets the sign of the velocity and also decreases its magnitude using the restitution coefficient $\lambda$. Note that the position of the ball does not change during jumps. Figure \ref{fig2} illustrates the behavior of the solutions generated by this ``bouncing'' seeking system when $\omega=100$ (colored trajectories), and also when $\gamma(\tau)=\overline{\gamma}$ is constant (black trajectories). \tcb{It can be observed that the highly oscillatory hybrid dynamics that use $\omega=100$ provide a suitable approximation of the behavior of the hybrid system that uses a constant acceleration $\gamma$.} This approximation is particularly accurate at the beginning of the simulation but, as shown on the inset, it eventually deteriorates, as seen by the bouncing times starting to deviate, as the number of jumps grows unbounded. Note that the hybrid system with constant acceleration $\dot{x}_2=-\overline{\gamma}$ exhibits Zeno behavior as time approaches approximately 28s, as the amount of time between jumps decreases to zero. Similar phenomena typically emerge in certain hybrid control systems such as in event-triggered controllers \cite{event1}. This behavior is also "approximately" inherited by the hybrid system \eqref{seekingballsystems} with oscillatory acceleration $\gamma(\tau)$. To study hybrid set-seeking systems, we will formalize this "approximation" property using averaging theory for hybrid systems. \QEDB 
\end{example}

Switching systems constitute a class of hybrid dynamical systems that arise frequently in applications. When such systems are modeled within a hybrid systems framework, their qualitative behavior is typically analyzed with respect to prescribed \emph{families} of switching signals \cite{Liberzon_Book,Average_Dwell_time}. In this case, set-valued dynamics can provide a convenient mathematical tool to model different types of switching signals. We illustrate this idea in the following example.
\begin{pullquote}
The primary aim of this paper is to provide an accessible introduction to non-smooth and hybrid set-seeking systems formulated as hybrid inclusions
\end{pullquote}

\begin{example}[Source-Seeking and Surveillance]\label{intermittentsourceseeking}
Consider a point-mass vehicle in the plane, with state $x\in\mathbb{R}^2$ modeling its position, and simple dynamics given by a single integrator
\begin{equation}\label{pointmassopen}
\dot{x}=\alpha,
\end{equation}
where $\alpha$ is the input to be designed. We study a typical application of ES in the context of source-seeking \cite{SourceSeekingWater,ESCsource,Poveda:20TAC}, where the vehicle \eqref{pointmassopen} seeks to stabilize its position at the location where a potential signal, only available via measurements, attains its maximum intensity. However, unlike standard source-seeking problems with a single static potential field, we consider the presence of multiple intermittent potential fields, which are "active" only during bounded periods of time and never simultaneously. In this case, the goal of the vehicle is to achieve \emph{persistent source-seeking and surveillance} assuming the intermittence of the sources is low, i.e., they remain individually active for "sufficiently" long periods of time. This problem can be modeled using the formalism of hybrid set-seeking systems. In particular, let $\mathcal{Q}:=\{0,1,\ldots,q_{\max}\}$ denote a finite set of logic modes,  where $q_{\max}\in\mathbb{Z}_{\geq1}$. For each $q\in\mathcal{Q}$, we let $J_q:\mathbb{R}^2\to\mathbb{R}$ denote a potential field that attains its maximum at the point $x_q^*$. For simplicity, we assume that each function $-J_q(\cdot)$ is strongly convex and has a globally Lipschitz gradient $\nabla J_q$. For each fixed potential field $J_q(\cdot)$, the vehicle can achieve source-seeking by implementing the following standard seeking feedback law \cite{ZhangAutomatica}:
\begin{figure}[t!]
\begin{tcolorbox}[colback=ivoryA, colframe=ivoryA]
\centering
\includegraphics[width=1.05\linewidth]{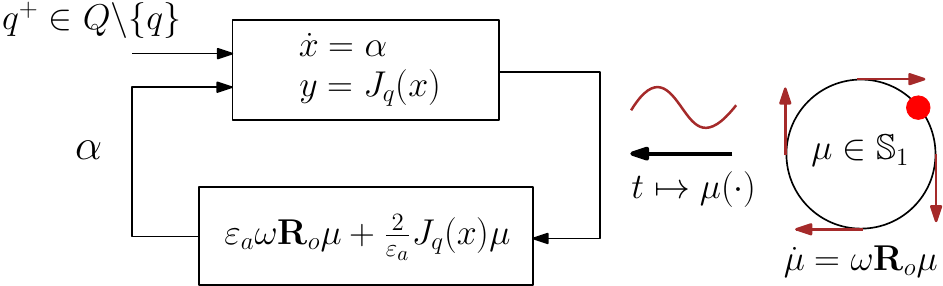}
\end{tcolorbox}
 \caption{Scheme of source-seeking dynamics for point-mass vehicle models under intermittent switching potentials with distinct optimal points $x_q^*$.\label{fig3cc}}
\end{figure}
\begin{equation}\label{feedbacklaw1234}
\alpha=\varepsilon_a\omega \mathbf{R}_o\mu+\frac{2}{\varepsilon_a}J_q(x)\mu,
\end{equation}
where $\varepsilon_a,\omega>0$ are tunable parameters and $\mu$ is a state that models a dithering signal with dynamics given by the following constrained differential equation:
\begin{equation}\label{dynamicoscillator}
\mu\in\mathbb{S}^1,~~~~\dot{\mu}=\omega \mathbf{R}_o\mu.
\end{equation}
The linear oscillator \eqref{dynamicoscillator} provides an alternative time-invariant approach to generate sinusoidal functions, such as those used in \eqref{ESCODEVanilla}, for the purpose of exploration. In particular, by direct computation, any solution of \eqref{dynamicoscillator} has the form
\begin{equation}
\mu(t)=\left[\begin{array}{cc}
\cos(\omega t) & \sin(\omega t)\\
-\sin(\omega t) & \cos(\omega t)
\end{array}\right]\left[\begin{array}{c}
\mu_1(0)\\
\mu_2(0)
\end{array}\right],
\end{equation}
where $\mu_1(0)^2+\mu_2(0)^2=1$. In this way, the feedback law \eqref{feedbacklaw1234} emulates the source-seeking algorithm studied in \cite{ZhangAutomatica} in the context of single-source seeking problems, but using a time-invariant model of the dynamics (with $\mu$ restricted to evolve in $\mathbb{S}^1$), see Fig. \ref{fig3cc} for a graphical representation of the system.

To incorporate the switching behavior of the logic state $q$ that indexes the potential fields $J_q$, we can consider the following set-valued hybrid dynamics with states $(q,\tau)\in\mathcal{Q}\times\mathbb{R}_{\geq0}$, which are parameterized by the constants $N_0\geq1$ and $\eta_d>0$:
\begin{subequations}\label{hybridtimer1}
\begin{align}
&(q,\tau)\in \mathcal{Q}\times[0,N_0],~~~~~\dot{q}=0,~~~~~~~~~~\dot{\tau}\in[0,\eta_d],\label{setvaluedtimer01}\\
&(q,\tau)\in \mathcal{Q}\times[1,N_0],~~~q^+\in \mathcal{Q}\backslash\{q\},~~~\tau^+=\tau-1.
\end{align}
\end{subequations}
Similar to Example \ref{example1}, in this system the state $\tau$ acts as a timer that regulates how frequently $q$ can jump to other values in the set $\mathcal{Q}$. However, in this case, the continuous-time dynamics of $\tau$ in \eqref{setvaluedtimer01} are set-valued and allow the derivative $\dot{\tau}$ to be any (measurable) selection from the set $[0,\eta_d]$, which includes keeping $\tau$ constant via $\dot{\tau}=0$, increasing $\tau$ at the maximum rate $\eta_d>0$ until $\tau=N_0$, or any other absolutely continuous function $\tau$ (not necessarily constant) having a derivative lower bounded by $0$ and upper bounded $\eta_d$. This flexibility allows the model \eqref{hybridtimer1} to capture a large class of switching signals $q$ acting on the feedback law \eqref{feedbacklaw1234}, including signals that never switch (corresponding to $\dot{\tau}=0$ for all time). In fact, for the case where $N_0=1$, \emph{all} switching signals generated by \eqref{hybridtimer1} satisfy the \emph{dwell-time} condition \cite{Dwell_time}. When $N_0\in\mathbb{Z}_{\geq2}$, \emph{all} switching signals generated by \eqref{hybridtimer1} satisfy an \emph{average dwell-time} condition \cite{Average_Dwell_time} (see also "Switching Set-Seeking Systems with Average
Dwell-time"). Both classes of switching signals are common in the modeling and stability analysis of switching systems, and by using \eqref{hybridtimer1} interconnected with \eqref{pointmassopen}- \eqref{dynamicoscillator} (with $x^+=x$), we obtain a hybrid system suitable for the study of the source-seeking feedback law \eqref{feedbacklaw1234} under average dwell-time switches of $J_q$. Figure \ref{fig3c} shows an example of this switching signal for the case when $\eta_d=1/100$, $N_0=1$, and $q_{\max}=1$, such that $\mathcal{Q}=\{0,1\}$. When $\eta_d$ is sufficiently small, the resulting slowly-switching set-seeking system generates the stable behavior shown in Figure \ref{fig3}. In the left plot, we show the trajectories of \eqref{pointmassopen} under the feedback law \eqref{feedbacklaw1234} with $\omega=1000$, $\epsilon_a=0.01$, and from different initial conditions shown with red dots. The trajectories converge to a neighborhood of the set $\Omega(K)$, shown in green. 

\begin{figure}[b]
\begin{tcolorbox}[colback=iceblue, colframe=iceblue]
\includegraphics[width=\linewidth]{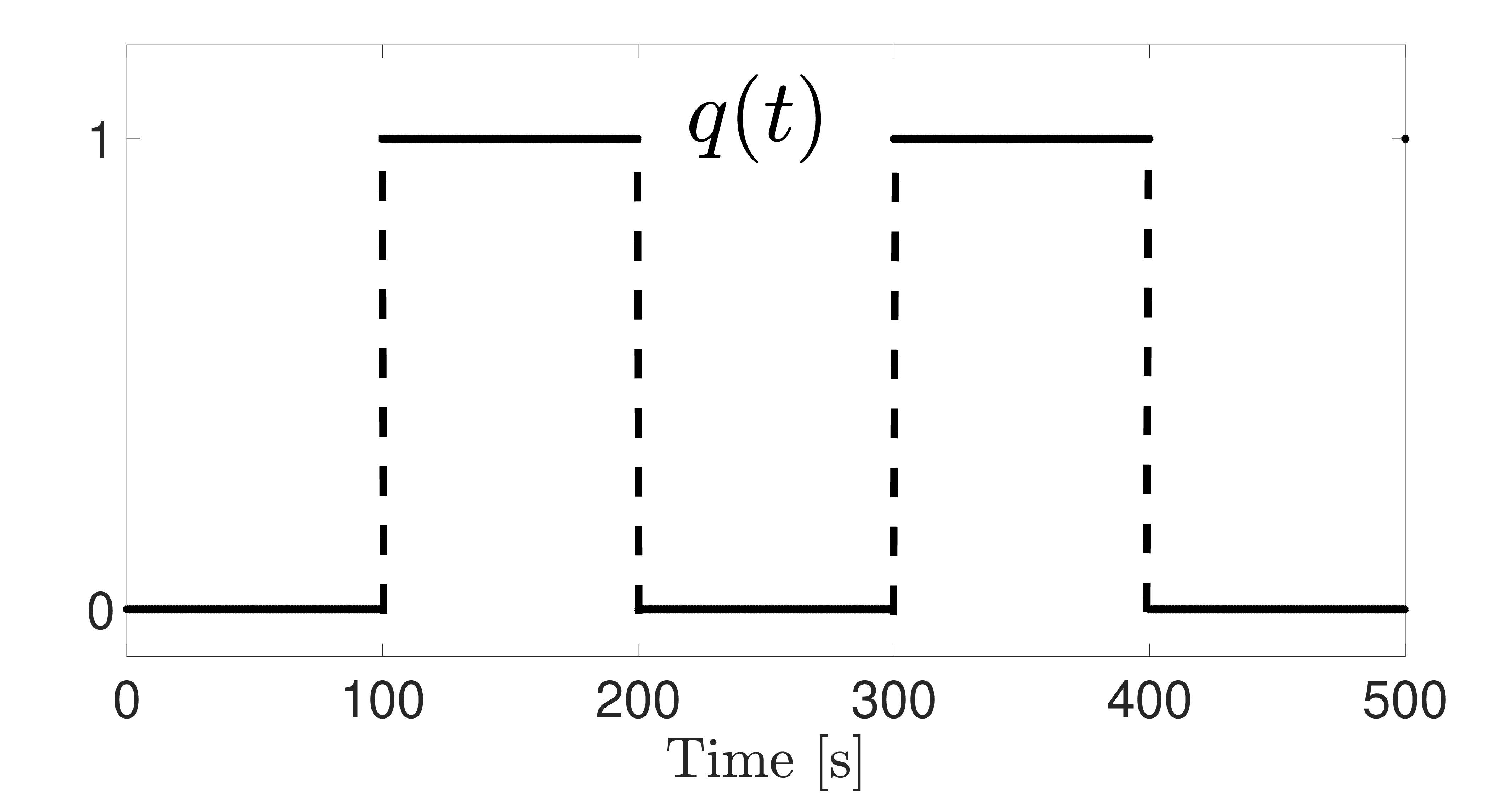}
\end{tcolorbox}
 \caption{Evolution in time of a switching signal $q$ taking values in the set $\{0,1\}$, indicating which potential field is active in the source-seeking problem.\label{fig3c}}
\end{figure}

\begin{figure*}
\begin{tcolorbox}[colback=iceblue, colframe=iceblue]
\includegraphics[width=0.49\linewidth]{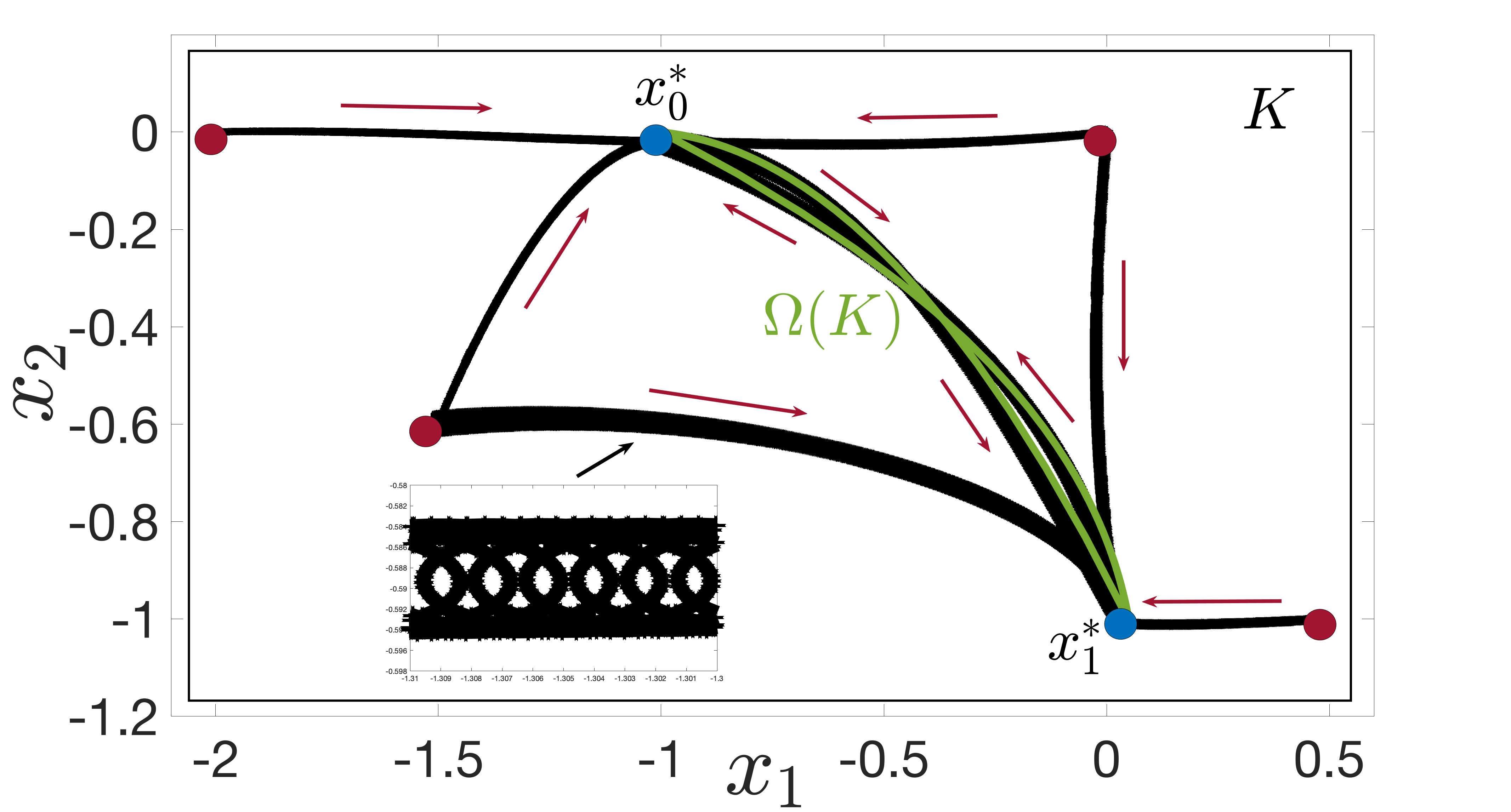}\includegraphics[width=0.49\linewidth]{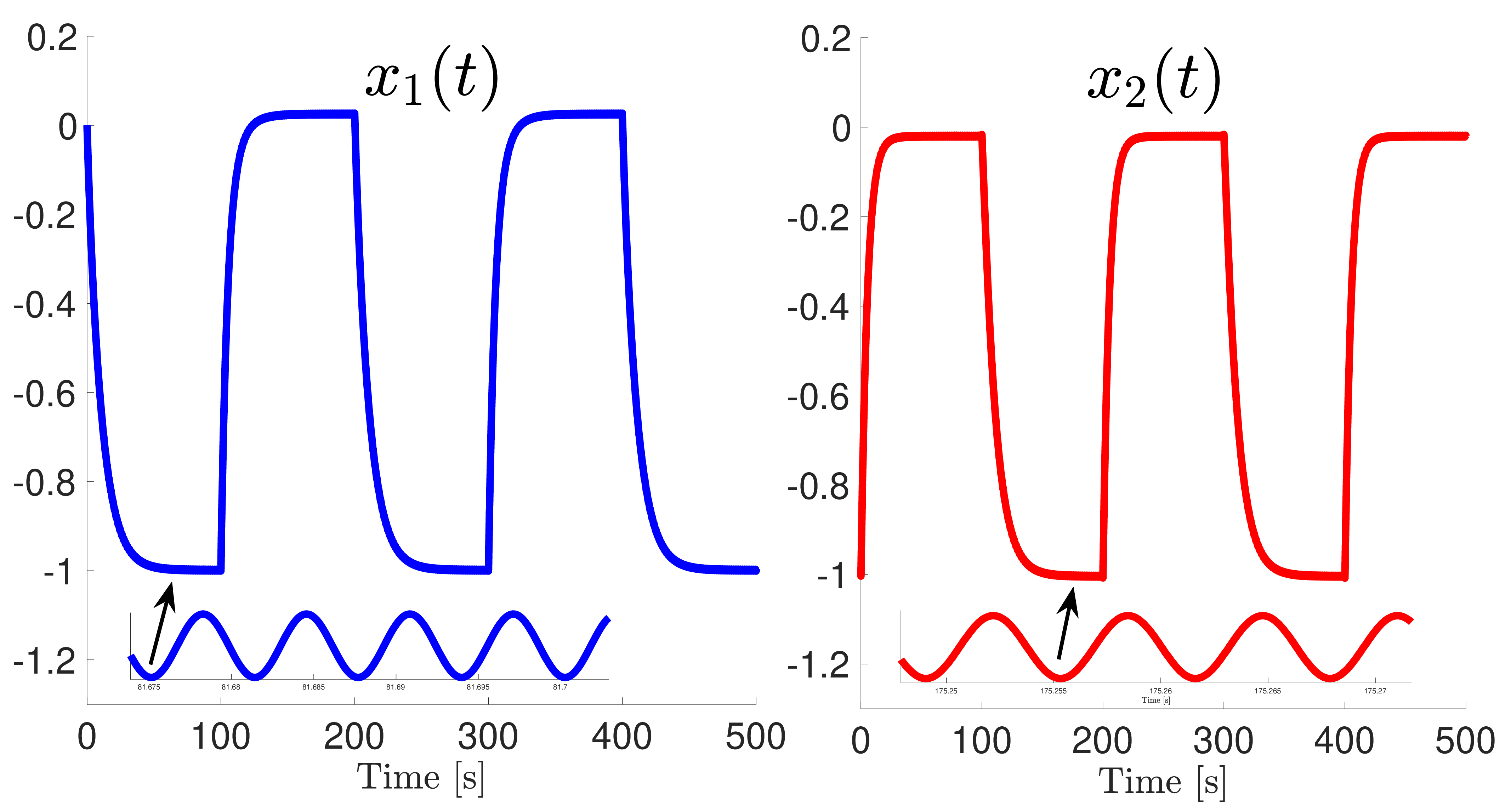}
\end{tcolorbox}
 \caption{Left: Trajectories of system \eqref{pointmassopen} (shown in color black) under the feedback law \eqref{feedbacklaw1234}, evolving on the plane, and converging to a neighborhood of the set $\Omega(K)$, shown in color green. Right: Evolution in time of the position of the vehicle, which oscillates between the two optimal points $x_0^*$ and $x_1^*$, which are contained in $\Omega(K)$. \label{fig3}}
\end{figure*}
This set, called the "Omega-limit set from $K$", described later in equation \eqref{omegalimitset}, turns out to be (semi-globally practically) asymptotically stable for the overall hybrid set-seeking system \eqref{pointmassopen}-\eqref{hybridtimer1} as $\eta_d\to0^+$, $\varepsilon_a\to0^+$, and $1/\omega\to0^+$. While in certain cases the set \( \Omega(K) \) can be computed explicitly \cite{baradaran2020omega}, its mere existence is typically guaranteed in slowly switching systems with individually asymptotically stable modes, under mild assumptions
 \cite[Sec. 7.4]{bookHDS}.\QEDB 
\end{example}
The preceding example underscores the advantages of considering hybrid set-seeking systems with \emph{set-valued} dynamics. Consequently, for the remainder of this paper, we will focus on hybrid systems that generalize \eqref{eq:HDSmodel01}, and are represented by the following \emph{hybrid inclusions}:
\begin{subequations}\label{eq:HDSmodel}
\begin{align}
&x\in C,~~~~~~~\dot{x}\in F(x),\label{eq:flow_map}\\
&x\in D,~~~~x^+\in G(x)\label{eq:jump_map},
\end{align}    
\end{subequations}
where the flow map $F:\mathbb{R}^n\rightrightarrows\mathbb{R}^n$ and the jump map $G:\mathbb{R}^n\rightrightarrows\mathbb{R}^n$ are in general set valued. Note that for particular cases in which such maps are single-valued, we can still use the formalism \eqref{eq:HDSmodel} by defining 
\begin{subequations}\label{setvaluedhds}
\begin{align}
F(x)&:=\left\{\begin{array}{ll}
\{f(x)\},&x\in C\\
\emptyset,&x\notin C
\end{array}\right.\\
G(x)&:=\left\{\begin{array}{ll}
\{g(x)\},&x\in D\\
\emptyset,&x\notin D.
\end{array}\right. 
\end{align}
\end{subequations}
The data of \eqref{eq:HDSmodel} is also denoted $\mathcal{H}:=\{C, F, D, G\}$, and, similarly to \eqref{eq:HDSmodel0}, for the case where \eqref{eq:HDSmodel} depends on exogenous time-varying signals, we can consider models of the form
\begin{subequations}\label{eq:HDSmodel2}
\begin{align}
&(x,\tau)\in C\times\mathbb{R}_{\geq0},~~~~~~~\dot{x}\in F(x,\tau),~~~~\dot{\tau}=\omega,\label{eq:flow_map2}\\
&(x,\tau)\in D\times\mathbb{R}_{\geq0},~~~~~x^+\in G(x),~~~~~\tau^+=\tau\label{eq:jump_map2},
\end{align}    
\end{subequations}
where $\omega>0$ dictates the rate of evolution of the auxiliary timer state $\tau$. Note that \eqref{eq:HDSmodel2} can be written as \eqref{eq:HDSmodel} using the overall state $y=(x,\tau)$.  Also, note that for hybrid systems with periodic flows, it is also possible to consider models with bounded timers by incorporating periodic resets in $\tau$, similar to Example \ref{example1}.
\subsection{The Hybrid Basic Conditions}
A crucial factor in establishing the stability properties of hybrid set-seeking systems is ensuring robustness to perturbations in the states and dynamics. To guarantee this property, we will impose the appropriate regularity conditions on the data of the hybrid system \eqref{eq:HDSmodel}. For differential equations, such conditions are reduced to the continuity of the vector field. However, to study dynamic inclusions we will consider set-valued maps that are outer-semicontinuous (OSC) and locally bounded (LB). A set-valued mapping $M:\mathbb{R}^p\rightrightarrows\mathbb{R}^n$ is OSC at $z\in\mathbb{R}^p$ if for each convergent sequence $\{z_i,s_i\}\to(z,s)\in\mathbb{R}^p\times\mathbb{R}^n$ satisfying $s_i\in M(z_i)$ for all $i\in\mathbb{Z}_{\geq0}$, we have $s\in M(z)$. A mapping $M$ is locally bounded (LB) at $z$ if there exists an open neighborhood $N_z\subset\mathbb{R}^p$ of $z$ such that  $M(N_z)$ is bounded. The mapping $M$ is OSC and LB relative to a set $K\subset\mathbb{R}^p$ if the mapping from $\mathbb{R}^p$ to $\mathbb{R}^n$ defined by $M(z)$ for $z\in K$ and $\emptyset$ for $z\notin K$, is OSC and LB at each $z\in K$. Using OSC and LB, as well as standard convexity and closeness notions for sets \cite{Rockafellar}, we consider hybrid systems \eqref{eq:HDSmodel} that satisfy the following  \emph{Hybrid Basic Conditions}:
\begin{figure}[b]
\begin{tcolorbox}[colback=iceblue!50, colframe=iceblue!50]
\includegraphics[width=\linewidth]{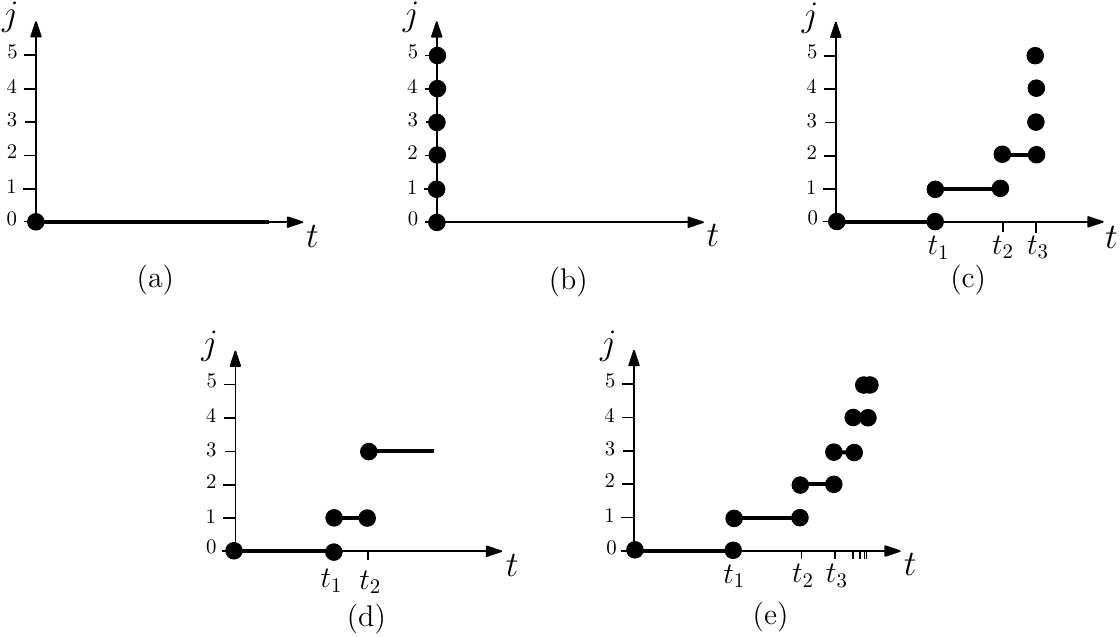}
\end{tcolorbox}
 \caption{Five different hybrid time domains illustrating different types of solutions in a hybrid system: (a) Purely continuous solution; (b) Purely discrete solution; (c) Eventually discrete solution; (d) Solution with average-dwell time; (e) Solution with Zeno behavior.   \label{HTD}}
\end{figure}
\begin{enumerate}[(a)]
\item  The sets $C$ and $D$ are closed.
\item  The set-valued mapping $F$ is OSC and LB relative to $C$, and for each $x\in C$ we have that $F(x)$ is non-empty and convex.
\item The set-valued mapping $G$ is OSC and LB relative to $D$, and for each $x\in D$ we have that $G(x)$ is non-empty.
\end{enumerate}
Outer-semicontinuity of a set-valued mapping $M:\mathbb{R}^p\to\mathbb{R}^n$ relative to a set $K\subset\mathbb{R}^p$ is equivalent to asking for closedness of the graph of $M$ with respect to $K$, i.e., the set $\text{gph(M)}:=\{(z\in \mathbb{R}^p,s\in\mathbb{R}^n):z\in K, s\in M(z)\}$ is relatively closed in $K\times\mathbb{R}^n$ \cite[Lemma 5.10]{bookHDS}. Throughout this paper, we will work only with hybrid systems that satisfy the Hybrid Basic Conditions. In the context of ES, the controllers will always be designed to satisfy these conditions.
\begin{figure}[t!]
\begin{tcolorbox}[colback=iceblue, colframe=iceblue]
\includegraphics[width=\linewidth]{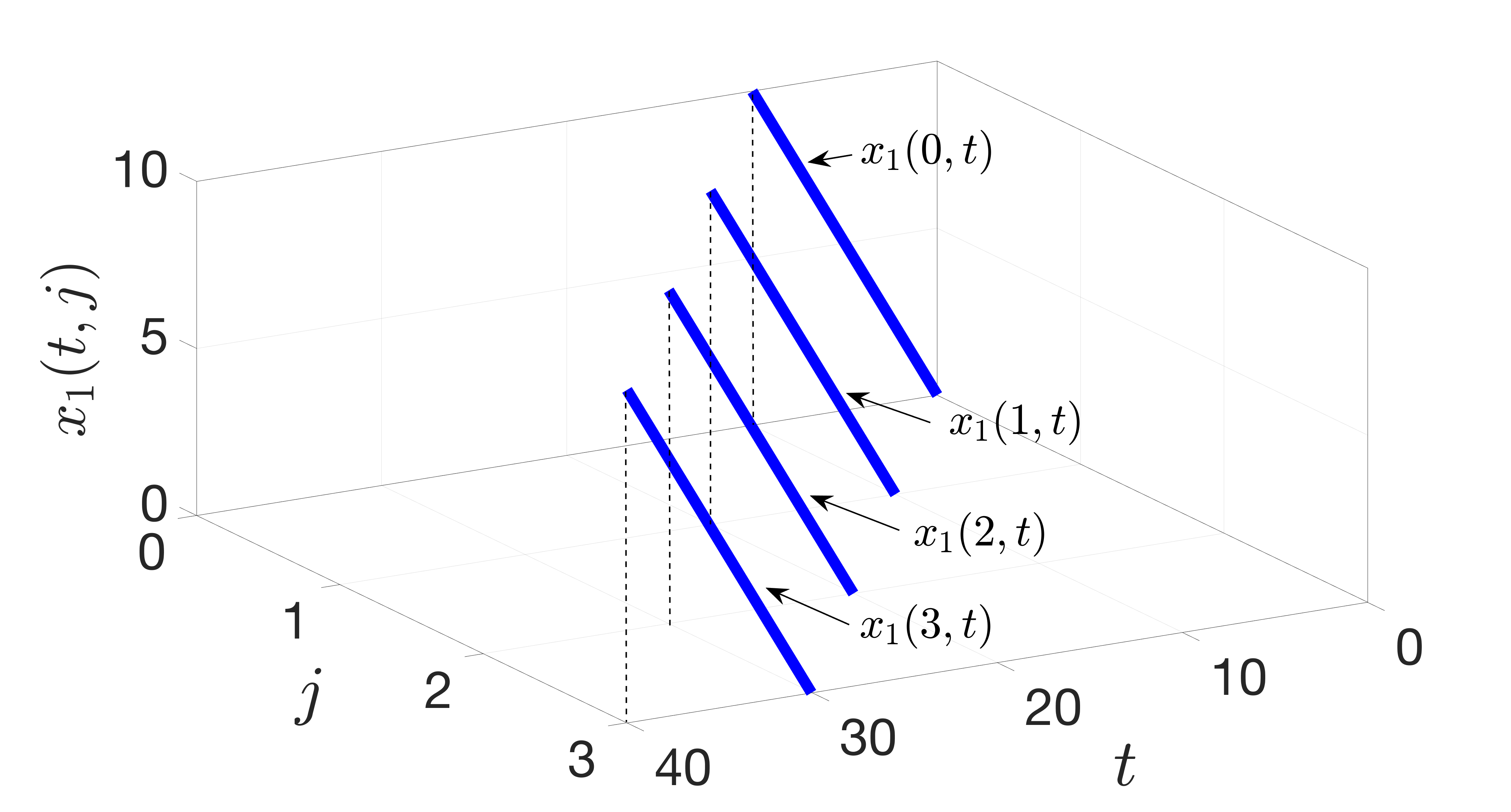}
\includegraphics[width=\linewidth]{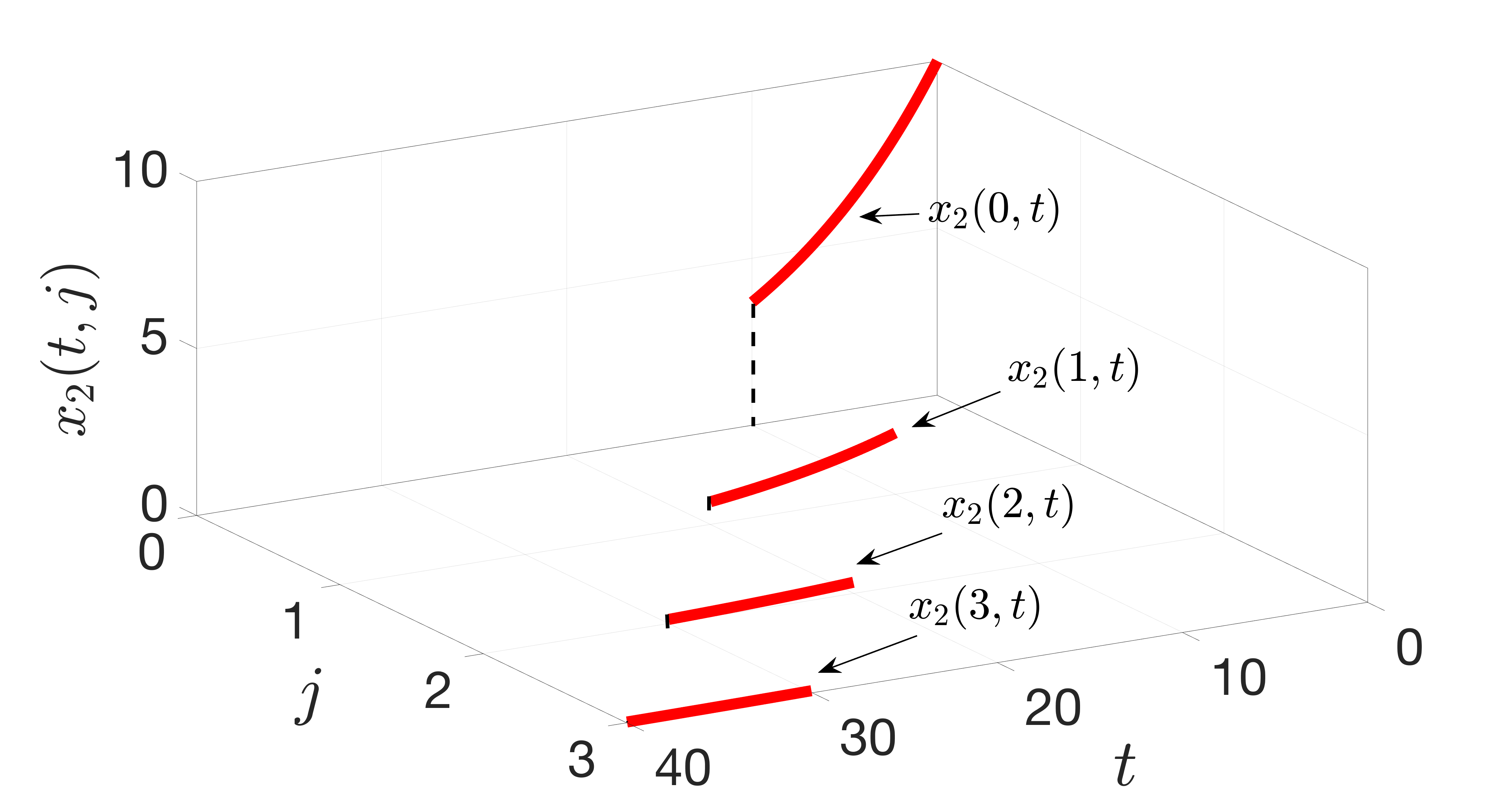}
\end{tcolorbox}
 \caption{Hybrid arcs of a solution of the system studied in Example 1, illustrated as a function of time $t$ in Figure \ref{fig1}.\label{fig5c}}
\end{figure}

\begin{sidebar}{Discrete-Time Set-Valued Dynamical Systems}
\sdbarinitial{D}iscrete-time dynamical systems modeled as \emph{difference inclusions} are common in optimization, estimation, and control algorithms. They are also relevant for the study of robustness properties in discrete-time systems with a discontinuous right-hand side. Here, we review some fundamentals on difference inclusions. %Additional material can be found in \cite{SAA1,SAA2}.

\setcounter{sequation}{0}
\renewcommand{\thesequation}{S\arabic{sequation}}
\setcounter{stable}{0}
\renewcommand{\thestable}{S\arabic{stable}}
\setcounter{sfigure}{0}
\renewcommand{\thesfigure}{S\arabic{sfigure}}

\section{Difference Inclusions}
Consider the constrained difference inclusion of the form
\begin{equation}\label{DisI0}
x\in D,~~~x^+\in G(x),
\end{equation}
where $D\subset\mathbb{R}^n$, and $G:\mathbb{R}^n\rightrightarrows\mathbb{R}^n$ is a set-valued mapping. System \eqref{DisI0} is said to satisfy the \emph{Basic Conditions} if $D$ is a closed set, $G$ is outer-semicontinuous and locally bounded relative to $D$, and $G(x)$ is not empty for each $x\in D$. Given $J\in\mathbb{Z}_{\geq1}$, a sequence $x:\{0,1,\ldots,J\}\to\mathbb{R}^n$ is said to be a solution of \eqref{DisI0} if it satisfies
\begin{equation}
    ~x(j)\in D,~~~\text{and}~~~x(j+1)\in G(x(j)),
\end{equation}
for all $j\in \{0,1,\ldots,J-1\}$. For system \eqref{DisI0}, the existence of solutions is always guaranteed provided $D\neq\emptyset$. In general, given a solution $x:\text{dom}(x)\to\mathbb{R}^{n}$, it will not be the unique solution from $x(0)$ if $G(x(j))$ is not single valued for some $j$ in $\text{dom}(x)$ such that $x(j)\in D$.

\section{Krasovskii Solutions to Discrete-Time Systems}
Difference inclusions of the form \eqref{DisI0} are relevant for the study of robustness properties in discontinuous discrete-time systems. In particular, given a set $D\subset\mathbb{R}^n$ and a function $g:\mathbb{R}^n\to\mathbb{R}^n$, a sequence $x:\{0,1,\ldots,J\}\to\mathbb{R}^n$ is said to be a Krasovskii solution of the constrained difference equation 
\begin{equation}\label{differenceequation}
x\in D,~~x^+=g(x), 
\end{equation}
if it satisfies:
\begin{equation}\label{disKrasovskii}
x(j)\in \overline{D},~~\text{and}~~x(j+1)\in G_K(x):=\bigcap_{\epsilon>0}\overline{g((x+\epsilon\mathbb{B})\cap D)},
\end{equation}
for all $j\in\{0,1,\ldots,J\}$ \cite[Def. 4.13]{bookHDS}. When $D$ is closed, the sets $\overline{D}$ and $D$ coincide. Similarly, when $g$ is continuous, the mappings $G_K$ and $g$ coincide. However, when $g$ is discontinuous, the map $G_K$ is set valued at the points of discontinuity. In particular, $G_K$ describes the outer-semicontinuous hull of $g$, namely the unique set-valued mapping $G_K$ satisfying $\text{graph}(G_K)=\overline{\text{graph}(g)}$. The construction in \eqref{disKrasovskii} also holds when $g$ is a set-valued mapping. In this case, by construction (see \cite[Ex.6.6]{bookHDS}), if $g$ is locally bounded relative to $\overline{D}$, then the difference inclusion \eqref{disKrasovskii} also satisfies the Basic Conditions.
\begin{example}[Regularization of Sign Function]
Consider the discrete-time system given by
\begin{equation}
x^+=x-\alpha\cdot \text{sign}(x),~x\in\mathbb{R}^n,
\end{equation}
for $\alpha>0$, were the sign function is defined as
\begin{equation}
\text{sign}(x)=\left\{\begin{array}{cl}
1 & \text{if}~~x>0\\
0 & \text{if}~~x=0\\
-1 & \text{if}~~x<0
\end{array}\right.
\end{equation}
To obtain $G_K$, in this case it suffices to close the graph of \text{sign}(x). We obtain the set-valued mapping:
\begin{equation}
\overline{\text{sign}}(x)=\left\{\begin{array}{cl}
1 & \text{if}~~x>0\\
\{-1,0,1\} & \text{if}~~x=0\\
-1 & \text{if}~~x<0
\end{array}\right.
\end{equation}
and the difference inclusion:
\begin{equation}
x^+\in G_K(x):=x-\alpha\cdot\overline{\text{sign}}(x),~~x\in\mathbb{R}^n.
\end{equation}
\end{example}
\section{Robustness: Discrete-Time Hermes Solutions}
As in the continuous-time case, the limiting effect of small additive state disturbances on difference equations can be captured using the notion of Hermes solutions. A function $x:\{0,1,\ldots,J\}\to\mathbb{R}^n$ is said to be a Hermes solution of the constrained difference equation \eqref{differenceequation} if there exists a sequence of functions $x_i:\{0,1,\ldots,J-1\}\to\mathbb{R}^n$ and $e_i:\{0,1,\ldots,J-1\}\to\mathbb{R}^n$ such that for all $j\in\{0,1,\ldots,J-1\}$:
\begin{align*}
&x_i(j)+e_i(j)\in D,~~~~x_i(j+1)=g(x_i(j)+e_i(j)),
\end{align*}
and the sequence $\{x_i\}_{i=1}^{\infty}$ converges uniformly to $x$ on $j\in\{0,1,\ldots,J\}$, and, moreover, the sequence $e_i:\{0,1,\ldots,J\}\to\mathbb{R}^n$ converges uniformly to the zero function on $\{0,1,\ldots,J\}$ \cite[Def. 4.12]{bookHDS}. As in the continuous-time case, there is a close connection between discrete-time Hermes and Krasovskii solutions, which is stated in the next lemma, corresponding to \cite[Thm. 4.17]{bookHDS}. 
\begin{lemma}[Discrete Time Krasovskii Solutions]\label{discretetimekrasovskii}
Consider the function $x:\{0,1,\ldots,J\}\to\mathbb{R}^n$, for some $J\in\mathbb{Z}_{\geq1}$, and suppose that $g$ is locally bounded. Then $x$ is a Hermes solution of the constrained difference equation \eqref{differenceequation} if and only if it is a Krasovskii solution of the same system.
\end{lemma}
The result of Lemma \ref{discretetimekrasovskii} indicates that in order to study the robustness properties of difference equations with a discontinuous right-hand side, we should focus on studying their Krasovskii solutions. For further details on the study of difference inclusions, we refer the reader to references \cite{SAA1} and \cite{SAA2}. Discrete-time extremum seeking algorithms based on difference equations and inclusions have been studied in \cite{PopovicPHD} and \cite{Khong:14} using finite-difference approximation methods to estimate the gradients of cost functions. In the context of hybrid set-seeking systems via averaging, difference inclusions can be used to model synchronization mechanisms in multi-agent systems \cite{TAC21Momentum_Nash}, switching systems under dwell-time or average dwell-time constraints \cite{PoTe17Auto}, etc. 

% Here is some random text showing how references are handled for a sidebar. Reference \cite{SAA1} gives one random reference. A basic equation is
% \begin{sequation}
% f = ma,\label{seq1}
% \end{sequation}
% and the impact of \eqref{seq1} is highlighted in Figure \ref{sfig1}. Reference \cite{SAA2} gives another random reference, but \cite{AA1} and \cite{BB1} are references in the main body of the text that can be cited in the sidebar.

%
\end{sidebar}
\subsection{Hybrid Time Domains, Hybrid Arcs, and Solutions}
For continuous-time systems of the form \eqref{eq:flow_map0} or \eqref{eq:flow_map}, solutions $x:\text{dom}(x)\to\mathbb{R}^n$ are given by absolutely continuous functions with respect to time $t\in\text{dom}(x)\subset\mathbb{R}_{\geq0}$ (see "Continuous-Time Set-Valued Dynamical Systems"). Similarly, for discrete-time systems of the form \eqref{eq:jump_map0} or \eqref{eq:jump_map}, solutions $x:\text{dom}(x)\to\mathbb{R}^n$ are defined as sequences parameterized by a discrete-time index $j\in\text{dom}(x)\subset \mathbb{Z}_{\geq0}$ (see "Discrete-Time Set-Valued Dynamical Systems"). Therefore, to define solutions for hybrid systems of the form \eqref{eq:HDSmodel01} or \eqref{eq:HDSmodel}, we consider a "hybrid" parameterization of the solutions that depends on both a continuous-time index $t\in\mathbb{R}_{\geq0}$, which increases continuously during flows, and a discrete-time index $j\in\mathbb{Z}_{\geq0}$, which increases by one during jumps. In this way, solutions to \eqref{eq:HDSmodel01} or \eqref{eq:HDSmodel} are defined in \emph{hybrid time domains}.
% Solutions $x$ to \eqref{eq:HDSmodel} are parameterized by a continuous-time index $t\in\mathbb{R}_{\geq0}$, which increases continuously during the flows \eqref{eq:flow_map}, and a discrete-time index $j$, which increases by one during the jumps \eqref{eq:jump_map}. As such, the notation $\dot{x}$ in \eqref{eq:flow_map} stands for $\dot{x}=\frac{\text{d}x(t,j)}{\text{d}t}$, and the notation $x^+$ in \eqref{eq:jump_map} stands for $x^+=x(t,j+1)$. When $D=\emptyset$, system \eqref{eq:HDSmodel} recovers a continuous-time system, and in this case the index $j$ can be omitted from the solutions. Additionally, when $F$ is singled-valued and continuous, system \eqref{eq:HDSmodel} reduces to a standard ODE. In this way, system \eqref{eq:HDSmodel} provides a unifying formalism to study both continuous-time systems and hybrid dynamical systems. In some cases, we will also use a different continuous-time scale (usually, denoted by $s$) to study the behaviors of systems \eqref{eq:HDSmodel} and \eqref{eq:HDSmodel2} under fast-varying signals or dynamics. 
%
% Since the solutions to system \eqref{eq:HDSmodel} are parameterized by both continuous-time and discrete-time indexes, they are defined on \emph{hybrid time domains}.
A set $E\subset\mathbb{R}_{\geq0}\times\mathbb{Z}_{\geq0}$ is called a \textsl{compact} hybrid time domain if $E=\cup_{j=0}^{J-1}([t_j,t_{j+1}],j)$ for some finite sequence of times $0=t_0\leq t_1\ldots\leq t_{J}$. A set $E \subset \mathbb{R}_{\geq0} \times \mathbb{Z}_{\geq0}$ is a hybrid time domain if it is the union of a non-decreasing sequence of compact hybrid time domains, namely, E is the union of compact hybrid time domains $E_j$ with the property that $E_0 \subset E_1 \subset E_2 \subset \ldots \subset E_j \ldots$, etc. 

Figure \ref{HTD} shows five different types of hybrid time domains, including those corresponding to Zeno behavior, average dwell time between jumps, purely discrete solutions, and purely continuous solutions.

Using the notion of hybrid time domain, the following definition formalizes what we call \emph{solutions} to hybrid systems of the form \eqref{eq:HDSmodel}.
\begin{definition}[Hybrid Arcs and Solutions]\label{defsolutionshybridystems}
A function $x:\text{dom}(x)\mapsto\mathbb{R}^n$ is a hybrid arc if $\text{dom}(x)$ is a hybrid time domain and $t\mapsto x(t,j)$ is a locally absolutely continuous function for each $j$ such that the interval $I_j:=\{t:(t,j)\in \text{dom}(x)\}$ has non-empty interior. A hybrid arc $x$ is a \emph{solution} to \eqref{eq:HDSmodel} if $x(0,0)\in C\cup D$, and the following two conditions hold:

\begin{figure*}
\begin{tcolorbox}[colback=iceblue, colframe=iceblue]
\centering
\includegraphics[width=0.8\linewidth]{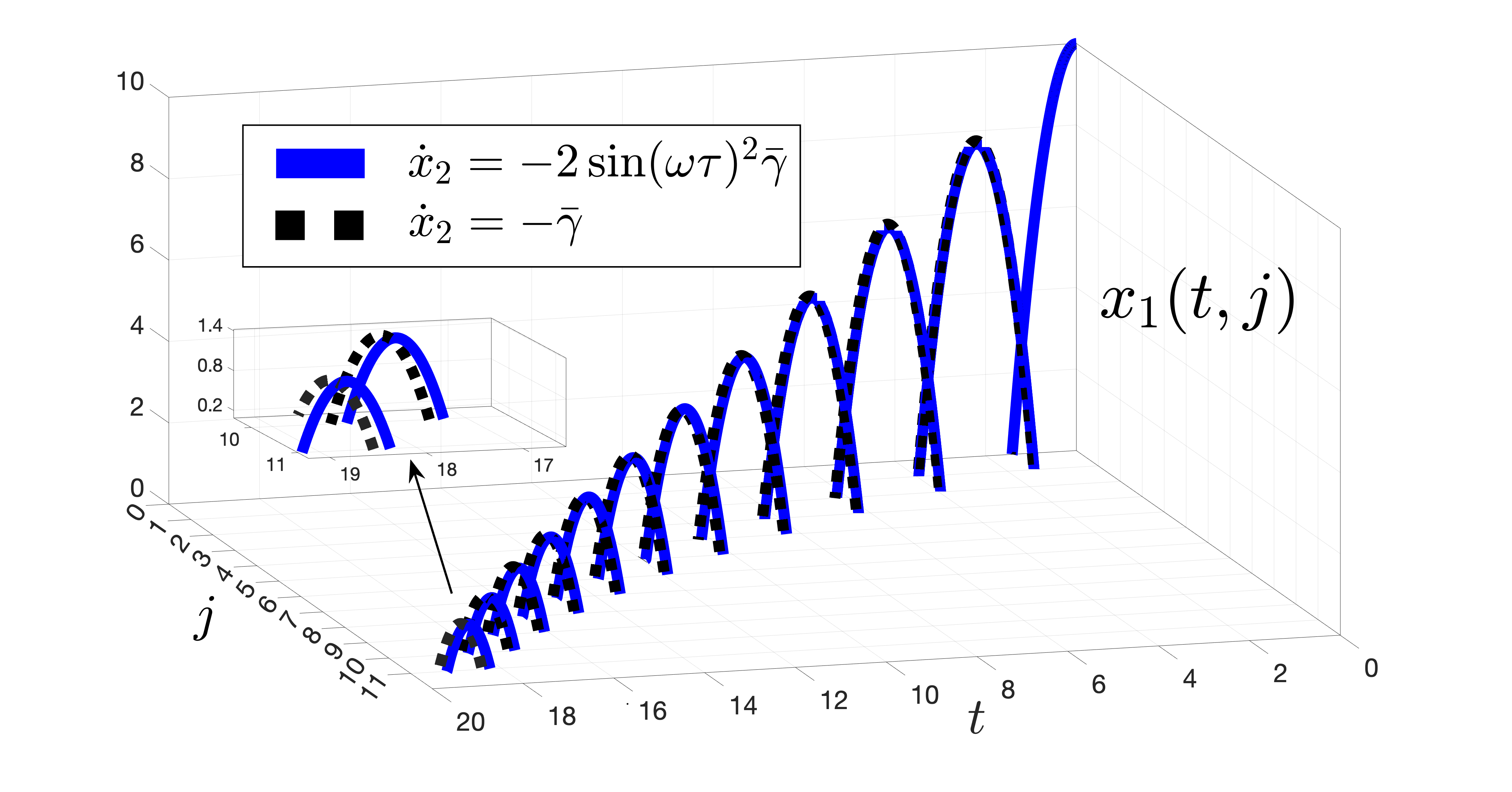}
\end{tcolorbox}
 \caption{Component \( x_1 \) of the hybrid arcs associated with the solutions of hybrid system \eqref{seekingballsystems}, as studied in Example 2, shown as a function of time \( t \) in Figure \ref{fig2}.
 The blue trajectory is generated using the highly-oscillatory acceleration $\gamma(\tau)$, with $\omega=100$. The dotted black trajectory is generated using the constant acceleration $\overline{\gamma}$, which is just the average of $\gamma(\tau)$.\label{fig5cZZ}}
\end{figure*}
\begin{enumerate}[(a)]
\item  For each $j\in\mathbb{Z}_{\geq0}$: 
\begin{align*}
&x(t,j)\in C,~ \text{for almost all}~t\in I_j\\
&\dot{x}(t,j)\in F(x(t,j)),~~\text{for almost all}~t\in I_j; 
\end{align*}
\item For each $(t,j)\in\text{dom}(x)$ such that $(t,j+1)\in \text{dom}(x)$:
\begin{equation}
x(t,j)\in D,~~\text{and}~~x(t,j+1)\in G(x(t,j)).
\end{equation}
\end{enumerate}
\end{definition}

A solution $x$ to system \eqref{eq:HDSmodel} is said to be \emph{forward pre-complete} if its domain is compact or unbounded. It is said to be \emph{forward complete} if its domain is unbounded, and it is said to be \emph{maximal} if there does not exist another solution $\psi$ to \eqref{eq:HDSmodel} such that $\text{dom}(x)$ is a proper subset of $\text{dom}(\psi)$, and $x(t,j)=\psi(t,j)$ for all $(t,j)\in\text{dom}(x)$. In general, given an initial condition $x(0,0)\in C\cup D$, the solutions to system \eqref{eq:HDSmodel} might not be unique. This will be the case if, for example, the intersection $C\cap D$ is not empty and from a point in this intersection it is possible to flow, or if the flow or jump maps in \eqref{eq:flow_map}-\eqref{eq:jump_map} admit non-unique solutions due to set-valuedness, or lack of Lipschitz continuity in the flow map.
\begin{example}[Periodic Hybrid System (Continued)] Recall the hybrid system discussed in Example 1, which implements periodic jumps in the state $x_2$ triggered by a periodic resetting timer $x_1$ with resetting period $T$. Figure \ref{fig5c} shows the hybrid arcs corresponding to the trajectories $x_1$ and $x_2$ shown in Figure \ref{fig1}, from the initial conditions $x_1(0,0)=0$ and $x_2(0,0)=10$. As observed, the values of the states $x_1$ and $x_2$ can be fully identified via the hybrid time indices for all $(t,j)\in\text{dom}(x)$. \QEDB 
\end{example}

\begin{example}[Bouncing Seeking System (Continued)]
Recall the ``bouncing'' seeking system discussed in Example 2. Figure \ref{fig5cZZ} shows the hybrid arcs corresponding to the trajectories $x_1$ shown in the left plot of Figure \ref{fig2}. The trajectory in color black is generated by the hybrid dynamics \eqref{seekingballsystems} using constant acceleration, i.e.,  $\dot{x}_2=-\overline{\gamma}$. On the other hand, the blue trajectory implements highly-oscillatory acceleration, i.e.,  $\dot{x}_2=-2\sin(\omega \tau)^2\overline{\gamma}$. \tcb{Both trajectories were generated from the same initial conditions as in Example 2.} As observed, both hybrid arcs remain ``close'' to each other during the simulation. \tcb{Specifically, the \emph{graphs} of the hybrid arcs are similar.} This hints to some closeness property, akin to the one discussed when studying system \eqref{ESCODEVanilla} and its average dynamics \eqref{perturbedgradientflow01} using inequality \eqref{closeness2}. In fact, by working with hybrid arcs, hybrid time domains, and a suitable metric for these objects, it is possible to naturally extend the notion of $(T,\epsilon)$-closeness of solutions to hybrid seeking systems, thus opening the door to the development and study of novel algorithms and applications with hybrid dynamics in the loop.\QEDB 
\end{example}
% \vspace{0.2cm}
% \noindent 
% \textbf{Example 3 (Continued)} For the multi-source seeking problem discussed in Example 3, jumps occur only sporadically and are separated by intervals of flow of at least.... The emerging behavior of the solutions of this system is shown in Figure ...
%
\begin{pullquote}
The use of hybrid time domains and graphical convergence to study solutions in hybrid systems is the key enabling mathematical tool for the development of a comprehensive theory in hybrid set-seeking control.
\end{pullquote}
\subsection{Graphical Convergence and $(T,\epsilon)$-Closeness in Hybrid Systems}
\begin{figure}[t!]
\begin{tcolorbox}[colback=iceblue, colframe=iceblue]
\includegraphics[width=\linewidth]{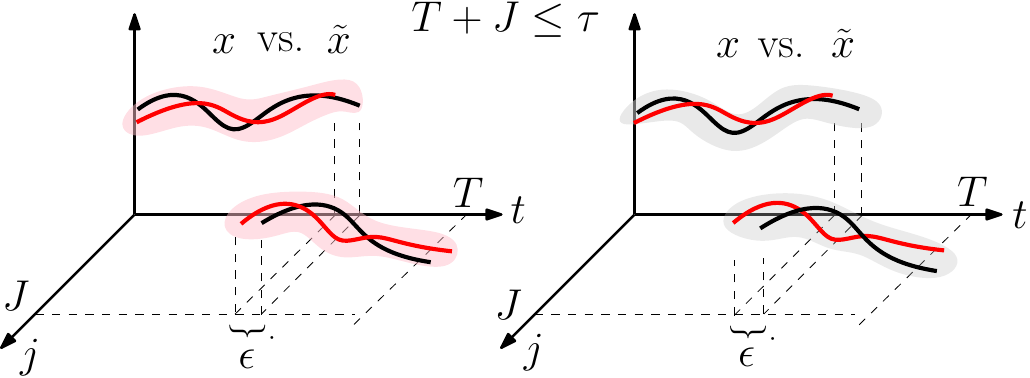}
\end{tcolorbox}
\caption{Graphical illustration of the property of $(\tau,\epsilon)$-closeness in two hybrid arcs $x$ and $\tilde{x}$. \label{fig10es}}
\end{figure}

As discussed in the analysis of the smooth ES dynamics \eqref{ESCODEVanilla}, using the average dynamics \eqref{perturbedgradientflow01} and the nominal  dynamics \eqref{gradientflow01}, the closeness-of-solutions property is crucial for the study of ES controllers based on periodic dithering. Specifically, inequalities \eqref{closeness1} and \eqref{closeness2} showed that the solutions of the seeking dynamics \eqref{ESCODEVanilla} behave similarly to those of a gradient flow, with this "similarity" being quantified, in compact sets and time domains, using the uniform distance between solutions. However, this approach does not extend directly to hybrid systems. The following example, which continues the discussion from Example \ref{seekingballexample1}, illustrates some of the challenges that arise when using the standard uniform distance to gauge the closeness between solutions of nominal and perturbed hybrid systems. These challenges can be addressed by working with solutions defined on hybrid time domains, as opposed to solutions defined on $\mathbb{R}_{\geq0}$, and by studying the closeness of solutions using the distance between the graphs of the hybrid arcs.
\begin{example}[Bouncing Seeking System (continued)]\label{example4}
Consider the ``bouncing'' seeking system with hybrid dynamics \eqref{seekingballsystems}, which experiences a high-frequency oscillating acceleration given by $-\gamma(\tau)$. Using the fact that the average of $\sin(\tau)^2$ over one period of oscillations is equal to $\frac{1}{2}$, the average dynamics of \eqref{flowsseekingball} can be explicitly computed (see "Averaging Theory: From ODEs to Hybrid Systems") to be
\begin{subequations}\label{averageballsystems}
\begin{align}
C&:=\{\bar{x}\in\mathbb{R}^2: \bar{x}_1\geq0\},~~\dot{\bar{x}}=
\bar{f}(\bar{x}):=\left(\begin{array}{c}
\bar{x}_2\\
-\overline{\gamma},
\end{array}\right),\\
D&=\{\bar{x}\in\mathbb{R}^2: \bar{x}_1=0,~\bar{x}_2\leq0\},~
\bar{x}^+=g(\bar{x}):=\left(\begin{array}{c}
0\\
- \lambda \bar{x}
\end{array}\right),
\end{align}
\end{subequations}
where $\bar{x}=(\bar{x}_1,\bar{x}_2)\in\mathbb{R}^2$ is the state of the average hybrid system and $\overline{\gamma}$ is the same gravitational constant of \eqref{oscillatorygravity}. We can argue, at least visually based on Figure \ref{fig2}, that there exist $T,\varepsilon>0$ such that the solutions of the oscillating hybrid dynamics \eqref{seekingballsystems} and the non-oscillating hybrid dynamics \eqref{averageballsystems} in general do not satisfy bounds of the form \eqref{closeness1} or \eqref{closeness2}, no matter how large we select $\omega$ in \eqref{flowsseekingball}. For example, this behavior is evident in the right plot of Figure~\ref{fig2}, where after approximately 15 seconds a significant mismatch emerges between the jump times of the red and black trajectories. Specifically, following a jump in the black trajectory, the red trajectory fails to closely track it due to the discontinuity, regardless of how large \( \omega \) is chosen. This occurs because the red trajectory has not yet undergone its corresponding jump at that moment, preventing the satisfaction of inequalities of the form \eqref{closeness1} or \eqref{closeness2} at that time.
\hfill \QEDB
\end{example}
To overcome the issues described in the previous example, and to study the closeness of solutions of hybrid systems of the form \eqref{eq:HDSmodel}, we consider an equivalent notion that compares the \emph{distance} between solutions defined under the same jump index $j$ for which both solutions are defined \cite[Def. 5.23]{bookHDS}, i.e., we use the distance between \emph{their graphs}. Specifically, given constants $\tau\geq 0$ and $\epsilon>0$, two hybrid arcs $x$ and $\tilde{x}$ are said to be $(\tau,\epsilon)$-close if the following conditions are satisfied:
\begin{enumerate}
\item For each $(t,j)\in\text{dom}(x)$ with $t+j\leq \tau$, there exists $s\in\mathbb{R}_{\geq0}$ such that $(s,j)\in\text{dom}(\tilde{x})$, $|t-s|\leq \epsilon$, and
\begin{equation}
|x(t,j)-\tilde{x}(s,j)|<\epsilon.
\end{equation}
\item For each $(t,j)\in\text{dom}(\tilde{x})$ with $t+j\leq \tau$, there exists $s\in\mathbb{R}_{\geq0}$ such that $(s,j)\in\text{dom}(x)$, $|t-s|\leq \epsilon$, and 
\begin{equation}
|\tilde{x}(t,j)-x(s,j)|<\epsilon.
\end{equation}
\end{enumerate}
Figure \ref{fig10es} presents an illustration of this property, which essentially asks that the \emph{graphs} of $x$ and $\tilde{x}$ are "close" to each other when truncated in (hybrid) time. This property can also be leveraged to study a type ``graphical'' convergence between hybrid arcs, akin to the standard uniform convergence properties studied for solutions of ODEs. In particular, given a sequence of hybrid arcs $\{x_i\}_{i=1}^{\infty}$ with $x_i:\text{dom}(x_i)\to\mathbb{R}^n$, and a hybrid arc $x$ with $x:\text{dom}(x)\to\mathbb{R}^n$, the property that for each $\tau\geq0$ and each $\epsilon>0$ there exists $i^*\in\mathbb{Z}_{\geq1}$ such that for all $i>i^*$ the hybrid arcs $x_i$ and $x$ are $(\tau,\epsilon)$-close implies the property that the sequence $\{x_i\}_{i=1}^{\infty}$ converges graphically to $x$, i.e., the sequence of sets $\{\text{gph}~x_i\}_{i=1}^{\infty}$ converges (in the sense of set convergence, see \cite[Def. 5.1]{bookHDS}) to the set $\{\text{gph}~x\}$. The converse implication is also true whenever the sequence $\{x_i\}_{i=1}^{\infty}$ is locally eventually bounded\footnote{A sequence of hybrid arcs is said to be locally eventually bounded if for any $m>0$, there exists $i_0^*\in\mathbb{Z}_{\geq1}$ and a compact set $K$ such that for all $i>i_0^*$ and all $(t,j)\in\text{dom}(x_i)$ with $t+j\leq m$, we have that $x_i(t,j)\in K$.} \cite[Thm. 5.25]{bookHDS}.

By working with hybrid time domains, we can now see (at least at this point, visually based on Figure \ref{fig5cZZ}) that the hybrid arcs describing the two solutions considered in Example \ref{seekingballexample1} are $(\tau,\epsilon)$-close when $\omega$ is sufficiently large. However, to formalize this property for a broad class of hybrid systems of interest, including those with additive perturbations on the states and dynamics (as in the perturbed gradient flow \eqref{perturbedgradientflow01}), we need to introduce an ``inflated'' version of hybrid systems of the form \eqref{eq:HDSmodel}, which will be particularly useful in the context of averaging (see "Averaging Theory: From ODEs to Hybrid Systems") and singular perturbation theory (see "Singularly Perturbed Hybrid Dynamical Systems").
\subsection{Inflated Hybrid Systems}
Given a hybrid system $\mathcal{H}$ of the form \eqref{eq:HDSmodel} with state $x\in\mathbb{R}^n$, a continuous function $\sigma:\mathbb{R}^n\to\mathbb{R}^n$, and $\rho>0$, consider the following "inflated" hybrid system $\mathcal{H}_{\rho}$ with state $x\in\mathbb{R}^n$, and data constructed from the data of \eqref{eq:HDSmodel} as follows:
\begin{subequations}\label{eqinflatedHDS}
\begin{align}
C_{\rho}&:=\left\{x\in\mathbb{R}^n: (x+\rho\sigma(x)\mathbb{B})\cap C\neq\emptyset\right\},\\
F_{\rho}(x)&:=\overline{\text{co}}~F\left((x+\rho\sigma(x)\mathbb{B})\cap C\right)+\rho\sigma(x)\mathbb{B},~\forall~x\in C_{\rho},\label{inflatedflowmapconstruction}\\
D_{\rho}&:=\left\{x\in\mathbb{R}^n: (x+\rho\sigma(x)\mathbb{B})\cap D\neq\emptyset\right\},\\
G_{\rho}(x)&:=\Big\{v\in\mathbb{R}^n:v\in g+\rho\sigma(g)\mathbb{B},\notag\\
&~~~~~~~~~g\in G((x+\rho\sigma(x)\mathbb{B})\cap D)\Big\},~\forall~x\in D_{\rho}.
\end{align}
\end{subequations}
This inflated hybrid system can capture different types of \emph{perturbations} acting on the flow and jump sets of the nominal hybrid dynamics, as well as in the nominal flow and jump maps. For example, the inflated flow map $F_{\rho}$ and the jump map $G_{\rho}$ can capture additive disturbances acting on the states and dynamics of the flow and jump maps of system \eqref{eq:HDSmodel}. Similarly, the inflated flow set $C_{\rho}$ and jump set $D_{\rho}$ can model measurement noise or other small additive disturbances acting on the states of the nominal system \eqref{eq:HDSmodel}. As shown in the following example, the inflated system \eqref{inflatedflowmapconstruction} can also be used for the study of continuous-time seeking systems, including system \eqref{ESCODEVanilla}.

\begin{example}[Inflations in Smooth Seeking Systems]\label{exampleperturbation}
\tcr{Consider the standard ES smooth dynamics \eqref{ESCODEVanilla}, which can be written as \eqref{eq:HDSmodel0} using an auxiliary state $\tau$, and dynamics:
\begin{equation}\label{example6detailed}
\dot{\hat{u}}=f_{\varepsilon}(\hat{u},\tau):=-\frac{2}{\varepsilon_a}J(\hat{u}+\varepsilon_a\sin(\tau))\sin(\tau),~~\dot{\tau}=\omega,
\end{equation}
where $\omega=\frac{2\pi}{\varepsilon_{\omega}}$, $\varepsilon=(\varepsilon_a,\varepsilon_{\omega})\in\mathbb{R}_{>0}$ are tunable parameters, and where the flow set is given by the condition $(\hat{u},\tau)\in C:=\mathbb{R}\times\mathbb{R}_{\geq0}$. For small values of $\varepsilon_a$, a Taylor expansion of $J$  leads to 
$$f_{\varepsilon}(\hat{u},\tau)=-\left(\frac{2}{\varepsilon_a}J(\hat{u})\sin(\tau)+2\frac{\partial J(\hat{u})}{\partial\hat{u}}\sin(\tau)^2+O(\varepsilon_a)\right),$$ 
where, for a given compact set $K\subset\mathbb{R}$, the notation $O(\varepsilon_a)$ indicates high-order terms $o_h(\tau,\hat{u})$ that satisfy, for some $k>0$, the inequality $|o_h(\tau,\hat{u})|\leq k\varepsilon_a$ for all $\tau\in\mathbb{R}_{\geq0}$ and all $\hat{u}\in K$. Since $\int_{0}^{\frac{2\pi}{\omega}}\sin(\tau)d\tau=0$ and $\frac{\omega}{2\pi}\int_{0}^{\frac{2\pi}{\omega}}\sin(\tau)^2d\tau=\frac{1}{2}$, it follows that the average map $f_{\text{ave}}$ of $f_{\varepsilon}$, given by \eqref{integralaverage}, is precisely system \eqref{perturbedgradientflow01}, which satisfies
\begin{align*}
f_{\text{ave}}(\bar{u})=-\frac{\partial J(\hat{u})}{\partial\hat{u}}+O(\varepsilon_a)&\in \overline{\text{co}}~\left(\frac{\partial J\left(\bar{u}+k\varepsilon_a\mathbb{B}\right)}{\partial\bar{u}}\right)+k\varepsilon_a\mathbb{B}\\
=F_{k\varepsilon_a}(\bar{u}),
\end{align*}
for all $\bar{u}\in K$, where $F_{k\varepsilon_a}$ is precisely the set-valued map defined in \eqref{inflatedflowmapconstruction} with $F=-\frac{\partial J}{\partial\bar{u}}$, $C=\mathbb{R}$, a constant function $\sigma=k$, and $\rho=\varepsilon_a$. The above inclusion implies that, on compact sets $K$, the set of solutions of system $\dot{\bar{u}}=f_{\text{ave}}(\bar{u})$ is contained in the set of solutions of the differential inclusion $\dot{\bar{u}}\in F_{k,\varepsilon_a}(\bar{u})$. Hence, by verifying appropriate properties for \emph{all solutions} of this differential inclusion (which is just an inflation of the model-based gradient flow), we can guarantee that these properties also hold for all solutions of $\dot{\bar{u}} = f_{\text{ave}}(\bar{u})$}

\QEDB
\end{example}
By using the construction of the inflated hybrid system \eqref{eqinflatedHDS}, the following key technical property, corresponding to \cite[Prop. 6.34]{bookHDS}, opens the door to a perturbation-based methodology to study hybrid set-seeking systems, thus generalizing the approach used to study the smooth ES algorithm \eqref{ESCODEVanilla}. We recall that throughout this paper we assume that the hybrid system $\mathcal{H}$ given by \eqref{eq:HDSmodel} always satisfies the Hybrid Basic Conditions by design.

\begin{lemma}[Key Technical Lemma]\label{keyclosenesshybrid}
Suppose that $x_0\in\mathbb{R}^n$ is such that each maximal solution of $\mathcal{H}$ from $x_0$ is complete or bounded. Then, for all $\tau\geq0$ and $\epsilon>0$ there exists $\rho>0$ such that for each solution $x_{\rho}$ to the perturbed HDS \eqref{eqinflatedHDS}, with $|x_{\rho}(0,0)-x_0|\leq\rho$, there exists a solution $x$ to the nominal hybrid system $\mathcal{H}$, given by \eqref{eq:HDSmodel}, with  $x(0,0)=x_0$, such that $x$ and $x_{\rho}$ are $(\tau,\epsilon)$-close.   \QEDB 
\end{lemma}
%
%By Lemma \ref{keyclosenesshybrid}, for each solution of the perturbed system \eqref{perturbedgradientflow01}
%we obtain $(\tau,\epsilon)$-closeness with respect to some solution of the nominal system \eqref{gradientflow01}, thus linking the qualitative behaviors of the solutions of both systems.

Lemma \ref{keyclosenesshybrid} also applies to differential equations simply by taking $D=\emptyset$, $C=\mathbb{R}^n$, and $F$ being a continuous function. \tcr{For instance, for the smooth ES system considered in Example \ref{exampleperturbation}, we can use Lemma \ref{keyclosenesshybrid} to directly  conclude that, on compact sets, for each solution $\bar{u}_{\varepsilon_a}$ of the differential inclusion $\dot{\bar{u}}_{\varepsilon_a}\in F_{k,{\varepsilon_a}}(\bar{u}_{\varepsilon_a})$, there exists a solution $\tilde{u}$ of the model-based gradient flow $\dot{\tilde{u}}=-\frac{\partial J(\tilde{u})}{\partial\tilde{u}}$ such that $\bar{u}_{\varepsilon_a}$ and $\tilde{u}$ are $(\tau,\varepsilon)$-close. Since, as shown before, on compact sets we have $f_{\text{ave}}\in F_{k,\varepsilon_a}$, it follows that the $(\tau,\varepsilon)$-closeness property also holds for the solutions of the average dynamics of \eqref{example6detailed} and the solutions of the target gradient flow. This analytical methodology extends beyond gradient flows and smooth
dynamics to general hybrid set-seeking systems, enabling the
establishment of closeness-of-solutions properties via perturbation
theory for hybrid systems. In this context, the use of hybrid time domains and graphical convergence in analyzing hybrid system solutions provides a foundational mathematical framework for developing a comprehensive theory of hybrid set-seeking control.} %We formalize this property in the following definition:

%When the nominal hybrid system \eqref{eq:HDSmodel} satisfies the hybrid basic conditions, and possess suitable uniform stability properties, Lemma \ref{keyclosenesshybrid} can be used to establish that similar stability properties will hold for the inflated hybrid system \eqref{eqinflatedHDS}.
%

\subsection{Stability Notions for Hybrid Set-Seeking Systems}
Since hybrid systems usually incorporate logic modes, timers, clocks, and other auxiliary states that do not settle into a point but rather evolve on bounded sets, we study the uniform stability properties of hybrid set-seeking systems with respect to general compact sets.
\begin{definition}[Uniform Global Asymptotic Stability]\label{Def:UGAS}
A compact set $\mathcal{A}\subset\mathbb{R}^n$ is said to be \emph{uniformly globally asymptotically stable} (UGAS) for system \eqref{eq:HDSmodel} (or \eqref{eq:HDSmodel2}) if there exists a class-$\mathcal{K}$$\mathcal{L}$ function $\beta$ such that any maximal solution of $\mathcal{H}$ satisfies the bound
\begin{equation}\label{KLbound1}
|x(t,j)|_{\mathcal{A}}\leq \beta(|x(0,0)|_{\mathcal{A}},t+j),~~\forall~(t,j)\in\text{dom}(x).
\end{equation}
If there exist positive constants $c_1,c_2>0$ such that $\beta(r,s)=c_1re^{-c_2s}$, then $\mathcal{A}$ is said to be \emph{uniformly globally exponentially stable} (UGES).
\QEDB 
\end{definition}
When \(\mathcal{A} = \{0\}\), the definition \ref{Def:UGAS} recovers the standard asymptotic stability notions studied for differential equations \cite[Ch. 4]{KhalilBook} and smooth ES systems \cite{Krstic2000,NesicChina}. The use of \(\mathcal{K}\mathcal{L}\) bounds in \eqref{KLbound1} provides uniformity in terms of rates of convergence and overshoots. Note that Definition \ref{Def:UGAS} requires that \emph{all} solutions to the HDS \eqref{eq:HDSmodel} satisfy the bound \eqref{KLbound1} at all times in the domain of the solution. However, it does not require that all solutions be complete. For hybrid systems, the completeness of solutions is usually established independently.
\begin{example}[Dynamics with Vanishing Coefficients]\label{example9}
\tcr{Under suitable assumptions, gradient flows and other optimization dynamics typically satisfy the UGAS property in Definition \ref{Def:UGAS}, relative to the set of minimizers of a cost function $J$. However, this assumption does not hold for all gradient systems, even if $J$ is strongly convex. For instance, systems with vanishing coefficients of the form
\begin{equation}
\dot{x}=-\alpha(\tau)\nabla J(x),~~~~ \text{or}~~~~\ddot{x}+\alpha(\tau)\dot{x}+\nabla J(x)=0,
\end{equation}
for which $\lim_{\tau\to\infty}\alpha(\tau)=0$, typically exhibit non-uniform convergence properties with respect to $\tau(0)$, thus precluding bounds of the form \eqref{KLbound1} \cite{PovedaTeelACC20}. However, under appropriate hybrid regularizations of $\tau$, e.g., time-triggered \cite{Poveda_Li:2019_CDC,PovedaNaliAuto20} or event-triggered restarting \cite{teel2019first}, the resulting hybrid dynamics can satisfy \eqref{KLbound1} for some $\beta\in\mathcal{K}\mathcal{L}$, see also ``Set-Seeking Dynamics with Momentum and Resets''.} \QEDB 
\end{example}
Since ES systems rely on perturbation theory, the following stability property is commonly used in the analysis of seeking algorithms.

\begin{definition}[Semi-Global Practical Stability]\label{Def:SGPAS}
Consider a HDS $\mathcal{H}_{\varepsilon}$, where $\varepsilon>0$ is a tunable parameter. System $\mathcal{H}_{\varepsilon}$ is said to render a compact set $\mathcal{A}$ semi-globally practically asymptotically stable (SGPAS) if there exists a class-$\mathcal{K}\mathcal{L}$ function $\beta$ such that for each $\Delta>\nu>0$ there exists $\varepsilon^*>0$ such that for all $\varepsilon\in(0,\varepsilon^*)$ and all initial conditions $|x(0,0)|\leq \Delta$, any maximal solution of $\mathcal{H}_{\varepsilon}$ satisfies the bound
\begin{equation}\label{KLbound2}
|x(t,j)|_{\mathcal{A}}\leq \beta(|x(0,0)|_{\mathcal{A}},t+j)+\nu,~~\forall~(t,j)\in\text{dom}(x).
\end{equation}
If there exist positive constants $c_1,c_2>0$ such that $\beta(r,s)=c_1re^{-c_2s}$, then $\mathcal{A}$ is said to be \emph{semi-globally practically exponentially stable} (SGPES).  
\QEDB 
\end{definition}
The SGPAS property can be extended to hybrid systems that depend on multiple parameters $\varepsilon=(\varepsilon_1,\varepsilon_2,\ldots,\varepsilon_m)$. In this case, we use the notation SGPAS as $(\varepsilon_1,\varepsilon_2,\ldots,\varepsilon_m)\to0^+$ to denote that \eqref{KLbound2} holds with the parameters $\varepsilon_i$ sufficiently small and being selected sequentially, starting from $\varepsilon_1$ and ending with $\varepsilon_m$. While the sequential tuning of certain parameters in ES can be relaxed in some cases \cite{TanAndNesic2006Local}, to simplify our presentation we do not study such relaxations in this paper. Instead, we observe that for most seeking dynamics based on averaging---though not all (see \cite{abdelgalil2024initialization,suttner2019extremum,weber2024inferring})---whether modeled as differential equations or hybrid systems, the property of SGPAS is typically the strongest guarantee one can expect. For smooth seeking systems such as in \eqref{ESCODEVanilla}, the small parameters in \( \varepsilon \) are related to the dither amplitude \( \varepsilon_a \), the inverse of the dither frequency \( \omega \), and, when additional dynamics such as filters are used, the corresponding inverse of the filter gain. For hybrid seeking systems, we will allow $\varepsilon$ to also include parameters that are intrinsic to the hybrid dynamics, e.g., $\eta_d$ in the "Source-Seeking and Surveillance" application discussed in Example \ref{intermittentsourceseeking}.

\begin{figure}[t!]
 \begin{tcolorbox}[colback=ivoryA, colframe=ivoryA]
\centering
\includegraphics[width=\linewidth]{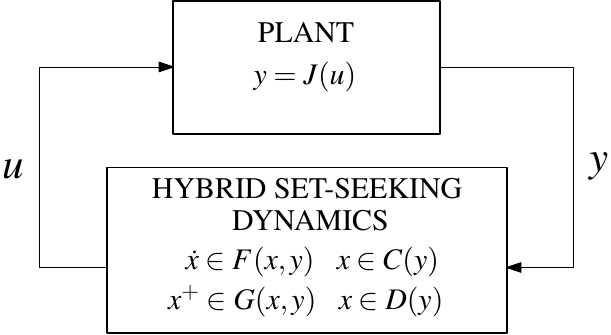}
\end{tcolorbox}
 \caption{High-level conceptual feedback scheme of hybrid set-seeking systems for static maps.\label{fig5cd0}}
\end{figure}
%

%

%\noindent\looseness=-1FIGURE 1 Title of the plot, not a sentence. Next, a sentence is given here to indicate the points that the figure is meant to highlight. Finally, you can include additional sentences to provide more detail about the meaning and importance of the figure. These sentences will greatly enhance the appeal of your article.  

\begin{sidebar}{Averaging Theory: From ODEs to Hybrid Systems}\label{sidebaraveraging}
\sdbarinitial{D}ynamical systems with time-varying oscillating vector fields are traditionally studied using averaging theory. Classic references on averaging theory for differential equations include \cite{averaging}, \cite{Verhulst_book}, \cite{TeelNesicAveraging} and \cite{KhalilBook}. For hybrid dynamical systems, averaging tools have also been studied in \cite{wang2012averaging}, \cite{AveragingHybridISS}, \cite{PovedaNaliAuto20}, \cite{abdelgalil2023multi}, and \cite{abdelgalil2023lie}. Below, we review some results on averaging theory for hybrid systems, adapted from \cite{TeelNesicAveraging}, \cite{abdelgalil2023multi}, and \cite{PovedaNaliAuto20}. 
\subsection{Hybrid Systems with Oscillating Flow Maps}
Consider a hybrid dynamical system of the form \eqref{eq:HDSmodel0}, with states $z=(x,p)\in\mathbb{R}^n\times\mathbb{R}^m$ and   $\tau\in\mathbb{R}_{\geq0}$, and dynamics 
\begin{subequations}\label{originalsystemaverage}
\begin{equation}\label{flows1average1}
    (x,p,\tau)\in C:=C_x\times C_p\times \mathbb{R}_{\geq0},~~~\left\{\begin{array}{l}
    \dot{x}=f(x,p,\tau)\\
    \dot{p}\in F_p(p)\\
     \dot{\tau}=\dfrac{1}{\varepsilon}
    \end{array}\right.
\end{equation}
\begin{equation}\label{jumps1a}
    (x,p,\tau)\in D:=D_x\times D_p\times\mathbb{R}_{\geq0},~~\left\{\begin{array}{l}
    x^+=g(x,p)\\
    p^+\in G_p(p)\\
    \tau^+=\tau
    \end{array}\right.
\end{equation}  
\end{subequations}
where $\varepsilon>0$ is a small parameter, $C_x,D_x\subset\mathbb{R}^n$, $C_p,D_p\subset\mathbb{R}^m$. The state $p$ can model logic states, timers, oscillators, as well as other auxiliary states. We assume that system \eqref{originalsystemaverage} satisfies the Hybrid Basic Conditions.  In \eqref{flows1average1}, the small parameter $\varepsilon$ introduces a time-scale separation between the dynamics of $\tau$ and those of $(x,p)$. When the vector field $f$ is periodic with respect to $\tau$, the rapid fluctuations in $\tau$ result in fast oscillatory behaviors in $\dot{x}$, which can be averaged out to approximate the dynamics of $x$ using a non-oscillating vector field. The following definition formalizes this notion of average, which for system \eqref{originalsystemaverage} involves only the function $f$.
\begin{definition}[The Average Map]\label{definitionaveragemap}
The function $f$ is said to have an \emph{average map} $\bar{f}(\cdot)$ if for each compact set $K_0\subset \mathbb{R}^n\times\mathbb{R}^m$ there exists a continuous, non-increasing function $\gamma$ satisfying $\lim_{\mathcal{T}\to\infty}\gamma(\mathcal{T})=0$, such that for all $(x,p,\tau,T)\in K_0\times\mathbb{R}_{\geq0}\times\mathbb{R}_{\geq0}$:
\begin{equation}\label{average0assumption}
\left|\frac{1}{\mathcal{T}}\int_{\tau}^{\tau+\mathcal{T}} \Big(f(x,p,s)-\bar{f}(x,p)\Big)\text{d}s\right|\leq \gamma(\mathcal{T}). 
\end{equation}
The function $\gamma(\cdot)$ is called the \emph{convergence function}. 
\end{definition}
When $f$ is periodic, as it is usually the case in ES, for each compact set $K_0$ there exists $c_0>0$ such that we can take $\gamma(\mathcal{T})=\frac{c_0}{1+\mathcal{T}}$, and it suffices to verify the integral \eqref{average0assumption} in the interval $[0,\mathcal{T}]$, where $\mathcal{T}$ is the period of the oscillations. In the rest of this section, we make the additional standing assumptions that $f$ is periodic with respect to $\tau$ and  that $\bar{f}$ is continuous. However, we note that periodicity of $f$ can be relaxed by working with generalized averages in \eqref{average0assumption}.

By using the average map $\bar{f}$, we can study \eqref{originalsystemaverage} by inspecting a "simpler" non-oscillating system, with states $(\bar{x},\bar{p})$, and the following \emph{average hybrid dynamics}:
\begin{subequations}\label{averagesystem}
\begin{equation}\label{flows1a}
    (\bar{x},\bar{p})\in C:=C_x\times C_p,~~~\left\{\begin{array}{l}
    \dot{\bar{x}}=\bar{f}(\bar{x},\bar{p})\\
    \dot{\bar{p}}\in F_p(\bar{p})
    \end{array}\right.
\end{equation}
\begin{equation}\label{jumps1a}
    (x,p)\in D:=D_x\times D_p,~~\left\{\begin{array}{l}
    \bar{x}^+=g(\bar{x},\bar{p})\\
    \bar{p}^+\in G_p(\bar{p})
    \end{array}\right..
\end{equation}  
\end{subequations}
Given that hybrid systems usually lack the uniqueness of solutions property, and might also exhibit solutions with discontinuities, standard results on closeness of solutions established in the literature of averaging for differential equations do not trivially extend to hybrid systems. Instead, we can establish that for every solution $(x,p,\tau)$ of the original hybrid system \eqref{originalsystemaverage} there exists some solution $(\bar{x},\bar{p})$ of the average hybrid dynamics \eqref{averagesystem} that is  $(T,\epsilon)$-close. In particular, the following theorem follows from the results in \cite{TeelNesicAveraging,abdelgalil2023multi,wang2012averaging}:

\vspace{-0.3cm}
\begin{theorem}[Closeness of Solutions via Averaging]\label{theoremclosenessvarying}
Consider the hybrid systems \eqref{originalsystemaverage} and \eqref{averagesystem}. Suppose that system \eqref{averagesystem} has no finite escape times. Then, for each compact set of initial conditions $K_0\subset\mathbb{R}^n$, $T\in\mathbb{R}_{>0}$ and $\epsilon\in\mathbb{R}_{\geq0}$, there exists $\varepsilon^*>0$ such that for all $\varepsilon\in(0,\varepsilon^*)$ and for all solutions $(x,p,\tau)$ to \eqref{originalsystemaverage} from $K_0$ there exists a solution $(\bar{x},\bar{p})$ to \eqref{averagesystem} such that $(\bar{x},\bar{p})$ and $(x,p)$ are $(T,\epsilon)$-close.
\end{theorem}
The result of Theorem \ref{theoremclosenessvarying} is a closeness of solutions property for hybrid systems that parallels similar properties used for the study of smooth extremum seeking systems. By using this property, as well as suitable uniform stability properties in the average hybrid dynamics, the following theorem can be established.
\vspace{-0.3cm}
\begin{theorem}[SGPAS via Averaging in HDS]\label{stabilityaveragingtime}
Consider the hybrid systems \eqref{originalsystemaverage} and \eqref{averagesystem}. Suppose that there exists a compact set $\mathcal{A}\subset\mathbb{R}^{n+m}$ that is UGAS for system \eqref{averagesystem}. Then, system \eqref{originalsystem} renders the set $\mathcal{A}\times\mathbb{R}_{\geq0}$ SGPAS as $\varepsilon\to0^+$.
\end{theorem}

\vspace{-0.3cm}
The stability result from Theorem \ref{stabilityaveragingtime} can be applied to analyze hybrid set-seeking dynamics of the form \eqref{originalsystemaverage}. This result remains valid even when the stability of the averaged dynamics \eqref{averagesystem} is only semi-globally practically asymptotically stable (SGPAS) with respect to a tunable parameter \( \delta \) \cite[Thm. 7]{PovedaNaliAuto20}. This is commonly the case in extremum-seeking control, where \( \delta \) can represent the residual terms that appear in the average dynamics, and which are of the same order as the amplitude of the dither signals. In the context of hybrid set-seeking dynamics, \( \delta \) may encompass additional parameters, offering greater flexibility in the design of the algorithms.

\vspace{-0.3cm}
\subsection{Averaging using Time-Invariant Models}
Theorems \ref{theoremclosenessvarying} and \ref{stabilityaveragingtime} can be extended to other types of models, including hybrid systems of the form
\begin{subequations}\label{originalsystem}
\begin{equation}\label{flows1a}
    (x,p,\tau)\in C:=C_x\times C_p\times \mathbb{R}_{\geq0},~~~\left\{\begin{array}{l}
    \dot{x}\in F_x(x,p)\\
    \dot{p}=f_p(p,x,\tau)\vspace{0.1cm}\\
     \dot{\tau}=\dfrac{1}{\varepsilon},
    \end{array}\right.
\end{equation}
\begin{equation}\label{jumps1a}
    (x,p,\tau)\in D:=D_x\times D_p\times\mathbb{R}_{\geq0},~~\left\{\begin{array}{l}
    x^+\in G_x(x,p)\\
    p^+\in G_p(x,p)\\
    \tau^+=\tau
    \end{array}\right.
\end{equation}  
\end{subequations}
where all the mappings and sets are selected to satisfy the Hybrid Basic Conditions. Here, the dynamics of $x$ are set-valued, but the average map can still be computed as in \eqref{average0assumption} using $f_p$ instead of $f$. 
\end{sidebar}
\section{Part 2: Hybrid Set-Seeking Systems for Static Maps}\label{ESCStatic}
Equipped with the modeling and analysis tools introduced in Part~I, we now
consider a broad class of hybrid set-seeking systems for solving model-free
decision-making problems. The primary objective in extremum seeking is to
regulate the plant input $u$ using only measurements of the output $y$, 
to drive $u$ toward the solution of an \emph{application-dependent} optimization or decision-making problem of the form:
\begin{equation}\label{optimization_problem} \text{optimize}~~J(u),~~~\text{s.t.}~~u\in \hat{\mathbb{U}}, \end{equation}
where the set $\hat{\mathbb{U}}\subset\mathbb{R}^n$ captures the constraints on the input and $J(\cdot)$ corresponds to the cost function, which is assumed to be sufficiently smooth. We assume that the set of solutions to the problem \eqref{optimization_problem} is not empty and we denote it by $\mathcal{O}\subset\mathbb{R}^n$, which is also assumed to be compact. The fact that $\mathcal{O}$ might not be a singleton further motivates the terminology "set-seeking" as opposed to the more traditional names "extremum-seeking"  or "equilibrium-seeking". 

\begin{remark}\label{remarkMAS}
The decision-making problem \eqref{optimization_problem} can also be formulated to study multi-agent networked systems (MAS) with individual outputs $y_i$ and cost functions $J_i$, for all $i\in\{1,\ldots,N\}$, where $N$ is the number of subsystems. This case is relevant for distributed optimization problems, as well as for Nash equilibrium seeking problems in game theoretic scenarios; see, for instance \cite{Frihauf12a}, \cite{Dither_ReUse}, \cite{StankovicNashSeeking}, \cite{Poveda:15}, \cite{poveda2022fixed} for different examples of seeking dynamics studied in these settings.
\end{remark}

For the remainder of this section, we assume that the plant to be optimized is characterized by a static map $J$, as illustrated in the scheme of Figure \ref{fig5cd0}. In this scheme, we use suggestive notation in which the data of the hybrid controller \( \mathcal{H} = \{C, F, D, G\} \) is allowed to depend on \( y \), the output signal of the plant. This notation indicates that the particular structure of the \emph{data} of \( \mathcal{H} \) is designed specifically for each application of interest. However, the general structure of hybrid set-seeking controllers remains largely consistent across applications.
\begin{figure*}[t!]
 \begin{tcolorbox}[colback=ivoryA, colframe=ivoryA]
\centering
\includegraphics[width=0.92\linewidth]{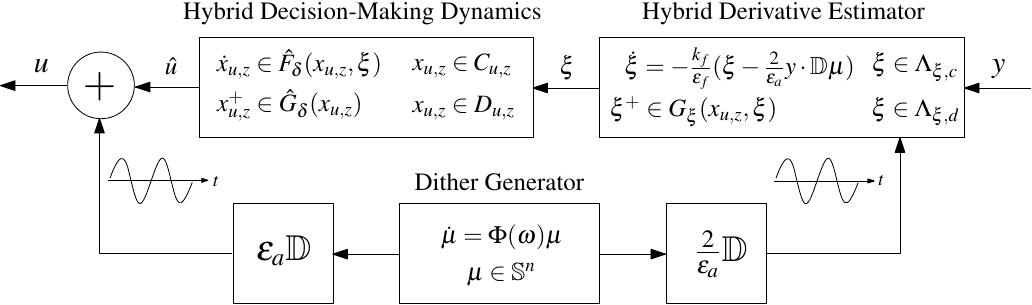}
\end{tcolorbox}
 \caption{\tcr{Main components of hybrid set-seeking systems with periodic probing.} \label{fig5cd}}
\end{figure*}
To investigate hybrid set-seeking systems with dynamic inclusions in the loop, we begin by extending the smooth seeking dynamics \eqref{ESCODEVanilla} to the setting of hybrid systems. Specifically, similar to Example \ref{intermittentsourceseeking}, and to avoid dealing with time-varying differential inclusions whose averaged maps may be difficult to compute explicitly, we consider a class of time-invariant hybrid seeking systems consisting of a \emph{dynamic derivative estimator} and \emph{dither generator}, with states \( (\xi, \mu) \in \mathbb{R}^n \times \mathbb{R}^{2n} \), and a \emph{hybrid decision-making algorithm} with internal state \( x_{u,z} := (\hat{u}, \hat{z}) \in \mathbb{R}^n \times \mathbb{R}^r \). Figure \ref{fig5cd} shows a block diagram representation of the set-seeking controllers. Our objective is to design the system’s core components so that, using only
output measurements, the plant input $u \in \mathbb{R}^n$ is driven to a small
neighborhood of the solution set of \eqref{optimization_problem}. The design
must also ensure that the resulting closed-loop system satisfies the Hybrid
Basic Conditions.

%%%%%%%%%%%%%%%%%%%%%%%%%%%%%%%%%%
\subsection{Autonomous Derivative Estimators}\label{Sec:Main}
%

%
%  \begin{figure}[t!]
%  \begin{centering}
%   \includegraphics[width=0.6\textwidth]{fig/HESC/SchemeHESC.eps}
%   \caption{Conceptual Scheme of HESC for set-valued plants.}
%   \label{scheme}
%   \end{centering}
% \end{figure}
%
%
%

The dynamic derivative estimators make use of $n$ linear oscillators of the form \eqref{dynamicoscillator}, with overall state $\mu:=(\mu_{1},\ldots,\mu_n)\in\mathbb{S}^{n}\subset\mathbb{R}^{2n}$, satisfying $\mu_i=(\mu_{i,1},\mu_{i,2})^\top\in\mathbb{S}^1$, for all $i\in\{1,\ldots,n\}$, and dynamics
\begin{equation}\label{oscillatordynamics}
\dot{\mu}=\Phi(\omega)\mu,~~~~\mu\in\mathbb{S}^n,
\end{equation}
where the $\omega$-parameterized  block matrix $\Phi_\omega:\mathbb{R}^{2n}\rightarrow\mathbb{R}^{2n}$ is given by
\begin{equation}\label{oscillator}
\Phi_\omega:=\left[\begin{array}{cccc}
\omega_1\mathbf{R}_o & \textbf{0} & \ldots & \textbf{0}\\
\textbf{0} & \omega_2\mathbf{R}_o & \ldots & \textbf{0}\\
\vdots & \vdots & \ddots & \vdots \\
\textbf{0} & \textbf{0} & \ldots & \omega_n\mathbf{R}_o\\
\end{array}\right],
\end{equation}\noindent
%
%and where the block components $\Omega_{\omega_i}$ are given by
%
%\begin{equation}\label{individual_oscillator}
%\Omega_{\omega_i}:=\omega_i\mathbf{R}_o,~~\forall~i\in\{1,\ldots,n\},
%\end{equation}\noindent
%
with $\omega=(\omega_1,\ldots,\omega_n)$. The entries of $\omega$ are selected as $\omega_i=\frac{\kappa_i}{\epsilon_{\omega}}$, where $\epsilon_{\omega}>0$ is a tunable constant, and the $\kappa_i$'s are rational numbers that satisfy $\kappa_i\neq\kappa_j$ for $i\neq j$ and $(i,j)\in\{1,\ldots,n\}$.  

The role of system \eqref{oscillatordynamics} is to generate a vector of $n$ periodic probing signals with different frequency for the purpose of real-time \emph{exploration}. This is achieved by extracting the odd entries of the vector $\mu$, given by $\mu_{i,1}(t) = \cos(\omega_i t) \mu_{2i-1}(0) + \sin(\omega_i t) \mu_{2i}(0)$, for each $i \in \{1, \ldots, n\}$, via multiplication with the matrix  $\mathbb{D}:=[e_1,\textbf{0},e_2,\textbf{0},\ldots,e_i,\ldots,\textbf{0},e_n,\textbf{0}]$, where $\textbf{0}\in\mathbb{R}^n$ corresponds to a column vector of zeros, and $e_i\in\mathbb{R}^n$ corresponds to the unitary vector with the $i^{th}$ entry equal to $1$. Since, by the definition of its flow set in \eqref{oscillatordynamics}, the oscillator is restricted to evolve in $\mathbb{S}^n$, we have $\mu_{2i-1}(0)^2+\mu_{2i}(0)^2=1$. Moreover, the set \( \mathbb{S}^n \) is forward invariant and, since the system \eqref{oscillatordynamics} has no solutions defined outside \( \mathbb{S}^n \), the set is also trivially UGAS for the dynamics \eqref{oscillatordynamics}. We note that although other types of dither signals and oscillators can be considered---including hybrid ones \cite{PovedaCDC18}---we focus here, for simplicity, on standard sinusoidal probing signals generated by linear oscillators.

The performance of the derivative estimators in ES systems can usually be enhanced by incorporating filters into the dynamics. For hybrid set-seeking, we can consider a simple low-pass filter with state $\xi\in\mathbb{R}^n$ and continuous-time dynamics
\begin{equation}\label{filterdynamics}
\dot{\xi}
=-\frac{k_f}{\varepsilon_{f}} \cdot\left(\xi-\frac{2}{\varepsilon_a} y\cdot\mathbb{D}\mu\right),
\end{equation}
where $k_f,\varepsilon_f>0$ are tunable parameters and which is allowed to evolve in a compact set $\Lambda_{\xi}:=\lambda_{\xi}\mathbb{B}\subset\mathbb{R}^n$, with $\lambda_{\xi}\in\mathbb{R}_{>0}$ selected sufficiently large to encompass all the complete solutions of interest. As shown later in Theorem \ref{theorem4} and in "Hybrid Set-Seeking Dynamics with Momentum and Resets", it is also possible to incorporate filters with hybrid dynamics to, for example, improve transient performance. However, we can simply take $\xi^+=\xi$ when no discrete-time updates are incorporated into the filters. The compactness assumption on $\Lambda_{\xi}$ can also be relaxed under additional structure on the dynamics. 
%
%Apart from the compactness assumption on the set $\Lambda_{\xi}$, system (\ref{estimator_dynamics_1}) essentially mirrors the behavior generated by the standard gradient estimators used in Lipschitz continuous extremum seeking controllers for standard optimization \cite{NesicDerivative}, \cite{GeneralFrameworkESC}, and for learning in games \cite{Dither_ReUse}, \cite{Poveda:15}. However, constraining the state $\xi$ to lie in a compact set is required in our case in order to apply averaging and singular perturbation results for hybrid systems. %Additionally, by construction, the closed-loop system will not exhibit purely discrete-time solutions.
%
%

%%%%%%%
\subsubsection{Hybrid Decision-Making Dynamics}
\label{sec:generalized_dynamics}
%%%%%%%
As mentioned in the introduction, the smooth ES algorithm \eqref{ESCODEVanilla} is designed to approximate ---on average--- the behavior of the model-based ``target'' gradient flow \eqref{gradientflow01}. In the context of hybrid set-seeking systems, we aim to consider more general hybrid target systems for optimization and decision making. Therefore, instead of \eqref{ESCODEVanilla}, we consider a class of hybrid seeking dynamics with state $x_{u,z}:=(\hat{u},\hat{z})\in\mathbb{R}^{n}\times\mathbb{R}^r$, $r\in\mathbb{Z}_{\geq0}$, where $\hat{u}$ is the main state and output of the system, and where $z\in\mathbb{R}^r$ is an auxiliary state that can model timers, logic modes, memory states, momentum, etc.

The state $x_{u,z}$ evolves according to the following \emph{hybrid decision-making dynamics}:
\begin{subequations}\label{hybrid_learning_dynamics}
\begin{align}
\dot{x}_{u,z}&\in \hat{F}_{{\delta}}\big(x_{u,z},\xi\big),~~~~x_{u,z}\in C_{u,z},\label{flow_learning}\\
x_{u,z}^+&\in \hat{G}_{{\delta}}\big(x_{u,z}\big),~~~~~~~x_{u,z}\in D_{u,z},
\end{align}
\end{subequations}\noindent
where the set-valued mappings $\hat{F}_{\delta}:\mathbb{R}^{n+r}\times\mathbb{R}^{n}\rightrightarrows\mathbb{R}^{n+r}$ and $\hat{G}_{\delta}:\mathbb{R}^{n+r}\rightrightarrows\mathbb{R}^{n+r}$ are allowed to depend on a small tunable parameter $\delta>0$ that provides flexibility for the design of the controller. It follows that the data of system \eqref{hybrid_learning_dynamics} is given by $\hat{\mathcal{H}}_{\delta}:=\left\{C_{u,z},\hat{F}_{\delta},D_{u,z},\hat{G}_{\delta}\right\}$, it receives as input the state $\xi$, generated by the filter \eqref{filterdynamics}, and it generates as output the state $\hat{u}$, which acts on the plant via the following control law: 
\begin{equation}\label{main_input}
u=\hat{u}+\varepsilon_a\mathbb{D}\mu.
\end{equation}
%
%
 
% In order to design the data of system $\hat{\mathcal{H}}_{\delta}$ we assume that $\xi:=\nabla J(\hat{u})$, being $\nabla J:\mathbb{R}^n\rightarrow\mathbb{R}^n$ the gradient of the response map $J$. The following regularity and stability assumptions characterize the data (\ref{hybrid_learning_dynamics}) of $\hat{\mathcal{H}}_{\delta}$.
%

The nominal hybrid decision-making dynamics \eqref{hybrid_learning_dynamics} are designed so that the following assumption holds:
\begin{assumption}[Hybrid Decision-Making Dynamics]\label{keyassumptions}
The following properties hold:
\begin{enumerate}[(a)]
\item There exists $\delta^*\in\mathbb{R}_{>0}$ such that for all $\delta\in(0,\delta^*]$:
\begin{enumerate}[1.]
\item the sets $C_{u,z}$ and $D_{u,z}$ are closed, and satisfy $P_u(C_{u,z})\subset\hat{\mathbb{U}}$, $P_u(D_{u,z})\subset\hat{\mathbb{U}}$, where $P_u$ denotes the projection on the $u$-component of the set.
\item $\hat{F}_{\delta}(\cdot,\cdot)$ is OSC, LB, and convex-valued relative to $C_{u,z}\times\mathbb{R}^n$, and $\hat{G}_{\delta}(\cdot)$ is OSC and LB relative to $D_{u,z}$.
\end{enumerate}
\item There exists a compact set $\Psi\subset \mathbb{R}^r$  such that the set $\mathcal{A}:=\mathcal{O}\times\Psi$ is SGPAS as $\delta\rightarrow0^+$ for the \emph{target hybrid decision-making system} corresponding to \eqref{hybrid_learning_dynamics} with $\xi:=\nabla J(\hat{u})$, where the set  $\mathcal{O}$ is the set of solutions to \eqref{optimization_problem}.
%\item For each compact set $K\subset C_{u,z}\cup D_{u,z}$ there exists $\epsilon^*\in\mathbb{R}_{>0}$ such that for each $x_{u,z}(0,0)\in K$ and bounded measurable function $e(\cdot)$ satisfying $\sup_{t\geq0}|e(t)|\leq \epsilon^*$, the HDS (\ref{hybrid_learning_dynamics}) with flow map $\hat{F}_{\delta}(x_{u,z},\nabla J(\hat{u})+e)$ generates at least one complete solution.
%
\end{enumerate}
\end{assumption}

\begin{figure}[t!]
\begin{tcolorbox}[colback=ivoryA, colframe=ivoryA]
\centering
\includegraphics[width=0.92\linewidth]{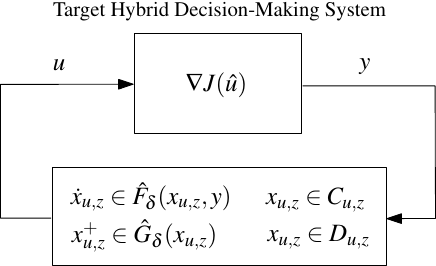}
\end{tcolorbox}
\caption{\tcr{Conceptual feedback scheme of the target hybrid decision-making system considered in Assumption \ref{keyassumptions}.} \label{fig5cd2}}
\end{figure}

While condition (a) provides regularity to the data of system $\hat{\mathcal{H}}_{\delta}$, item (b) indicates that $\hat{\mathcal{H}}_{\delta}$ should be designed to guarantee that the set of solutions to problem \eqref{optimization_problem} is stabilized \emph{under the assumption that the gradient $\nabla J$ is known by the hybrid decision-making dynamics}. Figure \ref{fig5cd2} shows a block diagram representation of the target hybrid decision-making system considered in item (b) of Assumption \ref{keyassumptions}. However, we stress that, in practice, the hybrid controller will use only real-time measurements of $J$, and the exact mathematical forms of $J$ and $\nabla J$ are assumed to be unknown. However, qualitative properties of $J$, e.g., quadratic-like structures, convexity, invexity, etc, are typically used to study stability in ES algorithms. 

If $\hat{\mathcal{H}}_{\delta}$ does not depend on any parameter $\delta$, then item (b) in Assumption \ref{keyassumptions} simply asks that the set $\mathcal{A}$ is UGAS. If, additionally, system \eqref{hybrid_learning_dynamics} does not depend on any auxiliary state $\hat{z}\in\mathbb{R}^r$,  then the set $\Psi$ can be neglected, and we can take $C_{u,z}=C_u\subset\mathbb{R}^n$ and $D_{u,z}=D_u\subset\mathbb{R}^n$. We will use the convention $r=0$ to specify this case. Finally, note that if, in addition, we have that $D_u=\emptyset$, $C_u=\mathbb{R}^n$, and $\hat{F}_{\delta}$ is a Lipschitz continuous function, then the above properties reduce to the standard assumptions considered in the literature of smooth extremum seeking \cite{DerivativesESC,MISONewton,Moase}. Local (practical) stability results can also be established when only local qualitative properties of \(\nabla J\) or \( J \) are known.

%  \begin{figure*}[t!]
% \centerline{\includegraphics[width=0.65\linewidth]{Figure13b.eps}}
%  \caption{(a) \tcr{Conceptual feedback scheme of the target hybrid decision-making dynamics considered in Assumption \ref{keyassumptions}. (b) Block diagram of hybrid set-seeking dynamics with growing timers and no filters.} \label{fig5cd2}}
% \end{figure*}
%
%%%%%%%%%%%%%%%%%%%%%%%%%%%%%%%%%%%%%%%%%
%
\subsubsection{Hybrid-Set Seeking: Closed-Loop System}
The overall hybrid set-seeking dynamics are obtained by interconnecting the static plant $y=J(u)$, the oscillator \eqref{oscillatordynamics}, the filter \eqref{filterdynamics}, and the hybrid decision-making dynamics \eqref{hybrid_learning_dynamics} using the feedback law \eqref{main_input}. The resulting hybrid system has an overall state $x:=(x_{u,z},\xi,\mu)\in\mathbb{R}^{r+4n}$, and data given by $\mathcal{H}:=\{C,F,D,G\}$, with
\begin{subequations}\label{originalsystemcomplete01}
 \begin{align}
 C:&=C_{u,z}\times \Lambda_{\xi}\times\mathbb{S}^n\label{flowsetcomplete}\\
 \dot{x}\in F(x):&= \left(  \begin{array}{c}
\hat{F}_{{\delta}}(x_{u,z},\xi)\\
 -\frac{k_f}{\varepsilon_f}\left(\xi-\frac{2}{\varepsilon_a} y\cdot\mathbb{D}\mu\right)\\
\Phi(\omega)\mu
\end{array}\right)\label{flowmapcomplete}\\
D:&=D_{u,z}\times \Lambda_{\xi}\times\mathbb{S}^n\label{jumpsetcomplete}\\
x^+\in G(x):&=\left(  \begin{array}{c} 
\hat{G}_{{\delta}}(x_{u,z})\\
\xi\\
\mu
  \end{array}\right).\label{jumpmapcomplete}
\end{align}
\end{subequations}\noindent
The stability properties of \eqref{originalsystemcomplete01} are captured by the following theorem, which follows by computing the average dynamics of \eqref{originalsystemcomplete01} with respect to $\mu$, and then applying singular perturbation arguments to study the average hybrid dynamics as a perturbed multi-time scale system with the state $\xi$ evolving on a faster time scale, see Theorem \ref{theorem3SP} and Remark 1 in ``Singularly Perturbed Hybrid Dynamical Systems''.
\begin{theorem}\label{theorem1}
Suppose that Assumption \ref{keyassumptions} holds. Then, the hybrid set-seeking system \eqref{originalsystemcomplete01} with $y=J(u)$ and the feedback control law \eqref{main_input} renders the set $\mathcal{A}\times\Lambda_{\xi}\times\mathbb{S}^n$ SGPAS as $(\delta,\varepsilon_f,\varepsilon_a,\varepsilon_{\omega})\to 0^+$. \hfill \QEDB 
\end{theorem}

The stability result of Theorem \ref{theorem1} states that the trajectories of \eqref{originalsystemcomplete01} satisfy bounds of the form \eqref{KLbound1} with respect to the set $\mathcal{A}\times\Lambda_{\xi}\times\mathbb{S}^n$. In particular, since $\Lambda_{\xi}$ and $\mathbb{S}^n$ are compact, the state $x_{u,z}$ in \eqref{originalsystemcomplete01} satisfies 
\begin{equation}\label{klbound2}
|x_{u,z}(t,j)|_{\mathcal{A}}\leq \beta(|x_{u,z}(0,0)|_{\mathcal{A}},t+j)+\nu,
\end{equation}
for all $(t,j)\in\text{dom}(x_{u,z})$. Note that Theorem \ref{theorem1} does not make assertions about completeness of solutions, since these properties need to be established for each particular application under an appropriate design of the data of the hybrid system. 
\begin{remark}
The SGPAS result of Theorem \ref{theorem1} also establishes a tuning order for the parameters $(\delta,\varepsilon_f,\varepsilon_a,\varepsilon_{\omega})$. The first parameter to be selected "sufficiently" small is $\delta$, such that item (b) in Assumption \ref{keyassumptions} holds. Subsequently, the parameters \( \varepsilon_f \), \( \varepsilon_a \), and \( \varepsilon_{\omega} \) should be selected to be sufficiently small, in that order. \QEDB 
\end{remark}
%
% For each compact set $\tilde{K}$ satisfying $\mathcal{O}\times\Psi\subset\text{int}~(\tilde{K})$ and each $\nu>0$, there exists $\delta^*\in\mathbb{R}_{>0}$ such that for each $\delta\in(0,\delta^*)$ there exists  $\omega_{L}^*\in\mathbb{R}_{>0}$ such that for each $\omega_L\in(0,\omega_L^*)$ there exists $a^*\in\mathbb{R}_{>0}$ such that for each $a\in(0,a^*)$ there exists $\varepsilon_{\omega}^*$, such that for each

% $\epsilon^*$ such that for each $\epsilon\in(0,\epsilon^*)$ there exists a compact set $\mathcal{A}_{\varepsilon}$ satisfying
% % 
% \begin{equation}
% \mathcal{A}_{\varepsilon}\subset \left((\mathcal{O}\times \Psi)+\varepsilon\mathbb{B}\right)\times\lambda_{\xi}\mathbb{B}\times\mathbb{S}^n\times\lambda_{\theta}\mathbb{B},
% \end{equation}
% %
% which is UGAS for the closed-loop system \eqref{originalsystemcomplete} with restricted flow set $\hat{C}\cap \tilde{K}\times \lambda_{\xi}\mathbb{B}\times \mathbb{S}^n\times\lambda_{\theta}\mathbb{B}$ and restricted jump set $\hat{D}\cap \tilde{K}\times \lambda_{\xi}\mathbb{B}\times \mathbb{S}^n\times\lambda_{\theta}\mathbb{B}$.
%
%
%The following corollary is a direct consequence of Theorem \ref{theorem1}.
%
\begin{corollary}\label{Corollary1}
Suppose that Assumption \ref{keyassumptions} holds and the sets $C_{u,z}$ and $D_{u,z}$ are bounded. Then, the hybrid dynamics \eqref{originalsystemcomplete01} with $y=J(u)$ and control law \eqref{main_input} renders the set $\mathcal{A}\times\Lambda_{\xi}\times\mathbb{S}^n$ Globally Practically Asymptotically Stable as $(\delta,\varepsilon_f,\varepsilon_a,\varepsilon_{\omega})\to 0^+$.  \QEDB 
\end{corollary}
Achieving \emph{global} stability results when the flow and jump sets are bounded is particularly relevant for ES problems defined on smooth compact manifolds. In such spaces, standard  ES algorithms with target decision-making dynamics characterized by smooth ODEs can achieve only local, or at best, almost global, asymptotic stability. This follows because such spaces are not contractible, a property that is needed for global stability in smooth dynamical systems \cite{sontag1999stability}. However, hybrid set-seeking systems can overcome these limitations and achieve global bounds of the form \eqref{klbound2}, see \cite{Strizic:17_CDC,OchoaPovedaManifolds}.

\begin{remark}
When $D_{u,z}:=\emptyset$, system \eqref{originalsystemcomplete01} does not experience jumps and the set-seeking dynamics reduce to the continuous-time system
\begin{subequations}\label{originalsystemcomplete03}
 \begin{align}
 C:&=C_{u,z}\times \Lambda_{\xi}\times\mathbb{S}^n,\label{flowsetcomplete02}\\
 \dot{x}\in F(x):&= \left(\begin{array}{c}
\hat{F}_{{\delta}}(x_{u,z},\xi)\\
 -\frac{k_f}{\varepsilon_f}\left(\xi-\frac{2}{\varepsilon_a} y\cdot\mathbb{D}\mu\right)\\
\Phi(\omega)\mu
\end{array}\right).\label{flowmapcomplete02}
\end{align}
\end{subequations}\noindent
Set-valued ES systems of the form \eqref{originalsystemcomplete03} can emerge in applications where the optimizing vector field is discontinuous. In such cases, $\hat{F}_{{\delta}}$ will correspond to the Krasovskii regularization of the vector field (see "Krasvoskii Solutions of ODEs" in "Continuous-Time Set-Valued Dynamical Systems"). If, on the other hand, $\hat{F}_{\delta}$ is a Lipschitz continuous function, and $C_{u,z}=\mathbb{R}^n$, then \eqref{originalsystemcomplete03} recovers the standard framework for the study of gradient-based smooth ES dynamics \cite{KrsticBookESC}, \cite{TanAndNesic2006Local}, \cite{DerivativesESC}. \QEDB
\end{remark}

One of the advantages of modeling seeking systems with autonomous dither generators of the form \eqref{oscillatordynamics} is that the resulting closed-loop system \eqref{originalsystemcomplete01} satisfies the Hybrid Basic Conditions and has suitable stability properties with respect to a compact set (see Theorem \ref{theorem1}). As a consequence, one can study the inflated hybrid system \eqref{eqinflatedHDS} obtained from the data of the hybrid seeking dynamics \eqref{originalsystemcomplete01} to establish the following robustness result via \cite[Thm. 7.21]{bookHDS}.
\begin{corollary}\label{Corollary1}
Suppose that Assumption \ref{keyassumptions} holds. For each compact set of initial conditions $K_0\subset\mathbb{R}^{n+r}$ let Theorem \ref{theorem1} generate sufficiently small parameters $(\delta,\varepsilon_f,\varepsilon_a,\varepsilon_{\omega})$ such that property \eqref{klbound2} holds. Then, the \emph{inflated} hybrid set-seeking system \eqref{eqinflatedHDS} with control law \eqref{main_input} renders the set $\mathcal{A}\times\Lambda_{\xi}\times\mathbb{S}^n$ SGPAS as $\rho\to0^+$ (relative to $K_0$).
\QEDB 
\end{corollary}
Since, as discussed in Part 1 of this article, the inflated system \eqref{eqinflatedHDS} can capture a variety of additive disturbances and perturbations acting on the states and dynamics of the algorithms, the result of Corollary \ref{Corollary1} essentially implies that the hybrid set-seeking algorithms \eqref{originalsystemcomplete01} will retain their (semi-global practical) stability properties under a broad class of small additive disturbances on the states and dynamics.
\begin{remark}\label{remarkunstable}
To account for situations where the plant output is influenced by external events, such as sporadic feedback measurements or external malicious signals, the filter in \eqref{flowmapcomplete} can be modified to adopt the more general structure
\begin{equation}
\dot{\xi}=-\frac{k_f}{\varepsilon_f}\left(\xi-\frac{2}{\varepsilon_a}\gamma(x_{u,z})yD\mu\right),
\end{equation}
where $\gamma(\cdot)$ is a continuous function. In this case, item (b) of Assumption \ref{keyassumptions} should hold with $\xi=\gamma(x_{u,z})\nabla J(\hat{u})$. \QEDB 
\end{remark}
\subsubsection{Hybrid-Set Seeking Dynamics with Hybrid Filters}
\tcr{By removing any explicit dependence of $\hat{F}_{\delta}$ on the oscillatory signal $\mu$, the non-hybrid low-pass filter considered in \eqref{originalsystemcomplete01} facilitates the computation of the average hybrid dynamics for a broad class of set-valued decision-making dynamics \eqref{hybrid_learning_dynamics}. However, this comes at the price of having an additional time scale in the system induced by selecting $\varepsilon_f>0$ sufficiently small, and by restricting $\xi$ to evolve in a compact set for the purpose of stability analysis. Nevertheless, it is possible to relax these assumptions by replacing the flow and jump sets of \eqref{originalsystemcomplete01} by the more general sets
\begin{subequations}\label{hybridfilteres}
\begin{align}
C:=C_{u,z}\times\Lambda_{\xi,c}\times\mathbb{S}^n,~~~~~D:=D_{u,z}\times\Lambda_{\xi,d}\times\mathbb{S}^n, 
\end{align}
where $\Lambda_{\xi,c},\Lambda_{\xi,d}\subset\mathbb{R}^n$ are closed sets, and by considering the following flow and jump maps:
\begin{align}
\dot{x}\in F(x):&= \left(  \begin{array}{c}
\hat{F}_{{\delta}}(x_{u,z},\xi)\\
 -k_f\left(\xi-\frac{2}{\varepsilon_a} y\cdot\mathbb{D}\mu\right)\\
\Phi(\omega)\mu
\end{array}\right)\label{flowmapcomplete09876}\\
x^+\in G(x):&=\left(  \begin{array}{c} 
\hat{G}_{{\delta}}(x_{u,z},\xi)\\
G_{\xi}(x_{u,z},\xi)\\
\mu
\end{array}\right),
\end{align}
\end{subequations}
where the set-valued maps $\hat{G}_{\delta}$ and $G_{\xi}$ are both assumed to be OSC and LB relative to $D_{u,z}\times\Lambda_{\xi,d}$. In this case, we consider the following assumption:}
\tcr{
\begin{assumption}\label{assumption2}
Items (a) of Assumption \ref{keyassumptions} holds, and for the target hybrid decision-making dynamics
\begin{subequations}\label{hybrid_learning_dynamics001}
\begin{align}
\left(\begin{array}{c}
\dot{x}_{u,z}\\
\dot{\xi}
\end{array}\right) &\in \left(\begin{array}{c}
\hat{F}_{{\delta}}\big(x_{u,z},\xi\big)\\
-k_f\left(\xi-\kappa\nabla J(\hat{u})\right)
\end{array}\right),(x_{u,z},\xi)\in C_{u,z}\times \Lambda_{\xi,c}\label{flow_learninghf001}\\
\left(\begin{array}{c}
x_{u,z}^+\\
\xi^+
\end{array}\right)
&\in \left(\begin{array}{c}
\hat{G}_{{\delta}}\big(x_{u,z},\xi\big)\\
G_{\xi}(x_{u,z},\xi)
\end{array}\right),~(x_{u,z},\xi)\in D_{u,z}\times\Lambda_{\xi,d},\label{jump_learninghf001}
\end{align}
\end{subequations}\noindent
where $k_f>0$, there exist $\kappa\in\mathbb{R}$ and compact sets $\Psi_z\subset \mathbb{R}^r$ and $\Psi_f\subset\mathbb{R}^n$ such that the set $\mathcal{A}:=\mathcal{O}\times\Psi_z\times\Psi_f$ is SGPAS as $\delta\rightarrow0^+$.
\end{assumption}
With Assumption \ref{assumption2} at hand, we can establish the following theorem.
\begin{theorem}\label{theorem4}
Suppose that Assumption \ref{assumption2} holds. Then, the hybrid set-seeking system \eqref{hybridfilteres} with $y=J(u)$ and control law \eqref{main_input} renders the set $\mathcal{A}\times\mathbb{S}^n$ SGPAS as $(\delta,\varepsilon_a,\varepsilon_{\omega})\to 0^+$.  \QEDB 
\end{theorem}
}
\tcr{Seeking systems with hybrid filters are useful for incorporating resets or restarting mechanisms to reduce overshoot, improve transient performance \cite{teel2019first,Candes_Restarting,PovedaNaliAuto20}, or regularize accelerated gradient flows with dynamic damping, which can speed up convergence when the cost function exhibits low curvature \cite{ODE_Nesterov,galarza2021self,PovedaNaliAuto20}. In such systems, the continuous-time dynamics of \( \xi \) can be appropriately modified without significantly altering the structure of \eqref{hybrid_learning_dynamics001}.}

\subsubsection{Hybrid Set-Seeking Dynamics with Growing Timers}
\tcr{We now consider the simpler hybrid set-seeking systems that dispense with the filters and for which the time variation in the dither signal is generated by an exogenous timer that grows unbounded, as in system \eqref{ESCODEVanilla} and \eqref{eq:HDSmodel2}, see also  ``Averaging Theory: From ODEs to Hybrid Systems''. In this case, we can consider hybrid set-seeking systems with control law \eqref{main_input} and dynamics
\begin{subequations}\label{timergrowsunbounded}
\begin{align}
(x_{u,z},\tau)&\in C_{u,z}\times\mathbb{R}_{\geq0},~~\left\{\begin{array}{l}
\dot{x}_{u,z}=\hat{F}_{\delta}(x_{u,z},y\cdot \mu_{\varepsilon_a}(\tau))\\
~~~\dot{\tau}=\frac{1}{\varepsilon_{\omega}}
\end{array}\right.\label{flowmapgrowing}\\
(x_{u,z},\tau)&\in D_{u,z}\times\mathbb{R}_{\geq0},~~
\left\{\begin{array}{l}
x^+_{u,z}=\hat{G}_{\delta}(x_{u,z})\\
~\tau^+=\tau
\end{array}\right.,
\end{align}
\end{subequations}
where $\mu_{\varepsilon_a}(\tau):=\frac{2}{\varepsilon_a}(\sin(2\pi\kappa_1\tau),\sin(2\pi\kappa_2\tau),\ldots,\sin(2\pi\kappa_n\tau))$, and the parameters $\kappa_i>0$ are rational and satisfy $\kappa_i\neq\kappa_j$ for all $i\neq j\in\{1,2,\ldots,n\}$. Note that the ``bouncing'' seeking system studied in Example \ref{seekingballexample1} has the form of \eqref{timergrowsunbounded}. When $D_{u,z}=\emptyset$, $C_{u,z}=\mathbb{R}^n$, and $\hat{F}_{\delta}=-y\mu_{\varepsilon_a}(\tau)$ we recover the most basic continuous-time ES scheme studied in the literature  \cite{KrsticBookESC,TanAndNesic2006Local}.}

\tcr{The stability properties of the hybrid set-system \eqref{timergrowsunbounded} are studied under the following assumption, which leverages the notion of \emph{average map} introduced in ``Averaging Theory: From ODEs to Hybrid Systems''.}

\begin{assumption}\label{assumptionseekinggrowing}
\tcr{There exists $\delta^*\in\mathbb{R}_{>0}$ such that for all $\delta\in(0,\delta^*]$ the following properties hold:} 
\begin{enumerate}[(a)]
\item \tcr{Item (a)-1. of Assumption \ref{keyassumptions} is satisfied, and the functions $\hat{F}_{\delta}$ and $\hat{G}_{\delta}$ are continuous.} 

\item \tcr{There exists a compact set $\Psi\subset \mathbb{R}^r$  such that the set $\mathcal{A}:=\mathcal{O}\times\Psi$ is SGPAS as $\delta\rightarrow0^+$ for the target hybrid decision-making dynamics}
\begin{subequations}\label{nominaldynamicsgrowingtimer}
\begin{align}
\dot{x}_{u,z}&=\hat{F}_{\delta}(x_{u,z},\nabla J(\hat{u})),~~~x_{u,z}\in C_{u,z}\\
x_{u,z}^+&=\hat{G}_{\delta}(x_{u,z}),~~~~~~~~~~~x_{u,x}\in D_{u,z}.
\end{align}
\end{subequations}
\item \tcr{For each compact set $K\subset\mathbb{R}^n$ there exists $\varepsilon_{a}^*>0$ such that for all $\varepsilon_a\in(0,\varepsilon_a^*)$ the average map of $\hat{F}_{\delta}$, denoted $\overline{F}_{\delta}$, satisfies 
$\bar{F}_{\delta}(\bar{x}_{u,z})\in F_{a}(\bar{x}_{u,z})$ for all $\bar{x}_{u,z}\in K$, where $F_{a}$ is the inflated set-valued map constructed as in \eqref{inflatedflowmapconstruction} from the data $F(x_{u,z}):=\hat{F}_{\delta}(x_{u,z},\nabla J(\hat{u}))$.}
\end{enumerate}
\end{assumption}
\tcr{The conditions of Assumption 3 essentially ask that the flow map in \eqref{timergrowsunbounded} is designed such that, as $\varepsilon_{\omega}\to0^+$, and on average, the hybrid dynamics behave as an $\varepsilon_a$-perturbed hybrid system for which the unperturbed nominal dynamics \eqref{nominaldynamicsgrowingtimer} are able to solve the decision-making problem \eqref{optimization_problem}.}
%

%~~~~~\includegraphics[width=0.59\linewidth]{schemeUpdated2.jpg}

\tcr{A direct application of averaging tools---see Theorems 1-2 in \emph{Averaging Theory: From ODEs to Hybrid Systems}---along with robustness results for hybrid systems satisfying the Hybrid Basic Conditions (see \cite[Thm.~7.21]{bookHDS}), leads to the following theorem concerning the closed set \( \mathcal{A} \times \mathbb{R}_{\geq 0} \):}

\begin{theorem}
\tcr{Suppose that Assumption \ref{assumptionseekinggrowing} holds. Then, the hybrid set-seeking system \eqref{timergrowsunbounded} with $y=J(u)$ and control law $$u=\hat{u}+\frac{\varepsilon_a^2}{2}\mu_{\varepsilon_a}(\tau),$$renders the set $\mathcal{A}\times\mathbb{R}_{\geq0}$ SGPAS as $(\delta,\varepsilon_a,\varepsilon_{\omega})\to 0^+$.} \QEDB
\end{theorem}
\subsection{Hybrid Set-Seeking Systems with Estimates of Higher Derivatives}
\tcr{A variety of decision-making dynamics benefit from incorporating higher-order derivatives of the cost function to improve transient performance. For example, Newton methods utilize the second derivative of \( J \) to eliminate convergence dependence on the cost function's curvature. In smooth ES systems, it has been shown \cite{NesicChina,DerivativesESC,MISONewton} that higher-order derivatives can also be embedded into the controller through suitable design and manipulation of the dither signals \( \mu \) and the output function \( y \). Naturally, similar extensions are possible in the context of hybrid set-seeking systems by appropriately modifying the flow maps in \eqref{flowmapcomplete} or \eqref{flowmapgrowing}. While we do not explore these manipulations in detail, we refer the reader to \cite{NesicChina,DerivativesESC,MISONewton,GalarzaPovedaDallanese,PovedaKrsticWC20}, as well as to an example in the next section, where an ES algorithm switches between Newton and gradient-based dynamics to enhance performance.}

%
%
%%%%%%%%%%%%%%%%%%%%%%%%%%%%%%%%%%%
\section{Examples and Applications of Hybrid Set-Seeking Systems}
\label{Sec:Examples}
 The study of adaptive seeking systems with logic modes, event-based rules, and set-valued dynamics opens new avenues for research at the intersection of real-time, model-free control, and computer science, with applications in cyber-physical systems. In particular, equipped with the previous theoretical tools, in this section, we now study different types of hybrid set-seeking algorithms, including systems that can be modeled as logic-based switching ES, state-based switched ES, and reset-based ES. By working with well-posed hybrid systems, we show how different "\emph{if-then}" rules can be systematically incorporated into set-seeking algorithms, while providing stability and robustness guarantees that parallel those existing for continuous-time systems.
\subsection{Constrained Seeking via Set-Valued Dynamics}
By considering \eqref{originalsystemcomplete01} with  $x_{uz}=\hat{u}$, $r=0$, $\xi=\nabla J$, and sets $C_{u,z}=\hat{\mathbb{U}}$, $D_{u,z}=\emptyset$, $C_z=D_z=\Psi=\hat{G}_{\delta}=\emptyset$, we can apply Theorem~\ref{theorem1} to establish stability properties for a broad class of novel continuous-time set-valued set-seeking systems, where the target decision-making dynamics evolve according to
\begin{equation}\label{setvaluedESexample}
\dot{\hat{u}}\in f\left(\hat{u},\nabla J(\hat{u})\right),~~~\hat{u}\in \hat{\mathbb{U}}.
\end{equation}
In particular, set-valued algorithms of the form \eqref{setvaluedESexample} are common in decision-making problems with constraints due to the use of inner parametric optimization loops \cite{Sandholm1, Clarke,Clarke2}, projections \cite{Projections1}, or state-dependent switching rules \cite{galarza2022sliding}. 

The example below considers a set-valued seeking algorithm inspired by population dynamics \cite{Sandholm1}, and also illustrates the discussion on Remark \ref{remarkMAS}.

\begin{figure}[t!]
\begin{tcolorbox}[colback=iceblue!80, colframe=iceblue!80]
\centering
\includegraphics[width=0.99\linewidth]{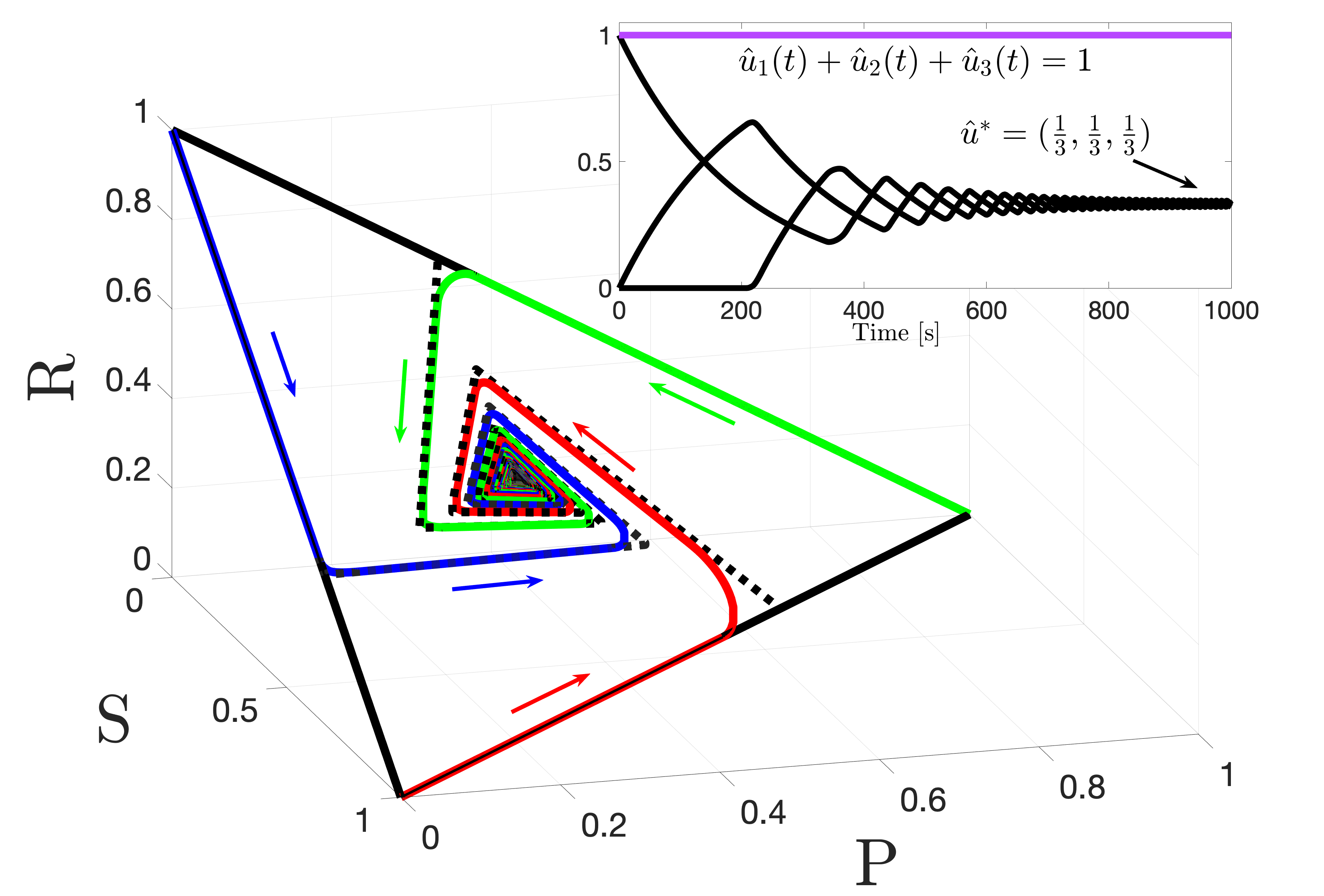}
\includegraphics[width=0.99\linewidth]{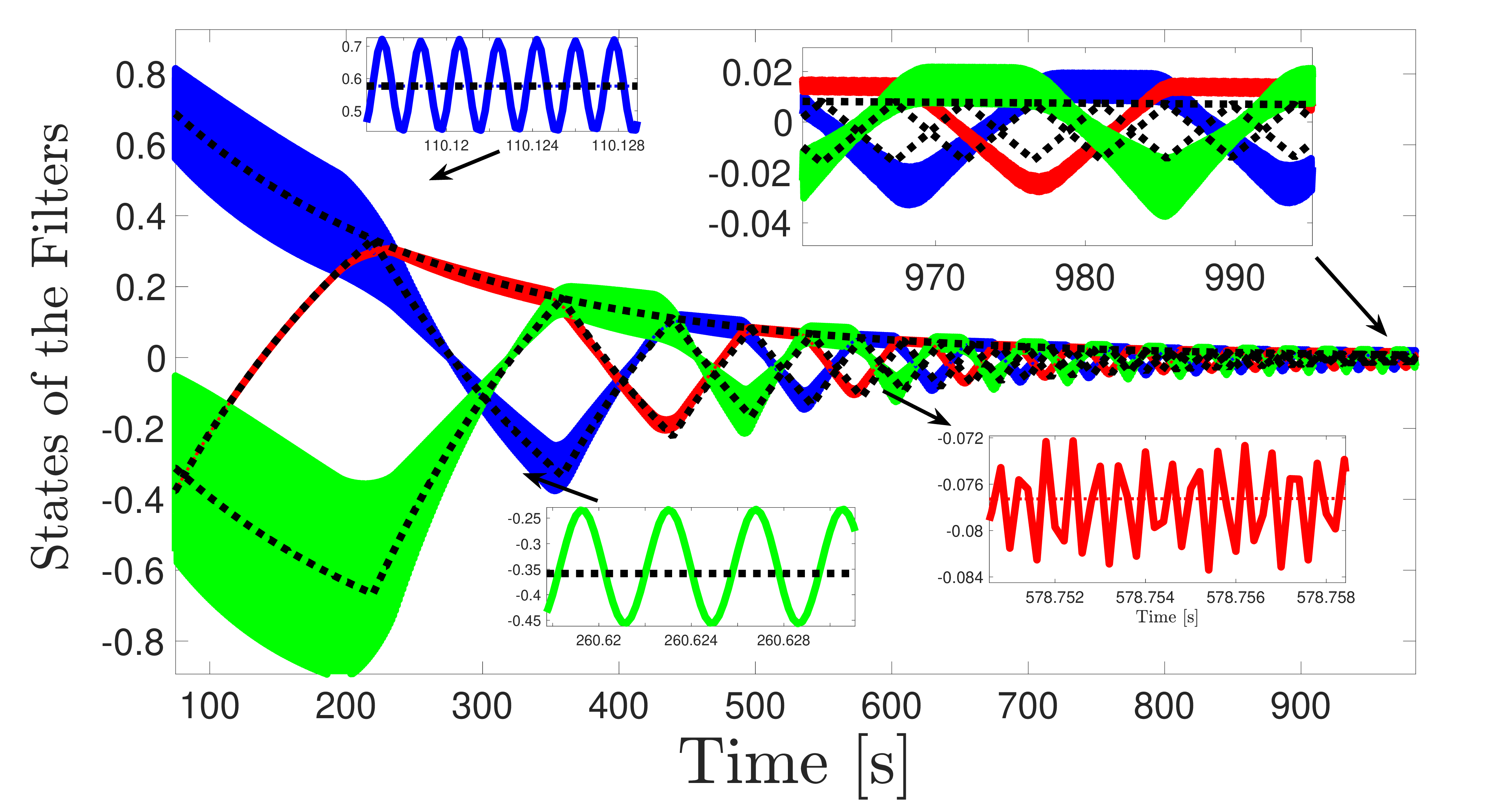}
\end{tcolorbox}
\caption{\tcb{Trajectories of the set-seeking system studied in 
Example~\ref{example10Pop}. \textbf{Top:} Three simulations with 
different initial conditions $\hat{u}(0)$, evolving on the simplex 
$\hat{\mathbb{U}}$ and converging to a neighborhood of $u^*$. 
\textbf{Bottom:} Time evolution of the filter states for one particular 
initialization of the algorithm.} \label{fig100098721345}}
\end{figure}

\tcs{\begin{example}[Nash-Seeking in Population Games]\label{example10Pop}
Consider a multi-agent system with three agents having payoffs of the form
\begin{equation}
J_1=u_1(-u_3+u_2),~
J_2=u_2(-u_3+u_1),~
J_3=u_3(-u_2+u_1),
\end{equation}
and actions $u=(u_1,u_2,u_3)\in\mathbb{R}^3$ restricted to evolve in the simplex set $\hat{\mathbb{U}}:=\{u\in\mathbb{R}_{\geq0}^3:u_1+u_2+u_3=1\}$, which is compact and convex. The agents seek to optimize their own payoff $J_i$ by controlling their own action $u_i$, subject to the coupled constraint imposed by the set $\hat{\mathbb{U}}$, resulting in a generalized Nash equilibrium (NE) seeking problem. Since the vector of partial gradients $\nabla J=(\nabla_1 J_1, \nabla_2 J_2, \nabla_3 J_3)$ has entries given by $\nabla_1 J_1=-u_3+u_2$, $\nabla_2 J_2=-u_3+u_1$, and $\nabla_3 J_3=-u_2+u_1$, the Nash-seeking problem can also be studied as a Rock-Paper-Scissors population game with three strategies \cite{Sandholm1}. The unique NE of this game is given by $u^*=(\frac{1}{3},\frac{1}{3}, \frac{1}{3})$. To converge to $u^*$ using only measurements of $J_i$, players can implement \eqref{originalsystemcomplete01} with dynamics
\begin{equation}\label{BRES1}
\dot{\hat{u}}\in -k\hat{u}+k\left\{w\in \hat{\mathbb{U}}:w^\top\xi=\max_{\hat{u}\in\hat{\mathbb{U}}}\hat{u}^\top\xi\right\},~~~\hat{u}\in\hat{\mathbb{U}},
\end{equation}
where $k>0$. By construction, the set-valued dynamics \eqref{BRES1} satisfy item (a) of Assumption \ref{keyassumptions}. Moreover, item (b) of Assumption \ref{keyassumptions} can be verified by studying the target system $\dot{\hat{u}}\in -\hat{u}+\{w\in \hat{\mathbb{U}}:w^\top\nabla J=\max_{\hat{u}\in\hat{\mathbb{U}}}\hat{u}^\top\nabla J\},~~\hat{u}\in\hat{\mathbb{U}}$, which renders UGAS the NE $u^*$. Indeed, this stability property holds for a broader class of population games that satisfy $(u_a-u_b)^\top (\nabla J(u_a)-\nabla J(u_b))\leq 0$ for all $u_a,u_b\in\hat{\mathbb{U}}$ \cite[Thm. 7.2.7]{Sandholm1}, which allows us to verify Assumption \ref{keyassumptions} without the explicit knowledge of the costs $J_i$. Figure \ref{fig100098721345} shows three different trajectories generated by the set-seeking dynamics \eqref{originalsystemcomplete01} using \eqref{BRES1}, evolving on the set $\hat{\mathbb{U}}$, and converging to a neighborhood of the NE $u^*$. Note that the slow-frequency oscillations shown in the inset upper plot are not induced by the proving signals $\mu$, but instead by the discontinuous nature of the dynamics \eqref{BRES1}. On the other hand, the oscillatory behavior on the state $\xi$, induced by the probing signals, is clearly shown on the lower plot. For this simulation, the parameters were selected as $\varepsilon_a=0.01$, $k=0.005$, $k_f/\varepsilon_f=10$, and $\omega=(18.5,20,25.3)\times 10^4$.
\end{example}}

\tcs{Set-valued dynamical systems are also common in the analysis of projected gradient-based algorithms. The following example illustrates this application in a seeking algorithm that satisfies Assumption \ref{assumption2}, where the filter and the decision-making dynamics are designed to operate on the same time scale.}

\tcs{\begin{example}[Set-Seeking with Projected Dynamics]\label{projectedexample}
Consider a model-free optimization problem of the form \eqref{optimization_problem}, where $J$ is strictly convex and $\hat{\mathbb{U}}$ is a closed and convex set with non-empty interior, describing safety constraints that need to be satisfied for all time. To tackle this problem, we can consider the set-seeking system \eqref{originalsystemcomplete01} with $D=\emptyset$ and
\begin{equation}\label{krasovskiiESE}
\dot{\hat{u}}=k\mathcal{P}_{T_{\hat{\mathbb{U}}}(\hat{u})}(-\xi),~~~\hat{u}\in\hat{\mathbb{U}},
\end{equation}
where $k>0$ and $\mathcal{P}_{T_{\hat{\mathbb{U}}}}(\cdot)$ projects the vector $-\xi$ onto the tangent cone of $\hat{\mathbb{U}}$ at $\hat{u}$. In general, this projection leads to discontinuous vector fields,
necessitating tools from nonsmooth analysis (see ``Continuous-Time Set-Valued Dynamical Systems''); namely, the use of Krasovskii regularizations \eqref{krasovskiicons}. Indeed, under mild regularity assumptions on $J$ and $\hat{\mathbb{U}}$, see \cite{chen2025continuous}, the Krasovskii solutions to \eqref{originalsystemcomplete01} and \eqref{krasovskiiESE} coincide with their standard (Caratheodory) solutions. Therefore, the resulting target average system \eqref{hybrid_learning_dynamics001}, given by
\begin{equation}
\dot{\hat{u}}=\mathcal{P}_{T_{\hat{\mathbb{U}}}(\hat{u})}(-\xi),~~~\dot{\xi}=-k_f\left(\xi-\nabla J(\hat{u})\right),~\hat{u}\in\hat{\mathbb{U}},
\end{equation}
renders the set $\mathcal{A}=\{u^*\}\times\{ \nabla J(u^*))\}\in\hat{\mathbb{U}}\times\mathbb{R}^n$ UGAS for $k_f$ sufficiently small and $\kappa=1$, see \cite[Lemma 10]{chen2025continuous}. It follows that Assumption \ref{assumption2} is satisfied with $\Lambda_{\xi,c}=\Lambda_{\xi,d}=\mathbb{R}^n$, and $k_f,\kappa>0$ sufficiently small \cite[Lemma 10]{chen2025continuous}. Therefore, by Theorem \ref{theorem4}, the model-free set-seeking dynamics \eqref{originalsystemcomplete01} assure convergence to a neighborhood of the optimal solution $u^*$ in $\hat{\mathbb{U}}$. Figure \ref{figprojected1} shows a trajectory generated by \eqref{originalsystemcomplete01} using the dynamics \eqref{krasovskiiESE}, approximately tracking the time-varying minimizer of a slowly-varying cost function of the form $J(u)=(u_1-q_1)^2+(u_2-q_2)^2$, where $t\mapsto q(t):=(q_1(t),q_2(t))$ is a slowly varying parameter, and where the feasible set $\hat{\mathbb{U}}$ corresponds to a 45$^\circ$-rotated box centered at the origin. As observed, the trajectories $\hat{u}$ generated by the set-seeking dynamics remain in $\hat{\mathbb{U}}$ for all time. The simulation used the parameters $\varepsilon_a=0.01$, $\omega=(100,150)$, $k_f=1$, $k=0.2$. The state $q$ was generated by a linear oscillator of the form \eqref{dynamicoscillator} with frequency equal to $0.001$. 
\begin{figure}[t!]
\begin{tcolorbox}[colback=iceblue!80, colframe=iceblue!80]
\centering
\includegraphics[width=0.9\linewidth]{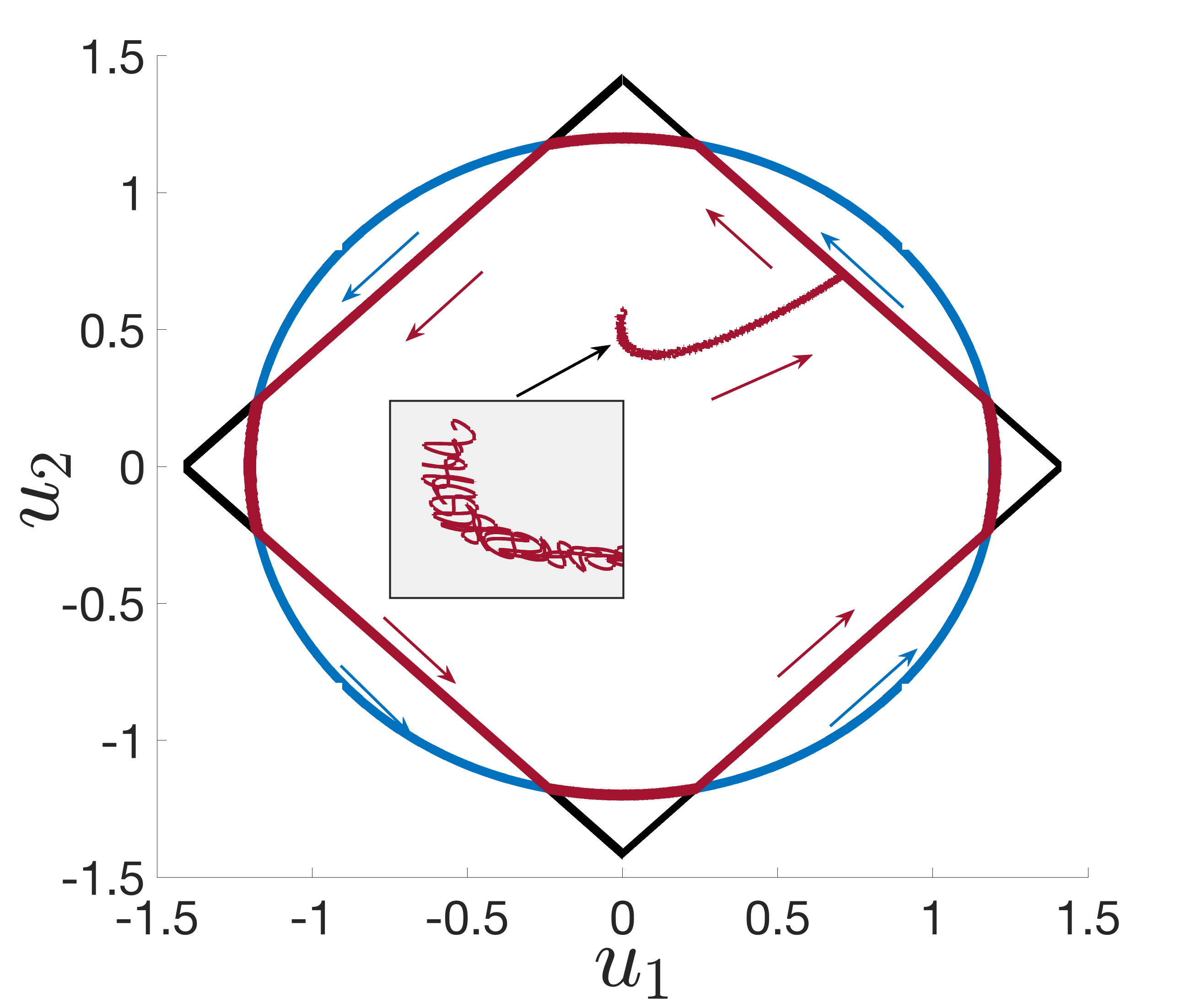}
\includegraphics[width=0.9\linewidth]{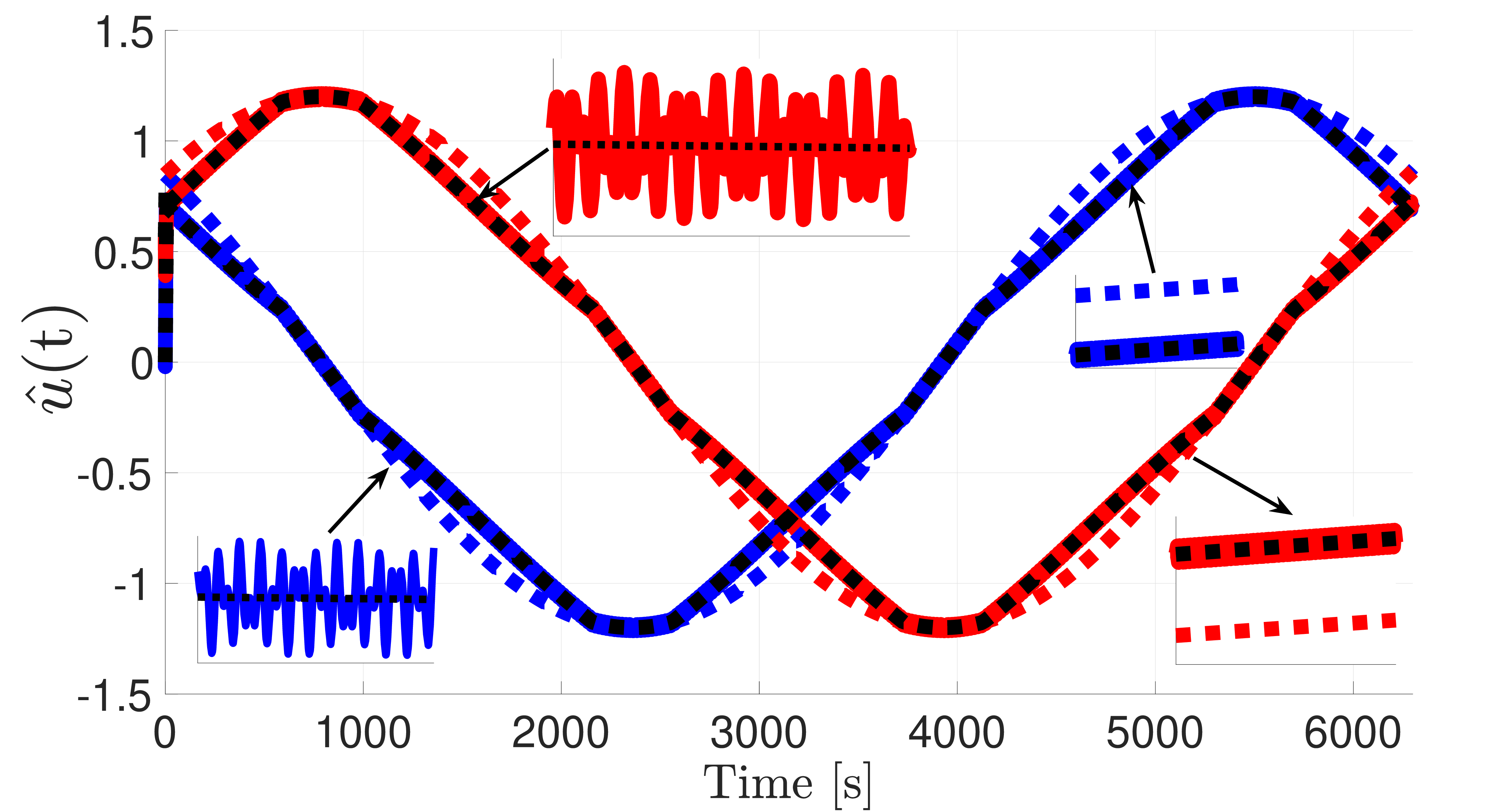}
\end{tcolorbox}
\caption{\tcb{Trajectories of the set-seeking system studied in 
Example~\ref{projectedexample}. \textbf{Top:} Evolution of the 
trajectories $\hat{u}$ (red) within the feasible set $\hat{\mathbb{U}}$, 
corresponding to the rotated square. The blue trajectory represents the 
slowly varying minimizer of the cost function $J$ in $\mathbb{R}^n$. 
As shown, the red trajectories closely follow the blue trajectory inside 
the feasible set. \textbf{Bottom:} Trajectories of $\hat{u}_1$ (blue) and 
$\hat{u}_2$ (red), the optimal solution $u^*$ (black dotted), and the 
time-varying maximizers in $\mathbb{R}^2$ (dotted red and blue). As shown in inset, the trajectory $\hat{u}$ closely tracks the constrained optimal solution $u^*$.} \label{figprojected1}}
\end{figure}
We note that seeking systems with projected dynamics can also be studied using Lipschitz vector fields \cite[Sec. IV.A]{chen2025continuous}, as well as smooth Riemannian projections for decision-making problems defined on manifolds \cite{Poveda:15,OchoaPovedaManifolds}. 
\end{example}}
\begin{example}[Seeking with Unknown Constraints]\label{examplesliding}
\tcs{In Examples \ref{example10Pop} and \ref{projectedexample}, the constraints that describe the decision-making problem \eqref{optimization_problem} were known \emph{a priori} by the practitioner. However, in certain applications, the precise mathematical form of the function that describes the constraints is also unknown. In particular, consider the constrained ES problem \eqref{optimization_problem} where $J$ is strictly convex, and
\begin{equation}
\hat{\mathbb{U}}=\{u\in\mathbb{R}^n:c(u)\leq0\},
\end{equation}
where $c:\mathbb{R}^n\to\mathbb{R}$ is an unknown smooth convex function that is accessible only via evaluations, and such that $\hat{\mathbb{U}}$ has non-empty interior. To tackle the seeking problem with unknown constraints we can consider the dynamics \eqref{originalsystemcomplete01} with $D=\emptyset$ and the set-valued rule 
\begin{equation}\label{slidingseekingcontrol}
\dot{\hat{u}}\in k\cdot\left\{\begin{array}{ll}
\{-\xi_J\}&\text{if}~c(\hat{u})<0\\
K(\xi_J,\xi_c)&\text{if}~c(\hat{u})=0\\
\{-\xi_c\}&\text{if}~c(\hat{u})>0,
\end{array}\right.
\end{equation}
where $k>0$, $K(\xi_J,\xi_c,):=\{z\in\mathbb{R}^n:z=\lambda\xi_J+(1-\lambda)\xi_c,\lambda\in[0,1]\}$, and where the states $\xi_J$ and $\xi_c$ are generated by low-pass filters of the form~\eqref{filterdynamics}, with inputs $y=J(u)$ and $y=c(u)$, respectively. It follows that item (a) of Assumption \ref{keyassumptions} holds \cite{galarza2022sliding}, and the target decision-making system \eqref{hybrid_learning_dynamics} associated with \eqref{slidingseekingcontrol} is given by:
\begin{equation}
\dot{\hat{u}}\in k\cdot \left\{\begin{array}{ll}
\{-\nabla J(\hat{u})\}&\text{if}~c(\hat{u})<0\\
K(\hat{u})&\text{if}~c(\hat{u})=0\\
\{-\nabla c(\hat{u})\}&\text{if}~c(\hat{u})>0
\end{array}\right.
\end{equation}
where $K(\hat{u}):=\{z\in\mathbb{R}^n:z=\lambda\nabla J(\hat{u})+(1-\lambda)\nabla c(\hat{u}),\lambda\in[0,1]\}$, which also satisfies item (b) of Assumption \ref{keyassumptions} \cite{galarza2022sliding}. Therefore, by Theorem \ref{theorem1}, the trajectories of \eqref{originalsystemcomplete01} with dynamics \eqref{slidingseekingcontrol} converge to a neighborhood of $u^*$, the optimizer of $J$ in $\hat{\mathbb{U}}$. Figure~\ref{fig5sliding} depicts a trajectory of the seeking dynamics
generated by~\eqref{slidingseekingcontrol}. The trajectory is initialized
outside the feasible set (shaded region) and tracks the minimizer of a
slowly varying cost function $J(u)$. In this case, $c(u)=a^\top u+b$,
where $a=[-1 -2]$, $b=2$, and $J$ is the same as in
Example~\ref{projectedexample}. The parameters used in the simulation were $\varepsilon_a=0.01$, $\omega=(100,250)$, $k_f=1$, and $k=0.1$. As shown, the trajectory exhibits two phases. Initially, it converges to the
feasible set while disregarding feedback from the cost function. Once inside
the feasible set, the seeking dynamics drive it to a neighborhood of the
optimal solution $u^*$. As is typical, the discontinuous nature of
\eqref{slidingseekingcontrol} induces chattering along the boundary of
the feasible set $\hat{\mathbb{U}}$. In the next sections, we show how to
avoid this chattering by introducing time regularization via resetting
timers, or spatial regularization via switching conditions that induce
hysteresis.}

\end{example}

\begin{figure}[t!]
\begin{tcolorbox}[colback=iceblue!80, colframe=iceblue!80]
\includegraphics[width=0.99\linewidth]{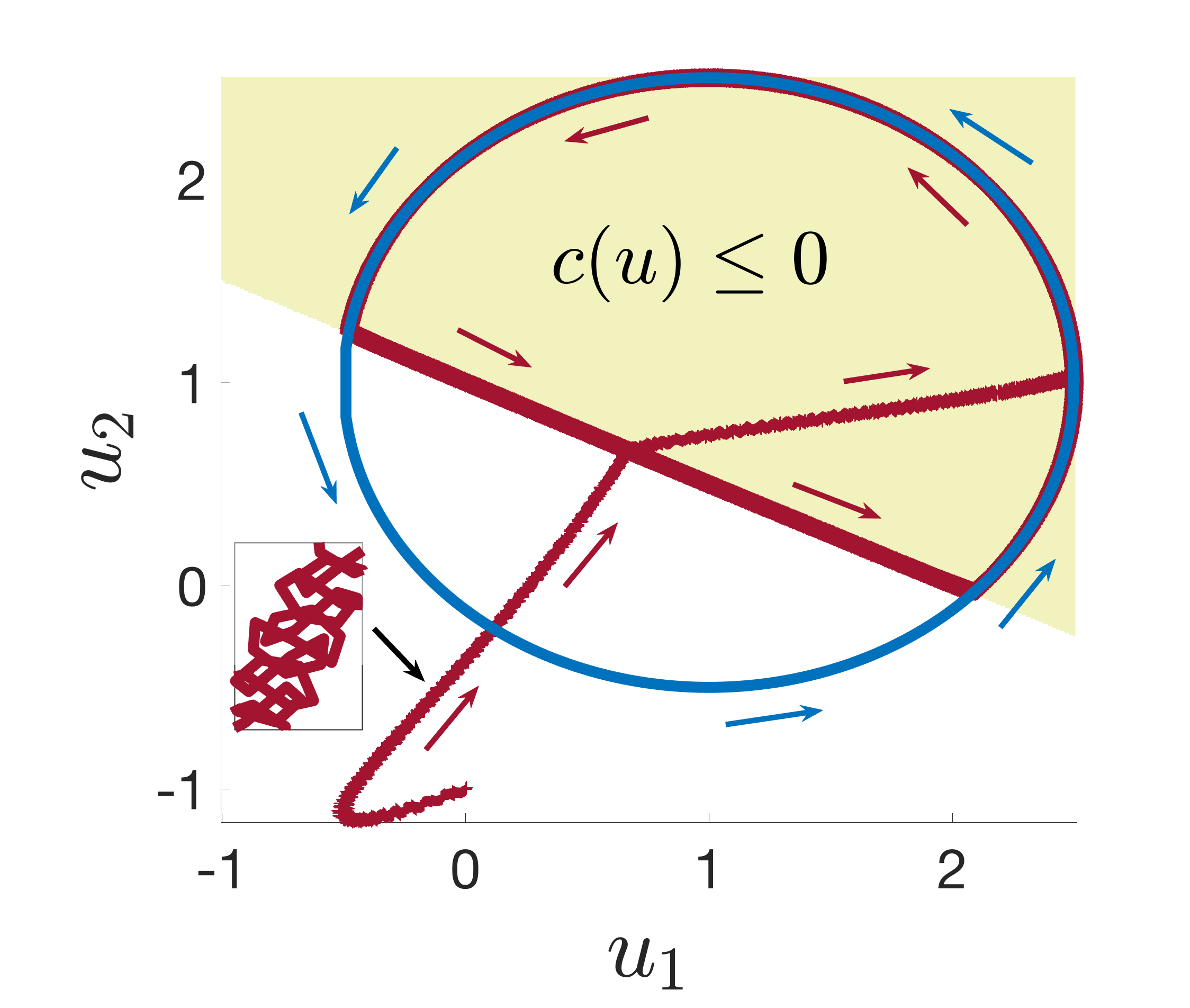}
\end{tcolorbox}
 \caption{Trajectories $\hat{u}$ of the seeking dynamics in Example~\ref{examplesliding}. 
After converging to the feasible set, the trajectories track the slowly varying 
optimizer $u^*$. \label{fig5sliding}}
\end{figure}

\subsection{Seeking with Arbitrarily Fast Switching Dynamics}
\begin{figure}[b!]
\begin{tcolorbox}[colback=ivoryA, colframe=ivoryA]
\centering
\includegraphics[width=0.85\linewidth]{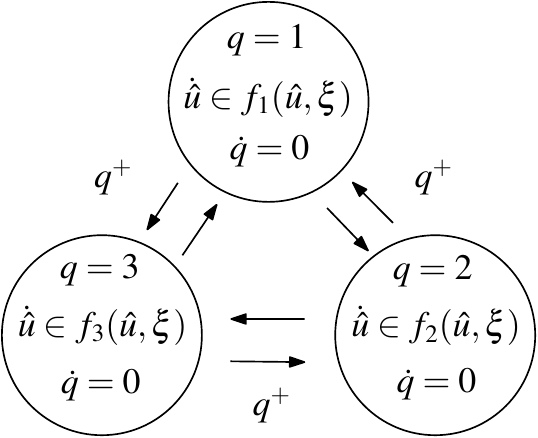}
\end{tcolorbox}
 \caption{Logic-based ES allows to switch between different algorithms or ``operating modes''. \label{fig5cde}}
\end{figure}

Let $Q:=\{1,\ldots,q_{\max}\}$, where $q_{\max}\in\mathbb{Z}_{\geq1}$ and consider a collection of set-valued mappings $f_q:\mathbb{R}^n\times\mathbb{R}^n\rightrightarrows\mathbb{R}^n$, with main state $\hat{u}\in\mathbb{R}^n$, auxiliary state $q\in\mathcal{Q}$, and target decision-making dynamics \eqref{hybrid_learning_dynamics} with $\xi=\nabla J$, given by
\begin{subequations}\label{arbitrary_switching}
\begin{align}
&\dot{\hat{u}}\in {f_q}\left(\hat{u},\nabla J(\hat{u})\right),~~\dot{q}=0,~~(\hat{u},q)\in \hat{\mathbb{U}}\times \mathcal{Q}\label{flow_optimizer0a1}\\
&\hat{u}^+=\hat{u},~~~~~~~~~~~q^+\in \mathcal{Q},~(\hat{u},q)\in \hat{\mathbb{U}}\times \mathcal{Q}\backslash\{q\}.\label{jump_q0a1}
\end{align}
\end{subequations}\noindent
The update rule \eqref{flow_optimizer0a1} describes a gradient-based differential inclusion parameterized by a constant $q$. On the other hand, the dynamics \eqref{jump_q0a1} model the switches of $q$ to a different mode in the set $\mathcal{Q}$, while keeping the main state $\hat{u}$ constant. 
Figure \ref{fig5cde} provides a common graphical representation of logic-based systems that switch between different modes. Such types of representation are connected to the well-known hybrid automata, which can also be written as a HDS \eqref{eq:HDSmodel}, see \cite[Ex. 1.10]{bookHDS}. In this sense, the tools presented in this article can also be used for the synthesis and analysis of seeking systems based on hybrid automata. Nonsmooth systems of the form \eqref{flow_optimizer0a1} can also be considered for the study of ES algorithms with finite-time or fixed-time (practical) convergence properties \cite{poveda2021fixed,poveda2020non}.

Let the data of \eqref{arbitrary_switching} satisfy the conditions of Assumption \ref{keyassumptions}-(a), so that it can be used as a suitable target hybrid decision-making system. The following proposition, adapted from \cite{PoTe17Auto,kamalapurkar2018invariance}, establishes that item (b) in Assumption \ref{keyassumptions} also holds  under arbitrarily fast switching of $q$, provided that there exists a common Lyapunov function that certifies the stability of \eqref{flow_optimizer0a1} with respect to the set $\mathcal{O}$ across all modes $q\in\mathcal{Q}$. 
%.
%
\begin{proposition}[Seeking under Fast Switching]\label{proposition_arbitrary}
Suppose that for each $q\in \mathcal{Q}$ the mapping $f_q$ in (\ref{flow_optimizer0a1}) is OSC, LB, and convex-valued. Suppose also that there exists $\alpha_1, \alpha_2 \in \mathcal{K}_{\infty}$,  a continuous positive definite function $W:\mathbb{R}_{\geq0}\rightarrow\mathbb{R}_{\geq0}$, and a continuously differentiable function $V:\text{dom}(V)\rightarrow\mathbb{R}_{\geq0}$, with $\hat{\mathbb{U}}\subset\text{dom}(V)$ such that for all $\hat{u}\in \hat{\mathbb{U}}$, all $\tilde{f}\in f_q$, and all $q\in \mathcal{Q}$:
\begin{subequations}
\begin{align}
&\alpha_1(|\hat{u}|_{\mathcal{O}})\leq V(\hat{u})\leq \alpha_2(|\hat{u}|_\mathcal{O}) \label{condition1}\\
&\langle\nabla V(\hat{u}),\tilde{f} \rangle\leq -W(|\hat{u}|_{\mathcal{O}}).
\end{align}
\end{subequations}
Then, the differential inclusion
\begin{equation}\label{arbitary_switching}
\dot{\hat{u}}\in \hat{F}(\hat{u},\nabla J(\hat{u})):= \overline{\text{co}}\bigcup_{q\in \mathcal{Q}}f_q(\hat{u},\nabla J(\hat{u})),~~~\hat{u}\in \hat{\mathbb{U}},
\end{equation}\noindent
has the structure of (\ref{hybrid_learning_dynamics}) with the state $x_{uz}=\hat{u}$, $r=0$, $\xi=\nabla J$, and sets $C_{u,z}=\hat{\mathbb{U}}$, $D_{u,z}=\emptyset$, $C_z=D_z=\Psi=\hat{G}_{\delta}=\emptyset$. Moreover, system \eqref{arbitary_switching} satisfies Assumption \ref{keyassumptions}, and for each non-purely discrete solution ($\hat{u}$, $q$) of the switching system \eqref{arbitrary_switching}, $\hat{u}$ is also a solution of \eqref{arbitary_switching}. \QEDB 
\end{proposition}

\begin{figure}
\begin{tcolorbox}[colback=iceblue!80, colframe=iceblue!80]
\centering
\includegraphics[width=0.85\linewidth]{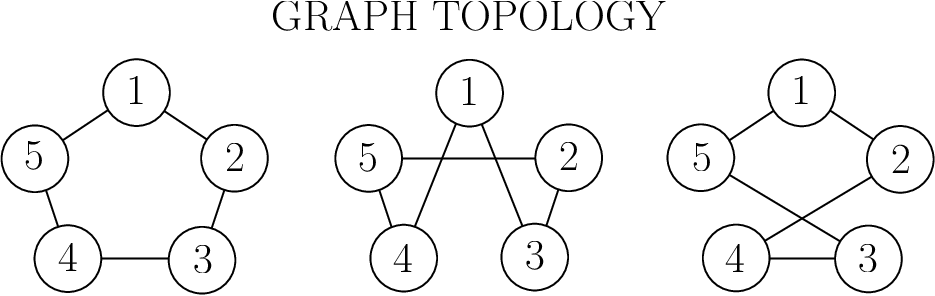}\vspace{0.25cm}
\includegraphics[width=0.95\linewidth]{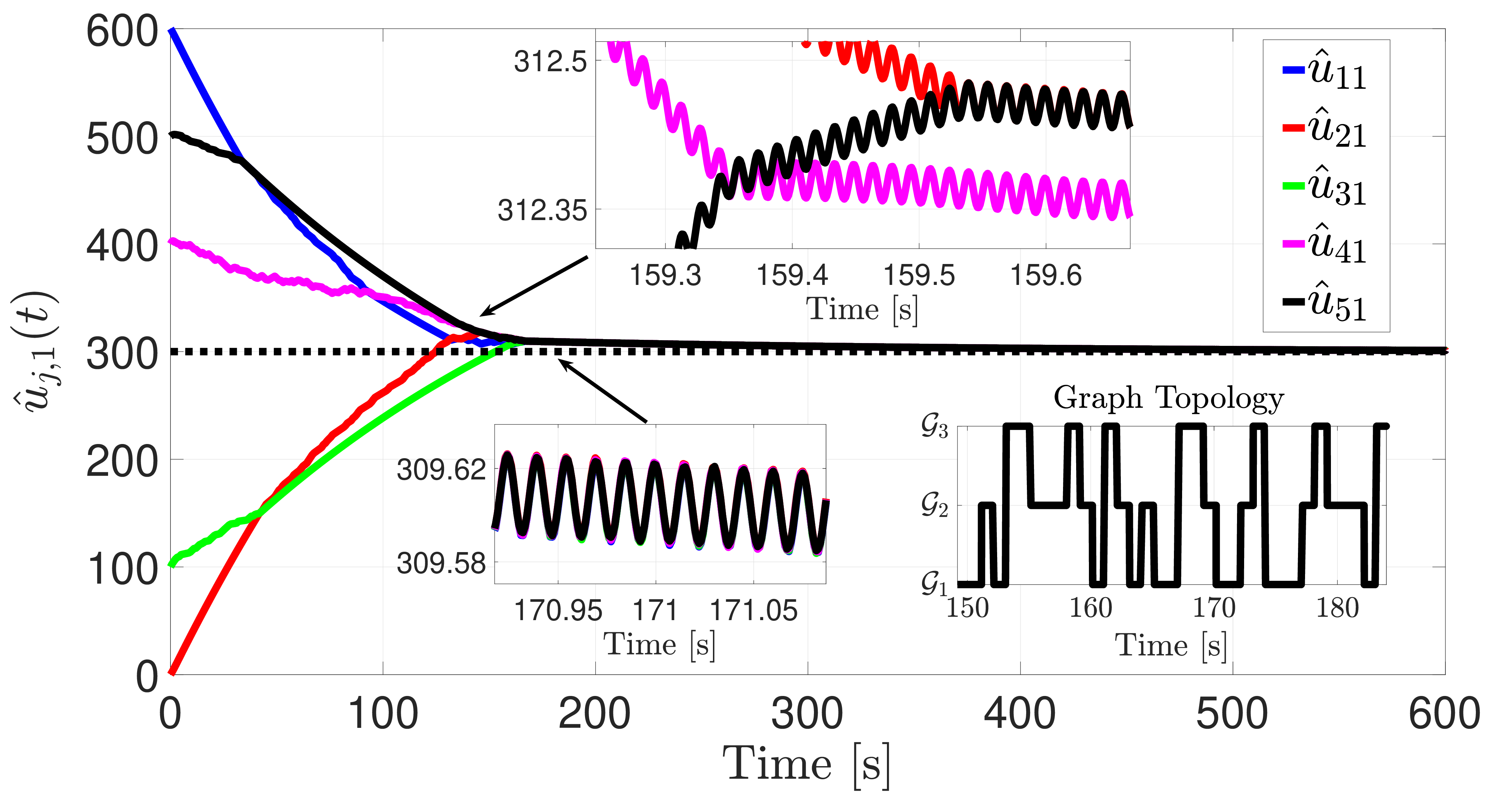}
\includegraphics[width=0.95\linewidth]{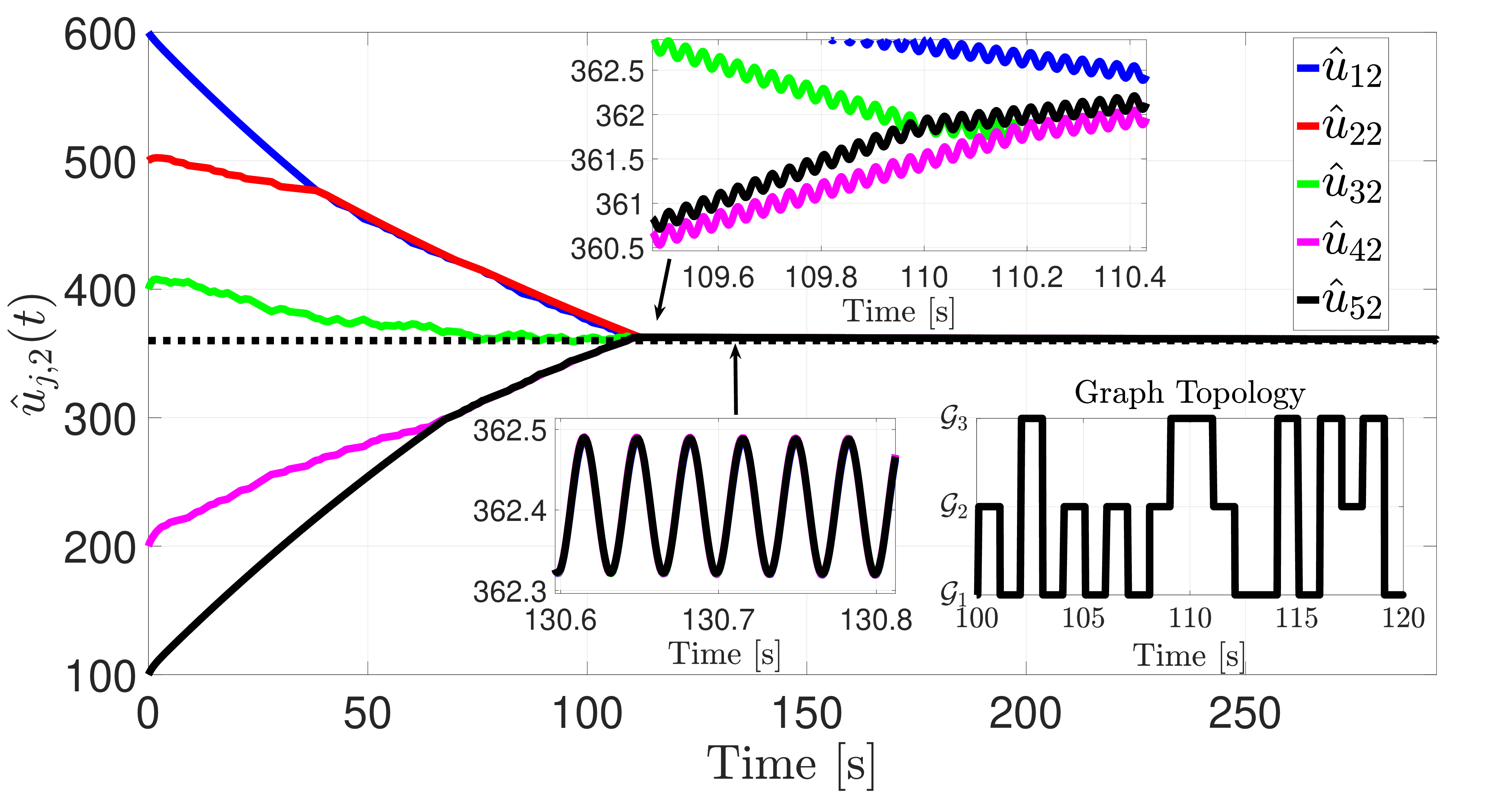}
\includegraphics[width=0.95\linewidth]{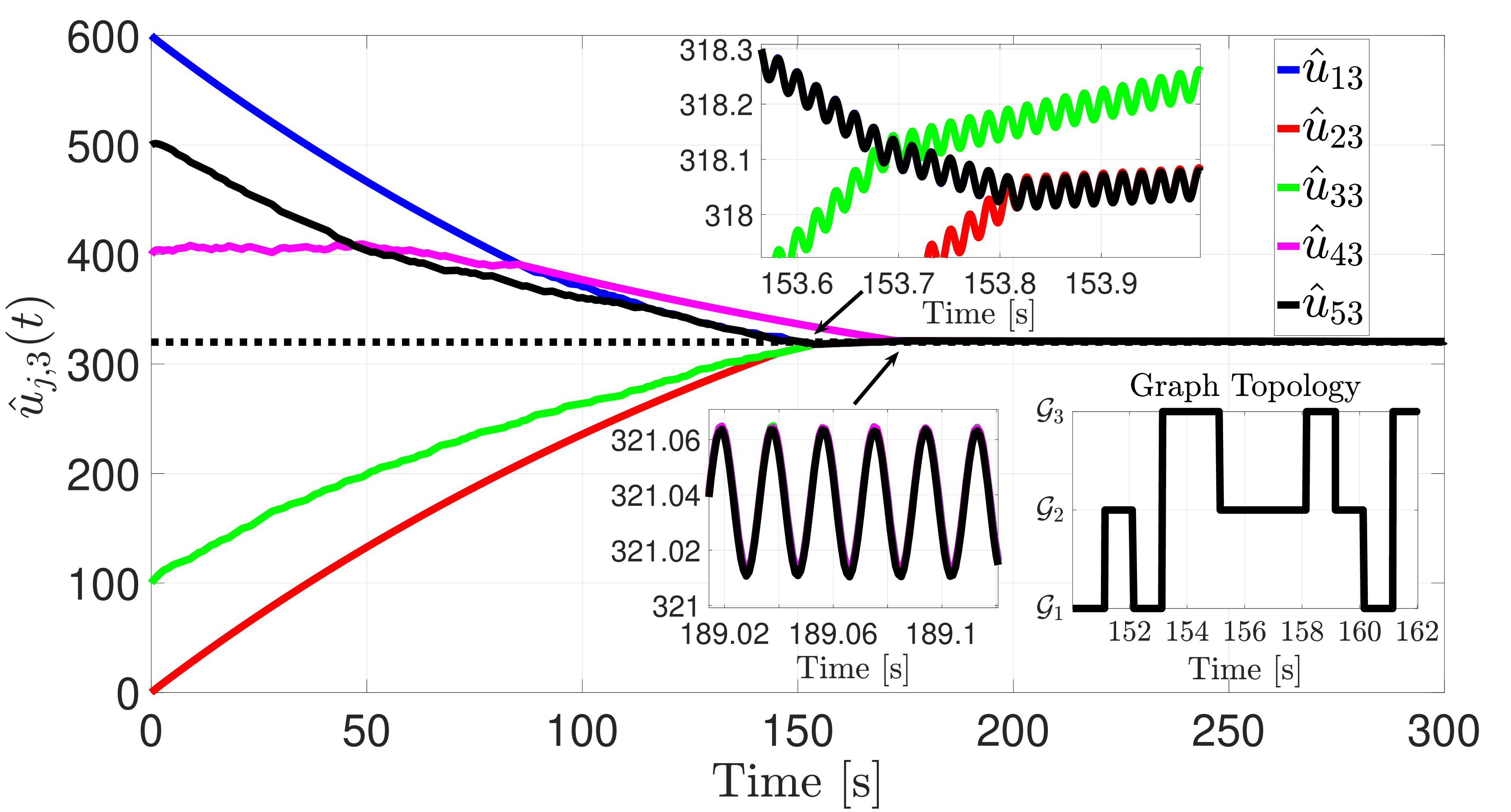}
\includegraphics[width=0.95\linewidth]{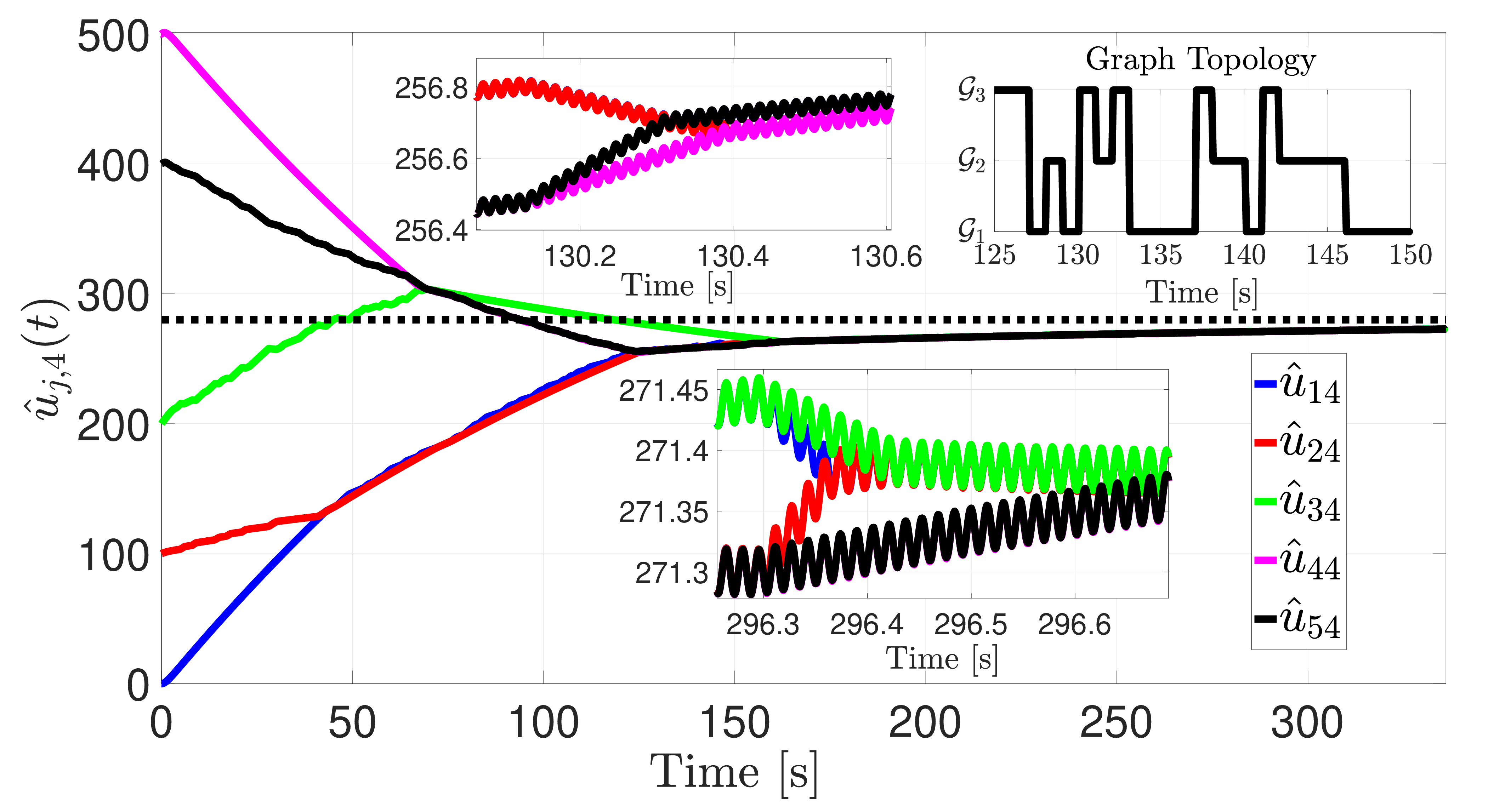}
\includegraphics[width=0.95\linewidth]{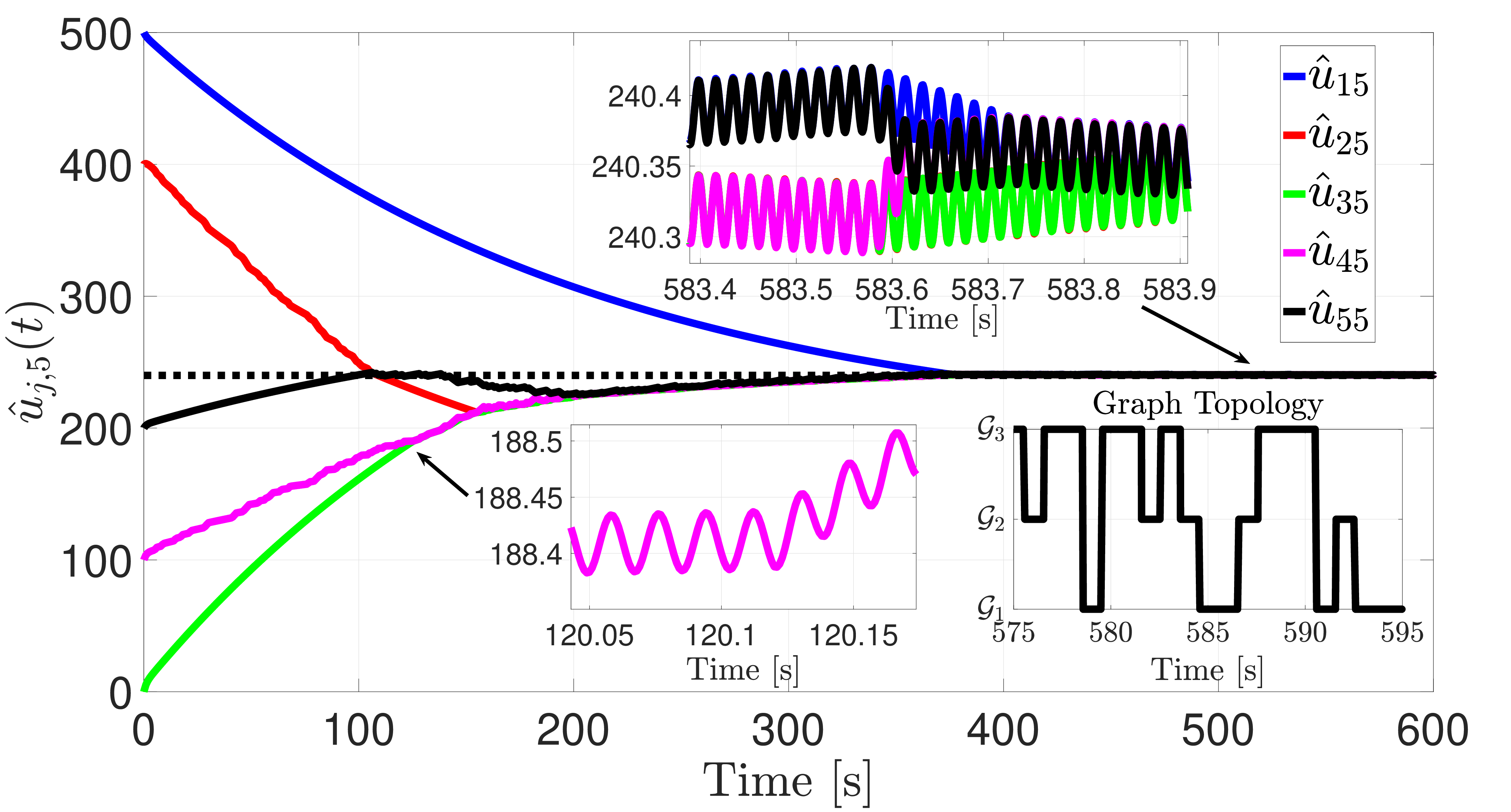}
\end{tcolorbox}
 \caption{Trajectories of hybrid set-seeking algorithm \eqref{udynamicsi} with communication graphs allowed to switch every 0.1 seconds. \label{fig11}}
\end{figure}

The main takeaway from Proposition \ref{proposition_arbitrary} is that seeking algorithms can be designed and analyzed to switch arbitrarily fast between a finite number of modes, provided that their continuous-time target decision-making dynamics share a common Lyapunov function.  In this case, all (non-discrete) solutions of the switching decision-making dynamics \eqref{arbitrary_switching} can be generated by the differential inclusion \eqref{arbitary_switching}. 
Note that common Lyapunov functions satisfying the conditions of Proposition \ref{proposition_arbitrary} emerge frequently in optimization algorithms, since the cost function $J$ usually serves as a common Lyapunov function for different optimization dynamics. For example, this is the case in the standard gradient ascent and Newton methods, as well as related dynamics such as the Newton variable structure algorithm or the continuous Jacobian matrix transpose algorithm, for quadratic cost functions (see \cite[Ch. 2]{book_algorithms}). A common Lyapunov function also emerges frequently in algorithms for distributed optimization, and learning in multi-agent systems with undirected switching communication topologies \cite{Olfati:04,Poveda_WC17}. %It also emerges frequently in many learning dynamics for \emph{potential games} \cite[Ch. 7.1]{Sandholm1}. 
\begin{example}[Switching Multi-Agent Seeking Systems]
Consider a multi-agent system where $N$ agents seek to cooperatively solve the following consensus-optimization problem:
\begin{equation}\label{optimizationdistributed}
\min~J(u)=\sum_{i=1}^N J_i(u_i),~~u_i\in\mathbb{R}^{n_i}.
\end{equation}
Agents can communicate only to neighbors, who are characterized by an undirected graph that is allowed to change topology arbitrarily fast, but which is connected at all times. Since \eqref{optimizationdistributed} is a well-studied problem in the literature whenever the gradients $\nabla J_i$ are known to the agents, there is a plethora of potential target nominal dynamics \eqref{hybrid_learning_dynamics} that could be considered for the design of an ES system \cite{lin2016distributed}. To illustrate the use of Krasovskii solutions, we consider the dynamics
\begin{equation}\label{udynamicsi}
\dot{\hat{u}}_i\in 
\hat{F}_{\delta,i}(\hat{u},\xi_i)=\frac{k}{\delta}\sum_{j\in\mathcal{N}_i(t)}\overline{\text{sign}}(\hat{u}_j-\hat{u}_i)-k\gamma \xi_i,
\end{equation}
where $\overline{\text{sign}}$ is the set-valued map defined in \eqref{krasovskiicons}, i.e.,
\begin{equation}
\overline{\text{sign}}(z)=\left\{\begin{array}{cl}
\{1\} &~\text{if}~z>0\\
{[}-1,1{]} &~\text{if}~z=0\\
\{-1\} & ~\text{if}~z<0
\end{array}\right..
\end{equation}
To generate $\xi_i$, each agent implements its own autonomous derivative estimator \eqref{oscillatordynamics} using its own cost function $J_i$ as input. In this way, agents share information with neighbors only via the dynamics \eqref{udynamicsi}. By indexing all possible graph topologies using a set of logic modes $\mathcal{Q}=\{1,2,\ldots,q_{\max}\}$, the resulting system has the form \eqref{arbitrary_switching}, where each fixed logic mode $q$ represents a particular graph topology. Figure \ref{fig11} presents a numerical example with 5 agents, showing the evolution in time of the states $\hat{u}_i$, which converge to the optimal solution. The cost functions used in the simulations are given by $J_i=\frac{1}{2}u_i^\top u_i-100b_i^\top u_i$, where $b_1=(1,2,3,4,5)^\top$, $b_2=(5,4,3,2,1)^\top$, $b_3=(2,3,4,5,1)^\top$, $b_4=(3,4,5,1,2)^\top$, $b_5=(4,5,1,2,3)^\top$, and the ES parameters are $\varepsilon_a=0.8$, $k=0.05$, $k_f=0.1$, $\varepsilon_l=0.1$, $\omega=500(5/6,3/8,2/3,4/5,7/10)$, $\alpha=15$, $\gamma=0.1$. It can be verified that the target nominal system \eqref{udynamicsi} with $\xi_i=\nabla J_i$ and sufficiently large $\delta$ admits a common Lyapunov function across all connected undirected graphs \cite{lin2016distributed}. It follows that Proposition \ref{proposition_arbitrary} holds and therefore, by Theorem \ref{theorem1}, the set-seeking dynamics achieve practical convergence to the set of solutions of \eqref{optimizationdistributed}.

\end{example}

\subsection{Seeking with Average Dwell-time Switching}
When a common Lyapunov function does not exist for the target switching system \eqref{arbitrary_switching}, it is still possible to stabilize $\mathcal{O}$ in a semi-global practical sense by restricting how fast the state $q$ switches between modes. For example, this can be achieved by imposing a dwell-time or an average dwell-time constraint on the switching signal. As discussed in Example \ref{intermittentsourceseeking}, systems with jumps that satisfy these types of constraints can be modeled by using a hybrid automaton with state $\tau_1\in\mathbb{R}_{\geq0}$ and dynamics
\begin{subequations}\label{dwell_time_dynamics}
\begin{align}
\dot{\tau}_1&\in[0,\eta_1],~~~~\tau_1\in [0,N_0]\\
\tau_1^+&=\tau_1-1,~~~\tau_1\in [1,N_0],
\end{align}
\end{subequations}\noindent
where $\eta_1\in\mathbb{R}_{>0}$ and $N_0\in\mathbb{R}_{\geq1}$. In particular, every solution of \eqref{dwell_time_dynamics} has a hybrid time domain $E$ with the property that for each of its elements $(s,i),(t,j)\in E$ with $s+i\leq t+j$, the following inequality holds:
\begin{equation}\label{averagedwellinequality}
N(t,s):=j-i\leq \eta_1(t-s)+N_0,
\end{equation}
where $N(t,s)$ denotes the number of jumps in the interval $[s,t]$, see  \cite{Average_Dwell_time}. Therefore, when considering the joint hybrid dynamics \eqref{arbitrary_switching} and \eqref{dwell_time_dynamics}, the hybrid time domains of the resulting system also satisfies inequality \eqref{averagedwellinequality}, thus effectively imposing a bound on how frequently $q$ can switch in \eqref{arbitrary_switching}. When $N_0=1$, inequality \eqref{averagedwellinequality} reduces to the well-known \emph{dwell-time} condition, see \cite{Dwell_time}. When $N_0\in\mathbb{Z}_{>1}$, inequality \eqref{averagedwellinequality} describes an \emph{average dwell-time} condition. It turns out that any hybrid time domain satisfying \eqref{averagedwellinequality} can be generated by the hybrid automaton  \eqref{dwell_time_dynamics}, see \cite[Proposition 1.1]{SmoothLyaHDS}. Therefore, these dynamics can be used to study switching set-seeking systems under dwell-time and average dwell-time constraints on the switching signal.
\begin{proposition}[Seeking with Slow Switching]\label{proposition_dwell}
Suppose that for each $q\in \mathcal{Q}$ the mapping $f_q$ in \eqref{flow_optimizer0a1} is OSC, LB, and convex-valued. Moreover, suppose that for each $q\in\mathcal{Q}$ the $\hat{u}$-dynamics in system \eqref{flow_optimizer0a1} render UGAS the set $\mathcal{O}$. Then, the HDS generated by Equations (\ref{arbitrary_switching})  and (\ref{dwell_time_dynamics}) has the structure of (\ref{hybrid_learning_dynamics}) with $r=2$, $\delta=\eta_1$, $\hat{z}:=(q,\tau)^\top$, $C_{u,z}:=\hat{\mathbb{U}}\times \mathcal{Q}\times[0,N_0]$, $D_{u,z}:=\hat{\mathbb{U}}\times \mathcal{Q}\times[1,N_0]$, $\Psi:= \mathcal{Q}\times [0,N_0]$, $\hat{F}_{\delta}:=f_q\times\{0\}\times[0,\eta_1]$, and $\hat{G}_{\delta}:=\{\hat{u}\}\times {Q\backslash\{q\}}\times \{\tau_1-1\}$. Moreover, this HDS satisfies all the items in Assumption \ref{keyassumptions}. \QEDB 
\end{proposition}

The result of Proposition \ref{proposition_dwell}, which follows from 
\cite[Corollary 7.28]{bookHDS}, enables the study of switching seeking dynamics where each seeking mode is able to individually stabilize the set $\mathcal{O}$, and where the switching between modes is "sufficiently slow" on average. Such types of problems might emerge in multi-agent systems with switching directed communication graphs for which a common Lyapunov function do not exist, as well as in seeking problems with slow variations in the landscape of the cost function, or problems that require a persistent switching of the adaptation gains to achieved a desired performance. 

The result of Proposition \ref{proposition_dwell} can also be extended to ES systems with switching cost functions, even when the cost functions do not share the same minimizer $u^*_q\in\mathbb{R}^n$, as in Example \ref{intermittentsourceseeking}. For instance, when $f_q(\cdot)$ is single-valued and continuous, for a sufficiently large compact set $K$ of initial conditions, the (practical) stability of the resulting switching  dynamics \eqref{arbitrary_switching} and \eqref{dwell_time_dynamics} can be studied with respect to the Omega-limit set from $K$, denote $\Omega(K)$, which is given by $\Omega(K)=\bigcup_{q\in\mathcal{Q}} \Omega_q$, where
\begin{equation}\label{omegalimitset}
\Omega_q=\bigcap_{j\in\mathbb{Z}_{\geq0}}\overline{R_0\left(\left(\{u^*_q\}+\frac{1}{j+1}\mathbb{B}\right)\times\{q\}\times[0,N_0]\right)},
\end{equation}
and where $R_0(K)$ denotes the reachable set from $K$. An illustration of this set is presented in Figure \ref{fig3} for the source-seeking problem with persistent surveillance. For further details on switching systems with distinct equilibrium points, we refer the reader to \cite{baradaran2020omega}.
\subsubsection{Seeking with Switching Unstable Modes}
\begin{figure}[t!]
\begin{tcolorbox}[colback=ivoryA, colframe=ivoryA]
\centering
\includegraphics[width=0.95\linewidth]{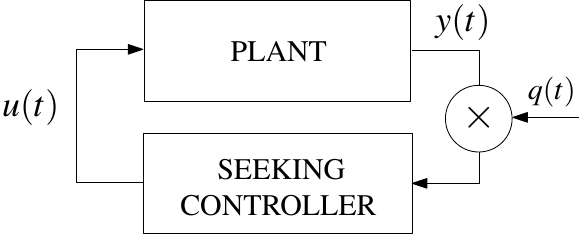}
\end{tcolorbox}
\caption{Scheme of ES controller under persistent multiplicative attacks ($q\in\{-1,1\}$) on the output measurements.}
\label{fig10023}
\end{figure}
In many applications, such as source-seeking missions, seeking algorithms are often implemented in highly adversarial and dynamic environments. In these settings, intermittent communication, sensing, or actuation failures are unavoidable due to external factors. Such intermittent failures can lead to unstable behaviors, which, however, can be corrected if the system is able to return to its nominal operating condition "sufficiently" often. This scenario can be analyzed using the framework of hybrid set-seeking systems by considering both stable and unstable modes in \eqref{arbitrary_switching}. Specifically, the compact set \(\mathcal{Q}:=\{1, \ldots, q_{\max}\}\) can be partitioned as \(\mathcal{Q} = \mathcal{Q}_u \cup \mathcal{Q}_s\), where the modes \(q \in \mathcal{Q}_s\) represent stable dynamics and the modes \(q \in \mathcal{Q}_u\) represent unstable dynamics. For this type of systems, stability can be guaranteed as long as the total activation time $T(s,t)$ of unstable modes during any time interval $(s,t)$ satisfies the bound
\begin{equation}\label{time_ratio_constraint}
T(s,t):=\int_{s}^t \mathbb{I}_{Q_u}(q(r,j(r)))dr\leq T_0+\eta_2(t-s),
\end{equation}\noindent
where $\eta_2\in[0,1)$, $T_0\in\mathbb{R}_{\geq0}$, and $j(r)$ is the minimum $j$ such that $(r,j)\in E$. Similarly to \eqref{dwell_time_dynamics} and \eqref{averagedwellinequality}, 
the average activation-time constraint \eqref{time_ratio_constraint} can also be imposed on the switching signal $q$ by using the following hybrid automaton with state $\tau_2\in\mathbb{R}_{\geq0}$ and dynamics
\begin{subequations}\label{time_ratio}
\begin{align}
\dot{\tau}_2&\in[0,\eta_2]-\mathbb{I}_{\mathcal{Q}_u}(q),~~~~~~\tau_2\in [0,T_0],\\
\tau_2^+&=\tau_2,~~~~~~~~~~~~~~~~~~~~~~~\tau_2\in [0,T_0],
\end{align}
\end{subequations}\noindent
where $\mathbb{I}_{\mathcal{Q}_u}(\cdot)$ corresponds to the standard indicator function. In particular, every solution of the HDS given by \eqref{arbitrary_switching} and \eqref{time_ratio} has a hybrid time domain $E$, such that for each of its elements $(s,i),(t,j)\in E$ with $s+i\leq t+j$ and signal $q:E\rightarrow \mathcal{Q}$ we have that the bound \eqref{time_ratio_constraint} holds. In fact, each hybrid time domain $E$ with elements $(s,i),(t,j)\in E$ with $s+i\leq t+j$, satisfying \eqref{time_ratio_constraint} can be generated by the hybrid dynamics \eqref{arbitrary_switching} and \eqref{time_ratio}. Therefore, by combining dynamics \eqref{arbitrary_switching}, \eqref{dwell_time_dynamics}, and \eqref{time_ratio}, we can design different types of hybrid algorithms \eqref{hybrid_learning_dynamics} that satisfy Assumption \ref{keyassumptions} even when the decision-making algorithms operate under a persistent loss of communication, sensing, or actuation.

% \begin{lemma}\label{Lemma_time_ratio_constraint}
% For each solution of the HDS (\ref{time_ratio}), corresponds a hybrid time domain $E$, such that for each of its elements $(s,i),(t,j)\in E$ with $s+i\leq t+j$ and signal $q:E\rightarrow Q$, with $Q_u\subset Q$, we have that
% %
% \begin{equation}\label{time_ratio_constraint}
% T(s,t):=\int_{s}^t \mathbb{I}_{Q_u}(q(r,j(r)))dr\leq T_0+\eta_2(t-s),
% \end{equation}\noindent
% %
% where $T(s,t)$ denotes the total activation time of modes $q\in Q_u$ between times $(s,t)$, and where $j(r)$ is the minimum $j$ such that $(r,j)\in E$. Moreover, each hybrid time domain $E$ with elements $(s,i),(t,j)\in E$ with $s+i\leq t+j$ satisfying the bound (\ref{time_ratio_constraint}), can be generated by  (\ref{time_ratio}).
% \end{lemma}
%

\begin{figure}[t!]
\begin{tcolorbox}[colback=iceblue!80, colframe=iceblue!80]
\centering
\includegraphics[width=0.99\linewidth]{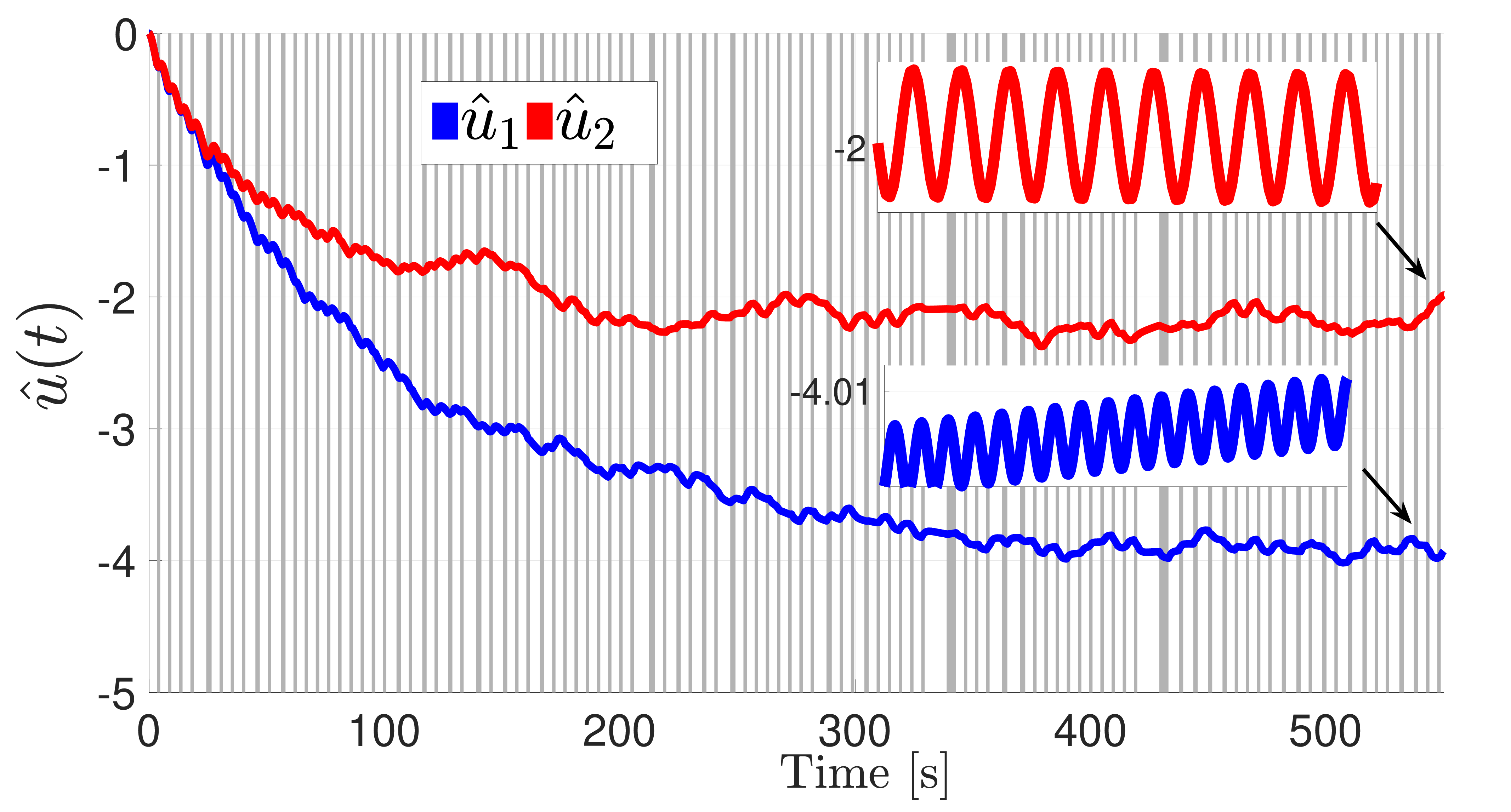}
\includegraphics[width=0.99\linewidth]{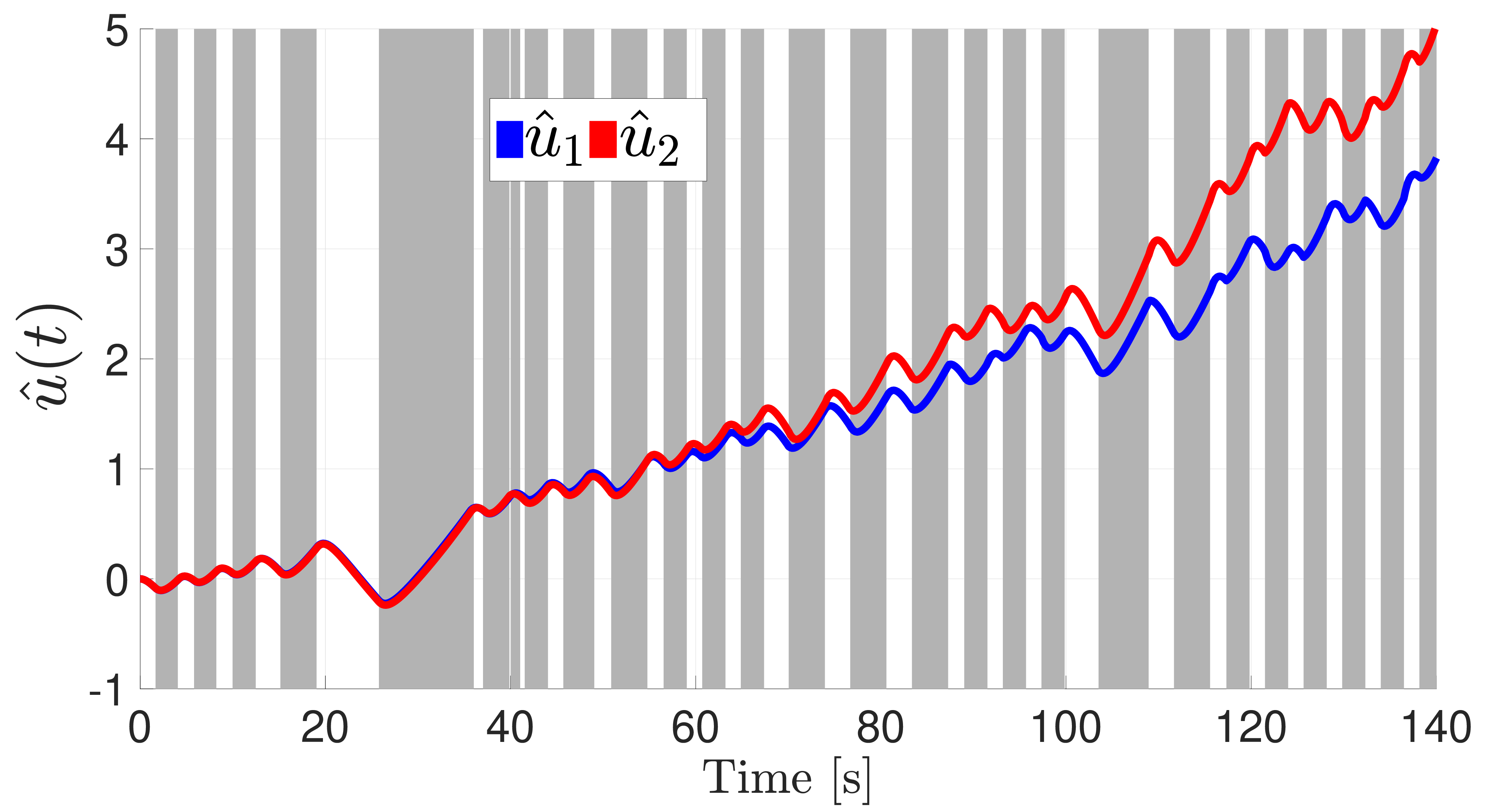}
\end{tcolorbox}
\caption{\tcb{ES dynamics in $\mathbb{R}^2$, under persistent attacks on the gradient estimation mechanism. The intensity of the attacks is modeled using \eqref{time_ratio}. As shown in the upper plot, when $\eta_2$ is sufficiently small, the ES dynamics preserve stability. The lower plot shows the unstable behavior that emerges when $\eta_2$ is moderate.} \label{fig1000987}}
\end{figure}

\begin{proposition}[Seeking with Unstable Dynamics]\label{time_ratio_constraint_stability}
Suppose that for each $q\in \mathcal{Q}$ the mapping $f_q$ in \eqref{arbitrary_switching} is OSC, LB, and convex-valued, and that there exist $\chi\in\mathbb{R}_{\geq1}$,  $\alpha_{1,q},\alpha_{2,q} \in\mathcal{K}_{\infty}$, $\lambda_s\in\mathbb{R}_{>0}$, $\lambda_u\in\mathbb{R}_{>0}$ and smooth functions $V_q:\hat{\mathbb{U}}\rightarrow\mathbb{R}_{\geq0}$, such that the following inequalities hold:
\begin{subequations}\label{conditions_unstability}
\begin{align}
\alpha_{1,q}(|\hat{u}|_{\mathcal{O}})\leq& V_q(\hat{u})\leq\alpha_{2,q}(|\hat{u}|_{\mathcal{O}}),~~\forall~(q,\hat{u})\in \mathcal{Q}\times \hat{\mathbb{U}},\label{condition1_time_ratio}\\
\langle \nabla V_{q_s}, \tilde{f} \rangle \leq& -\lambda_sV_{q_s}(\hat{u}),~\forall~(q_s,\hat{u},\tilde{f})\in \mathcal{Q}_s\times \hat{\mathbb{U}}\times f_q,\label{condition3_time_ratio}\\
\langle \nabla V_{q_u}, \tilde{f} \rangle \leq& \lambda_uV_{q_u}(\hat{u}),~~~\hspace{0.05cm}\forall~(q_u,\hat{u},\tilde{f})\in \mathcal{Q}_u\times \hat{\mathbb{U}}\times f_q.\label{condition4_time_ratio}\\
V_p(\hat{u})\leq&\chi V_q(\hat{u}), ~~~~~~\forall~(q,p),\hat{u}\in \mathcal{Q}\times \hat{\mathbb{U}}.\label{condition2_time_ratio}
\end{align}
\end{subequations}\noindent
If the above parameters satisfy the following inequality:
\begin{equation}\label{conditionstability}
\lambda_s>\eta_1\log(\mu)+\eta_2(\lambda_s+\lambda_u),
\end{equation}
then, the HDS composed by equations (\ref{arbitrary_switching}), (\ref{dwell_time_dynamics}), and (\ref{time_ratio}), has the structure of (\ref{hybrid_learning_dynamics}) with $r=3$ $\hat{z}:=(q,\tau_1,\tau_2)$, $C_u=D_u=\hat{\mathbb{U}}$, $C_z=\mathcal{Q}\times[0,N_0]\times[0,T_0]$, $D_z=\mathcal{Q}\times[1,N_0]\times[0,T_0]$, and $\Psi:=\mathcal{Q}\times[0,N_0]\times[0,T_0]$. Moreover, this HDS satisfies the items in Assumption \ref{keyassumptions}, and the solutions' hybrid time domains satisfy the bounds \eqref{averagedwellinequality} and \eqref{time_ratio_constraint}.\QEDB 
\end{proposition}

The Lyapunov-based conditions \eqref{conditions_unstability} and the inequality \eqref{conditionstability} are common in the literature of switching systems with unstable modes. In this sense, Proposition \ref{time_ratio_constraint_stability} simply recasts such conditions using the hybrid automaton \eqref{dwell_time_dynamics} and the hybrid monitor \eqref{time_ratio} such that, when interconnected with the target nominal dynamics \eqref{arbitrary_switching}, the resulting system has the form \eqref{hybrid_learning_dynamics} and satisfies all the conditions needed for the analysis of switching set-seeking systems with unstable modes. We note that verification of the assumptions in Proposition \ref{time_ratio_constraint_stability} can be straightforward for some problems and algorithms, particularly when the cost $J$ is quadratic, which leads to linear gradients $\nabla J$.  In fact, conditions (\ref{condition1_time_ratio}) and (\ref{condition3_time_ratio}) simply ask for UGAS of $\mathcal{O}$ for each mode $f_{q}$, where $q\in \mathcal{Q}_s$. Similarly, condition (\ref{condition4_time_ratio}) simply asks that solutions of the differential inclusions associated to the unstable modes $q\in \mathcal{Q}_u$  have no finite escape times. Additionally, condition (\ref{condition2_time_ratio}) is easily satisfied if the Lyapunov functions are quadratic. In this case, conservative estimates of the parameters $\lambda_s$, $\lambda_u$, and $\mu$, can be used to design the parameters $\eta_1$ and $\eta_2$ such that \eqref{conditionstability} holds. 

To illustrate Proposition~\ref{time_ratio_constraint_stability} as well as the discussion of Remark \ref{remarkunstable}, we consider an ES algorithm with standard gradient dynamics \( \dot{\hat{u}} = -\xi \), operating in an environment where the sign of the measurements of \( y \) is corrupted in real time. Specifically, we consider a quadratic cost function $J(u)=0.01u^\top Q u +b^\top u$, where  $Q=[1,1/2;1/2,3/2]$, $b=(0.1,0.1)$, and we define \( q \in \mathcal{Q}=\{-1, 1\} \), where \( q = -1 \) indicates that the output has been corrupted. See Figure~\ref{fig10023} for a block-diagram representation. Figure~\ref{fig1000987} shows the time evolution of the resulting hybrid set-seeking dynamics with two different average activation times of the mode \( q = -1 \), using $\omega=(810,420)$, $\varepsilon_a=0.1$ and $\varepsilon_f=1$. As observed, when the activation time is sufficiently short (i.e., \( \eta_2 \) is small), the seeking system successfully converges to the optimal point. In contrast, when \( \eta_2 \) is too large, the system becomes unstable. For related applications of ES in the context of deceptive attacks on gradient estimates, we refer the reader to \cite{GalarzaAttacks}. Note that by redefining the set \( \mathcal{Q} \) as \( \mathcal{Q} = \{0, 1\} \), the previous model can also capture ES systems operating under sporadic feedback. In this case, source-seeking can still be achieved, provided that the condition \eqref{time_ratio_constraint} is satisfied.
 Multi-agent seeking problems involving persistent adversarial agents can also be analyzed using similar techniques; see \cite[Sec.~6.3]{PoTe17Auto}.

\begin{pullquote}
The study of adaptive seeking systems with logic modes and event-based rules opens new avenues for research at the intersection of real-time, model-free control and computer science, with applications in cyber-physical systems.
%The study of switching ES systems with both stable and unstable modes paves the way for research on cyber-security, resilience, and robustness analysis of current methods. 
\end{pullquote}

%%%%%%%%%%%%%%%%%%%%%%%%%%%%
\subsubsection{Seeking with Slowly Varying or Jumping Parameters}
%%%%%%%%%%%%%%%%%%%%%%%%%%%&
%
In many practical seeking problems, the parameters associated with the plant or the algorithm vary slowly over time. For instance, this occurs when the response map \( J \) is time varying or exhibits mild jumps on some of its parameters. Although slowly-varying seeking problems have been studied under the assumption that the target decision-making algorithm  is input-to-state stable (ISS) with respect to the time derivatives of the varying signals in the cost \cite{scheinker2012extremum,poveda2021fixed,labar2022iss}, the ISS property may not necessarily hold for a broad class of decision-making algorithms \eqref{hybrid_learning_dynamics}. Nevertheless, such problems can still be addressed using the hybrid set-seeking systems framework studied in this paper. In particular, by letting \( q \) represent the time-varying parameter (or its derivative, depending on the application), redefining \( \mathcal{Q} \) in \eqref{arbitrary_switching} as a compact but not necessarily finite set, and by allowing the flow dynamics in the hybrid system \eqref{arbitrary_switching} to depend explicitly on \( q \), i.e., \( \dot{\hat{u}} \in f(\hat{u}, \nabla J(\hat{u}), q) \), the slowly varying, weakly jumping behavior of \( q \) within \( \mathcal{Q} \) can be effectively modeled via the inclusions
 \begin{equation}\label{q_dynamics}
 \dot{q}\in\eta_3\mathbb{B}, ~~~~\text{and}~~~~q^+\in q+\eta_3\mathbb{B},~~\text{where}~~~\eta_3\in\mathbb{R}_{>0},
 \end{equation}\noindent
which replaces the $q$-dynamics in (\ref{arbitrary_switching}). In this scenario, we assume that for each \( q \in \mathcal{Q} \), there exists an optimal compact set \( \mathcal{O}_q \subset \hat{\mathbb{U}} \) that minimizes $J$, and we employ a simple average dwell-time automaton of the form \eqref{dwell_time_dynamics} to rule out purely discrete-time solutions in the hybrid system. The following proposition is a direct consequence of \cite[Corollary 7.27]{bookHDS}.
\begin{proposition}[Seeking under Slowly-Varying Parameters]\label{time_varying_parameter}
Let $\mathcal{O}_q\subset \hat{\mathbb{U}}$ and suppose that for each $q\in \mathcal{Q}$ the mapping $f_q$ in \eqref{flow_optimizer0a1} satisfies item (a) in Assumption \ref{keyassumptions} with the flow map given by $f(\hat{u},\nabla J (\hat{u}), q)$ and the $q$-dynamics given by (\ref{arbitrary_switching}). Suppose also that for each $q\in \mathcal{Q}$  the dynamics $\dot{\hat{u}}\in f(\hat{u},\nabla J (\hat{u}), q)$ render the set $\mathcal{O}_q$ UGAS. Then, this HDS combined with the dynamics (\ref{dwell_time_dynamics}) has the structure of (\ref{hybrid_learning_dynamics}) with $r=2$, $\delta=\eta_3$, $\hat{z}:=(q,\tau_1)$, $C_u=D_u=\hat{\mathbb{U}}$, $C_z=\mathcal{Q}\times[0,N_0]$, $D_z:=\mathcal{Q}\times[1,N_0]$, $\mathcal{O}:=\{(\hat{u},q):\hat{u}\in\mathcal{O}_q, q\in \mathcal{Q}\}$, and $\Psi:=\mathcal{Q}\times [0,N_0]$. Moreover, this HDS satisfies all the items in Assumption \ref{keyassumptions}. \QEDB 
\end{proposition}

An application of Proposition \ref{time_varying_parameter} was already illustrated in Example \ref{projectedexample} and in Figure \ref{figprojected1}, where a slowly varying parameter $q$ defines the cost function $J$. In particular, when $\dot{q}=0$, the state $q$ remains constant and the trajectories $\hat{u}$ of the system converge to a neighborhood of the optimal point $u^*(q)$. Conversely, when $\dot{q}\in \eta_3 \mathbb{B}$, with $\mathcal{Q}$ compact and $\eta_3>0$ sufficiently small, the \emph{tracking error} $\tilde{u}=\hat{u}-u^*(q)$ converges to a neighborhood of zero. %Moreover, if the target decision-making dynamics exhibit stronger stability properties, such as input-to-state stability, one can also establish suitable convergence bounds for any $\eta_3>0$, see \cite{chen2025continuous}.
\subsection{Seeking in Noncontractible Spaces: Source-seeking with Obstacles}
Smooth controllers, and dynamical systems in general, face fundamental limitations in stabilization problems when operating in complex spaces, particularly when obstacles are present. This is a consequence of the fact that obstacles usually preclude contractibility of the operational space, which is a necessary condition for \emph{global} stabilization \cite{sontag1999stability} under smooth feedback laws. For example, it has been shown that in obstacle avoidance problems, different types of controllers, including those based on navigation functions or CBFs can generate ``deadlocks'' \cite{jankovic2023multiagent} and undesirable spurious equilibria in the closed-loop system \cite{marley2024hybrid}. While it is well known that such undesirable points—often forming a set of measure zero—can be ``pushed out" to the boundary of the feasible set~\cite{koditschek1990robot}, introducing even arbitrarily small perturbations from measurement noise or external disturbances in the closed-loop system can cause these problematic points in the state space to no longer be of measure zero, thereby precluding (almost) global convergence guarantees. However, these obstructions to \emph{robust} global stabilization under obstacles, which also apply to optimization dynamics, can be overcome using hybrid set-seeking systems \cite{ObstacleAvoidance,Strizic:17_CDC,OchoaPovedaManifolds,Poveda:20TAC}. For example, consider a vehicle operating on the plane according to the simple point-mass dynamics
\begin{equation}\label{plant_obstacle}
\dot{p}=\alpha,~~~~p\in\mathbb{R}^2,
\end{equation}
\begin{figure*}[t!]
\begin{tcolorbox}[colback=iceblue!80, colframe=iceblue!80]
\centering
\includegraphics[width=0.47\linewidth]{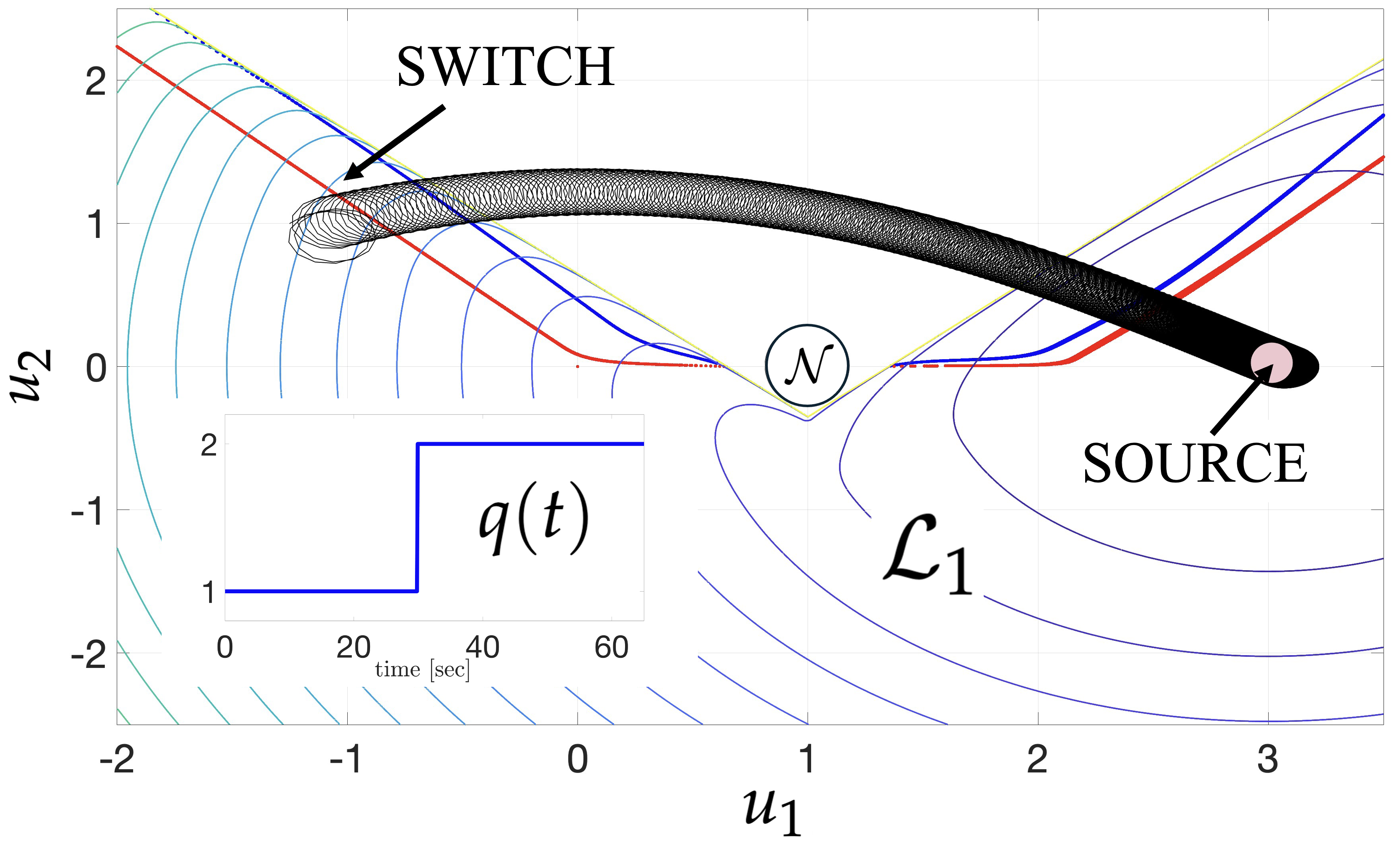}\includegraphics[width=0.506\linewidth]{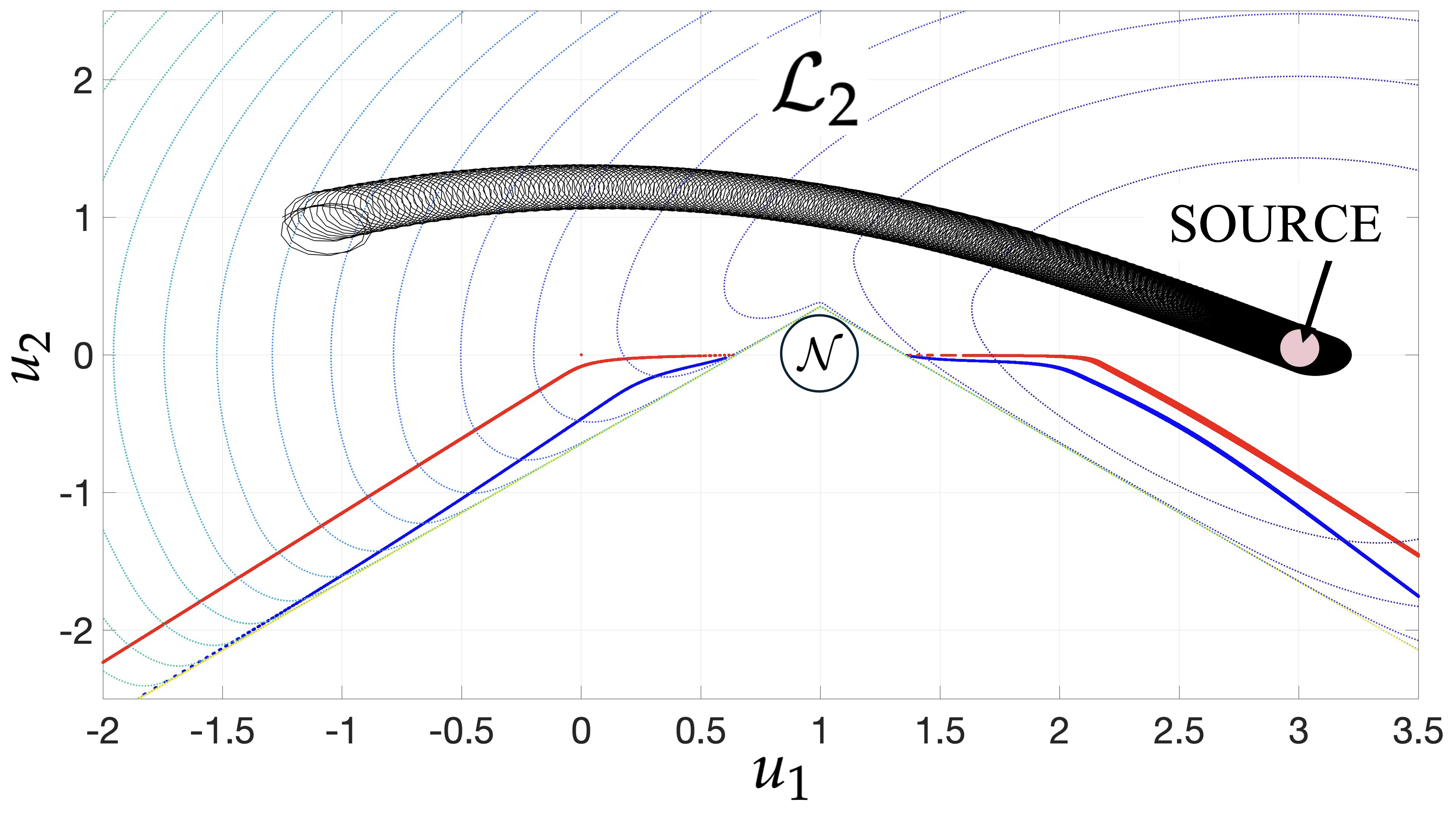}
\end{tcolorbox}
\caption{Trajectories of vehicle in $\mathbb{R}^2$, controlled via a hybrid source-seeking algorithm for obstacle avoidance. The trajectories are shown over the virtual partitions $\mathcal{L}_1$ and $\mathcal{L}_2$. The system is initialized in \( \mathcal{L}_1 \) with \( q = 1 \), but eventually switches to mode \( q = 2 \) to evolve in \( \mathcal{L}_2 \), until it converges to a neighborhood of the source.
 \label{fig11203}}
\end{figure*}
where $p=(p_1,p_2)$ is the position and $\alpha\in\mathbb{R}^2$ is the control input. The main goal of the vehicle is to "discover" in real time the location $\mathcal{O}=\{u^\star\}\in\mathbb{R}^2$ where a potential field $J$ attains its maximum value, using only intensity measurements, while simultaneously avoiding an obstacle  $\mathcal{N}\subset\mathbb{R}^2$ in the plane. As discussed in \cite{sontag1999stability} and \cite{ObstacleAvoidance},  this problem cannot be robustly solved using target decision-making algorithms \eqref{hybrid_learning_dynamics} characterized by smooth feedback dynamics due to the topological obstructions introduced by the obstacle. Therefore, in order to achieve source-seeking and obstacle avoidance, we consider a hybrid source-seeking controller with feedback law given by
%
% \begin{figure*}[t!]
%     \centering
%     \includegraphics[width=1\columnwidth]{figs/O1.eps}
%     \includegraphics[width=1\columnwidth]{figs/O2.eps}
%     \caption{The left plot shows the trajectories of the vehicle \eqref{plant_obstacle} controlled via \eqref{hybridcontroller01}, evolving over the operational space defined by the logic state $q=1$. Similarly, the right plot corresponds to $q=2$. The trajectory with the color blue corresponds to $b=15$. The green trajectory corresponds to $b=-15$, and the black trajectory is obtained when $b(t)=15\sin(t)$. \tcb{The purple trajectory shows the instability obtained when using the model-based hybrid controller of \cite{ObstacleAvoidance} (which assumes $b>0$) with $b=-15$.} } \label{fig:O1}
% \end{figure*}
%
\begin{equation}\label{hybridcontroller01}
\alpha=\varepsilon_a\omega \mathbf{R}_0\mu+\xi.
\end{equation}
Assuming that we have access to measurements of $J$ at the points $p$ (to be used by the filter \eqref{filterdynamics}) and $p-\varepsilon_a\mu$ (to be used in the flow and jump sets), which is possible via collocation of sensors, we consider the change of variable $\hat{u}=p-\varepsilon_a\mu$ to obtain the dynamics $\dot{\hat{u}}=\xi$, which puts the system into the form \eqref{flowmapcomplete02} with $\hat{F}=\xi$ and $y=J(\hat{u}+\varepsilon_a\mu)$. To design the target hybrid decision-making dynamics \eqref{hybrid_learning_dynamics}, and inspired by \cite{ObstacleAvoidance}, we let $q\in \mathcal{Q}=\{1,2\}$ be a logic state that remains constant during flows and which is used to construct the following mode-dependent function: 
\begin{equation}\label{hatJq}
\hat{J}_{q}(\hat{u})=-J(\hat{u})+B_{q}(\hat{u}),
\end{equation}
where $B_{q}$ is a mode-dependent barrier function to be designed. We restrict our attention to bounded obstacles $\mathcal{N}\subset\mathbb{R}^n$ for which there exists $p_0=(p_{0,1},p_{0,2})\in\mathbb{R}^2$, $\rho\in\mathbb{R}_{>0}$ and $\delta\in\mathbb{R}_{>0}$ such that $\mathcal{N}\subset p_0+\rho\mathbb{B}$ and $(p_0+2\rho\sqrt{2}\mathbb{B})\cap (u^{\star}+\delta\mathbb{B})=\emptyset$, that is, the obstacle $\mathcal{N}$ is contained in a ball of radius $\rho$, centered at the point $p_0$, and located sufficiently far away from the target $p^{\star}$. To achieve obstacle avoidance, we divide the space into two different regions, each region indexed by a mode $q\in\{1,2\}$, and having a corresponding mode-dependent signal $\hat{J}_{q}(\cdot)$, defined in \eqref{hatJq}. The main goal is to design the regions so that the vehicle can synergistically use the two potential fields $\hat{J}_{q}(\cdot)$ to safely navigate the space. To formalize this idea, consider the set 
\begin{equation}\label{ball_obstacle}
\mathcal{B}_{p_0,\rho}:=\left\{p:~\|p-p_0\|_{1}\leq 2\rho\sqrt{2}\right\},
\end{equation}
which satisfies $\{p_0\}+\rho\mathbb{B}\subset\mathcal{B}_{p_0,\rho}\subset \{p_0\}+2\rho\sqrt{2}\mathbb{B}$. Additionally, consider the following subsets of $\mathbb{R}^2:$
\begin{align*}
L_{1a}&:=\left\{p:~p_2< -p_1+p_{0,2}+p_{0,1}-2\rho\sqrt{2}\right\},\\
L_{1b}&:=\left\{p:~p_2< p_1+p_{0,2}+p_{0,1}-2\rho\sqrt{2}\right\},\\
L_{2a}&:=\left\{p:~p_2> p_1+p_{0,2}+p_{0,1}+2\rho\sqrt{2}\right\},\\
L_{2b}&:=\left\{p:~p_2> -p_1+p_{0,2}+p_{0,1}+2\rho\sqrt{2}\right\},
\end{align*}
and let $\mathcal{L}_1:=L_{1a}\cup L_{1b}$, $\mathcal{L}_2:=L_{2a}\cup L_{2b}$, $\mathcal{U}:=\mathcal{L}_1\cup\mathcal{L}_2$. Figure \ref{fig11203} illustrates the construction of the sets $\mathcal{L}_1$ and $\mathcal{L}_2$, which satisfy $u^{\star}\in \mathcal{L}_1\cap\mathcal{L}_2$, and also $\mathcal{L}_1\cap\mathcal{N}=\emptyset$, $\mathcal{L}_2\cap\mathcal{N}=\emptyset$. In fact, $\mathcal{U}=\mathbb{R}^2\backslash\mathcal{B}_{p_0,\rho}$. The mode-dependent barrier function $B_{q}$ is designed such that the function $W(\hat{u},q):=\hat{J}_{q}(\hat{u})+J(\hat{u}^*)$ satisfies the following assumption:
\begin{assumption}[Mode-Dependent Barrier Functions]\label{barrierassumption}
\begin{enumerate}[(a)]
\item For each $q\in \mathcal{Q}$, the map $\hat{J}_{q}:\mathbb{R}^2\to\mathbb{R}_{\geq0}\cup\{\infty\}$ is continuously differentiable in $\mathcal{L}_{q}$, and as $\hat{u}\to\infty$ or $\hat{u}\to \text{bd}(\mathcal{L}_{q})$ we have that $\hat{J}_{q}(\hat{u})\to\infty$. Moreover, for every $\hat{u}\in\mathbb{R}^2\backslash\mathcal{L}_{q}$, we define $J_{q}(\hat{u}):=\infty$.
\item There exist functions $\alpha_{1},\alpha_{2}\in\mathcal{K}_{\infty}$, and a proper indicator\footnote{A function $\tilde{\varpi}:\mathcal{U}\to\mathbb{R}_{\geq0}$ is a proper indicator on the open set $\mathcal{U}$ if it is continuous and $\tilde{\varpi}(x_i)\to\infty$ when $i\to\infty$ if either $|x_i|\to\infty$ or the sequence $\{x_i\}_{i=1}^{\infty}$ approaches the boundary of $\mathcal{U}$.} $\tilde{\varpi}$ of $u^{\star}$ on $\mathcal{U}$, such that
\begin{equation*}
\alpha_{1}(\tilde{\varpi}(\hat{u}))\leq \min_{s\in \mathcal{Q}}W(\hat{u},s)\leq\alpha_{2}(\tilde{\varpi}(\hat{u})),~~\forall~\hat{u}\in \mathcal{L}.
\end{equation*}
\item There exists a positive definite function $\rho$, such that for each $q\in \mathcal{Q}$:
\begin{equation*}
|\nabla \hat{J}_{q}(\hat{u})|^2 \geq \rho(W(\hat{u},q)),
\end{equation*}
for all $\hat{u}\in \mathcal{L}_{q}$. \QEDB
\end{enumerate}
\end{assumption}
The above properties essentially guarantee that in each region $\mathcal{L}_{q}$, a gradient-based feedback law of the form $\dot{\hat{u}}=-\nabla \hat{J}_{q}(\hat{u})$ can steer the vehicle towards points where the virtual potential $\hat{J}_{q}$ has smaller values, without leaving the region $\mathcal{L}_{q}$. Since it is not possible to satisfy property (c) in the entire space using only one barrier function, we ask that the condition holds only in each of the regions $\mathcal{L}_{q}$ via a different potential $\hat{J}_{q}$.

By using the measurements of $\hat{J}_{q}$ at the points $\hat{u}$, we can now introduce the flow and jump sets for the hybrid set-seeking system, with state $(u,q)$:
\begin{subequations}\label{flowobstacleset}
\begin{align} 
C_{u,q}&:=\left\{(\hat{u},q)\in\overline{\mathcal{U}}\times \mathcal{Q}:\hat{J}_{q}(\hat{u})\leq \chi \hat{J}_{3-q}(\hat{u})\right\}, \\
D_{u,q}&:=\left\{(\hat{u},q)\in \overline{\mathcal{U}}\times \mathcal{Q}:\hat{J}_{q}(\hat{u}) \geq (\chi-\lambda)\hat{J}_{3-q}(\hat{u})\right\},
\end{align}
\end{subequations}
where $\chi\in(1,\infty)$ and $\lambda\in(0,\chi-1)$ are tunable parameters. The set $C_{u,q}$ describes the points in the space where the controller uses the potential field $\hat{J}_{q}$, while the set $D_{u,q}$ indicates the "switching zones" for the controller, where the controller toggles the logic state using $q^+=3-q$. By construction, this switching behavior will take place whenever the value of the current function $\hat{J}_{q}$ exceeds a threshold compared to the other function $\hat{J}_{3-q}$. In particular, the construction of the flow and jump sets imposes a \emph{hysteresis} property on the switching controller, which precludes Zeno behavior. The red and blue lines in Figure \ref{fig11203} illustrate the boundaries of $D_{u,q}$ and  $C_{u,q}$ respectively. The following proposition shows that the resulting nominal target hybrid controller has the form \eqref{hybrid_learning_dynamics} and satisfies all the properties of Assumption \ref{keyassumptions} relative to $\overline{\mathcal{U}}\times \mathcal{Q}$.

\begin{proposition}[Source-Seeking $\&$ Obstacle Avoidance]
Consider the hybrid set-seeking system \eqref{originalsystemcomplete01}, with sets $C_{u,q}$ and $D_{u,q}$ given by \eqref{flowobstacleset}, function $\hat{F}_{\delta}=\xi$, and input $y=\hat{J}_{q}$. Suppose that Assumption \ref{barrierassumption} holds. Then, system \eqref{originalsystemcomplete01} satisfies item (a) in Assumption \ref{keyassumptions} with $r=1$, $\delta=0$, $\hat{z}=q$, $\Psi=\mathcal{Q}$, $\mathcal{O}=\{u^*\}$. Moreover, item (b) holds with basin of attraction $\mathcal{U}\times \mathcal{Q}$.    \QEDB 
\end{proposition}
To illustrate the stability properties of the resulting hybrid source-seeking system for obstacle avoidance, we show in Figure \ref{fig11203} the constructions of the partitions $\mathcal{L}_1$ and $\mathcal{L}_2$, as well as the trajectory of the vehicle \eqref{plant_obstacle} under the control law \eqref{hybridcontroller01}, evolving over the two virtual partitions. As observed, the hybrid controller can (practically) stabilize the source of the signal $J$ by switching the potential field after approximately 30s of flow, while avoiding the obstacle. We note that this simple example can be extended to other vehicle dynamics, including non-holonomic vehicles \cite{Poveda:20TAC}, and to the multi-obstacle case by considering other types of target hybrid systems \cite{Poveda:20TAC,berkane2021obstacle,sawant2024hybrid}. In such cases, the target hybrid decision-making dynamics \eqref{hybrid_learning_dynamics} will be modified, but the structure and stability properties of the hybrid set-seeking dynamics \eqref{originalsystemcomplete01} can be studied using similar tools.
%
% \begin{figure}[t!]
%     \centering \includegraphics[width=0.92\columnwidth]{figs/SchemeUnicycle.pdf}
%     \caption{Scheme of the source-seeking controller under sporadic measurements and persistent spoofing.} \label{fig:sourceSeekingScheme}
% \end{figure}
%
\begin{figure}[t!]
 \begin{tcolorbox}[colback=ivoryA, colframe=ivoryA]
\centering
\includegraphics[width=1\linewidth]{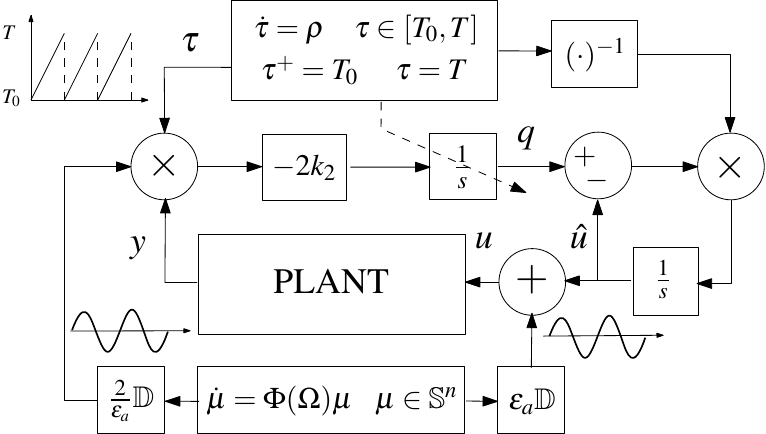}
\end{tcolorbox}
\caption{Block diagram representation of the hybrid set-seeking system with momentum, with  \eqref{HDS0121} as target hybrid decision-making dynamics with $k_1=0$. Here, the low-pass filter is omitted. \label{fig11203asd}}
\end{figure}
\subsection{Seeking with Momentum and Resets}
Reset control \cite{nonconvexopti,nesic2011stability,zhao2019overcoming} is a well-known technique used to improve transient performance in dynamical systems that incorporate integrators in the feedback loop. The main idea behind the incorporation of resets is to dissipate energy at appropriate times to decrease overshoots and improve transient performance \cite{zhao2019overcoming}. Resets can also be used to model restarting methods \cite{Candes_Restarting,ODE_Nesterov}, commonly employed in optimization routines. By leveraging tools from hybrid systems, set-seeking algorithms with momentum and restarting can also be designed to reduce overshoots and improve transient performance \cite{PovedaNaliAuto20,TAC21Momentum_Nash,GalarzaPovedaDallanese}. To illustrate this idea, we can consider the following target hybrid decision-making dynamics with state $(\hat{u},q,\tau)\in\mathbb{R}^n\times\mathbb{R}^n\times\mathbb{R}_{\geq0}$ and data
\begin{subequations}\label{HDS0121}
\begin{align}
C_{u,z}&=\mathbb{R}^n\times\mathbb{R}^n\times[T_0,T]\label{flowsetmomentum}\\
\hat{F}(\hat{u},q,\tau)&=\left(\begin{array}{c}
\frac{2}{\tau}(q-\hat{u})-k_1\nabla J(\hat{u})\\
-2k_2\tau\nabla J(\hat{u})\\
\rho
\end{array}\right),\label{flowmomentum1}\\
D_{u,z}&=\mathbb{R}^n\times\mathbb{R}^n\times\{T\}\label{jumpsetmomentum}\\
\hat{G}(\hat{u},q,\tau)&=\left(\begin{array}{c}
\hat{u}\\
\alpha\hat{u} +(1-\alpha)q\\
T_0
\end{array}\right),\label{discretemomentum}
\end{align}
\end{subequations}
where $k_1\geq0$ and $T_0,T,k_2,\rho>0$ are tunable parameters that satisfy $T>T_0$. The parameter $\alpha\in\{0,1\}$ is selected \emph{a priori} to indicate whether the state $q$ is reset during jumps. In particular, when the timer $\tau$ satisfies $\tau=T$, the discrete-time dynamics \eqref{discretemomentum} reset the timer back to $T_0$. If, additionally, the parameter $\alpha$ satisfies $\alpha=1$, then the state $q$ is also reset to the current value of $\hat{u}$. Figure \ref{fig11203} presents a block-diagram representation of the hybrid set-seeking dynamics  \eqref{originalsystemcomplete01} associated with \ref{HDS0121} for the case $k_1 = 0$. The filter is omitted for simplicity.

Using the change of coordinates $s=\hat{u}$ and $q=s+\frac{1}{2}\tau(\dot{s}+k_1\nabla J)$, the flows of the above system are related to the second-order ordinary differential equation
\begin{equation*}
\ddot{s}+\frac{(2+\dot{\tau})}{\tau}\dot{s}+4k_2\nabla J(s)+k_1\left(\nabla^2 J(s)^\top \dot{s}+\frac{\dot{\tau}}{\tau}\nabla J(s)\right)=0,
\end{equation*}
which has been studied in the context of accelerated optimization algorithms under different choices of $\dot{\tau}$ and parameters $k_1,k_2$ \cite{ODE_Nesterov,zhang2018direct,michalowsky2014multidimensional,laborde2020lyapunov,Poveda_Li:2019_CDC}. For example, when $k_1=0$ the above dynamics can also be written in state-space representation using $s_1=s$ and $s_2=\dot{s}$, leading to
\begin{equation}\label{ODESSDSF}
\dot{s}_1=s_2,~~\dot{s}_2=-\frac{c_1}{\tau}s_2-c_2\nabla J(s_1),~~\dot{\tau}=\rho,
\end{equation}
for suitable constants $c_1,c_2>0$. Note that when $\nabla J=\xi$, the system \eqref{ODESSDSF} has a structure similar to \eqref{flow_learninghf001} when $\hat{F}_\delta=\xi$. To avoid the issues discussed in Example \ref{example9}, we can consider the flow and jump sets described in \eqref{jumpsetmomentum} and \eqref{flowsetmomentum}, respectively,  and the jump map
\begin{equation}\label{JUMPSSSDSF}
s_1^+=s_1,~~~s_2^+=(1-\alpha) s_2,~~\tau^+=T_0,
\end{equation}
which has a similar structure as \eqref{jump_learninghf001}, with $s_1=\hat{u}$, $\hat{z}=\tau$, $s_2=\xi$, and $\hat{G}_{\delta}=(s_1,T_0)$ and $G_{\xi}=(1-\alpha)s_2$. Therefore, this system can be viewed as a target hybrid decision-making system equipped with a hybrid filter that incorporates jumps to periodically reset the momentum to zero.  Under a suitable choice of parameters, this hybrid set-seeking system can stabilize the set of minimizers of a convex functions $J$ using only measurements. In particular, the following proposition follows from the results in \cite{PovedaNaliAuto20}, and illustrates the role of Theorem \ref{theorem4} in a slightly modified version of system \eqref{hybrid_learning_dynamics001} with nonlinear $\tau$-dependent damping.
\begin{figure}[t]
\begin{tcolorbox}[colback=iceblue!80, colframe=iceblue!80]
\centering
\includegraphics[width=\linewidth]{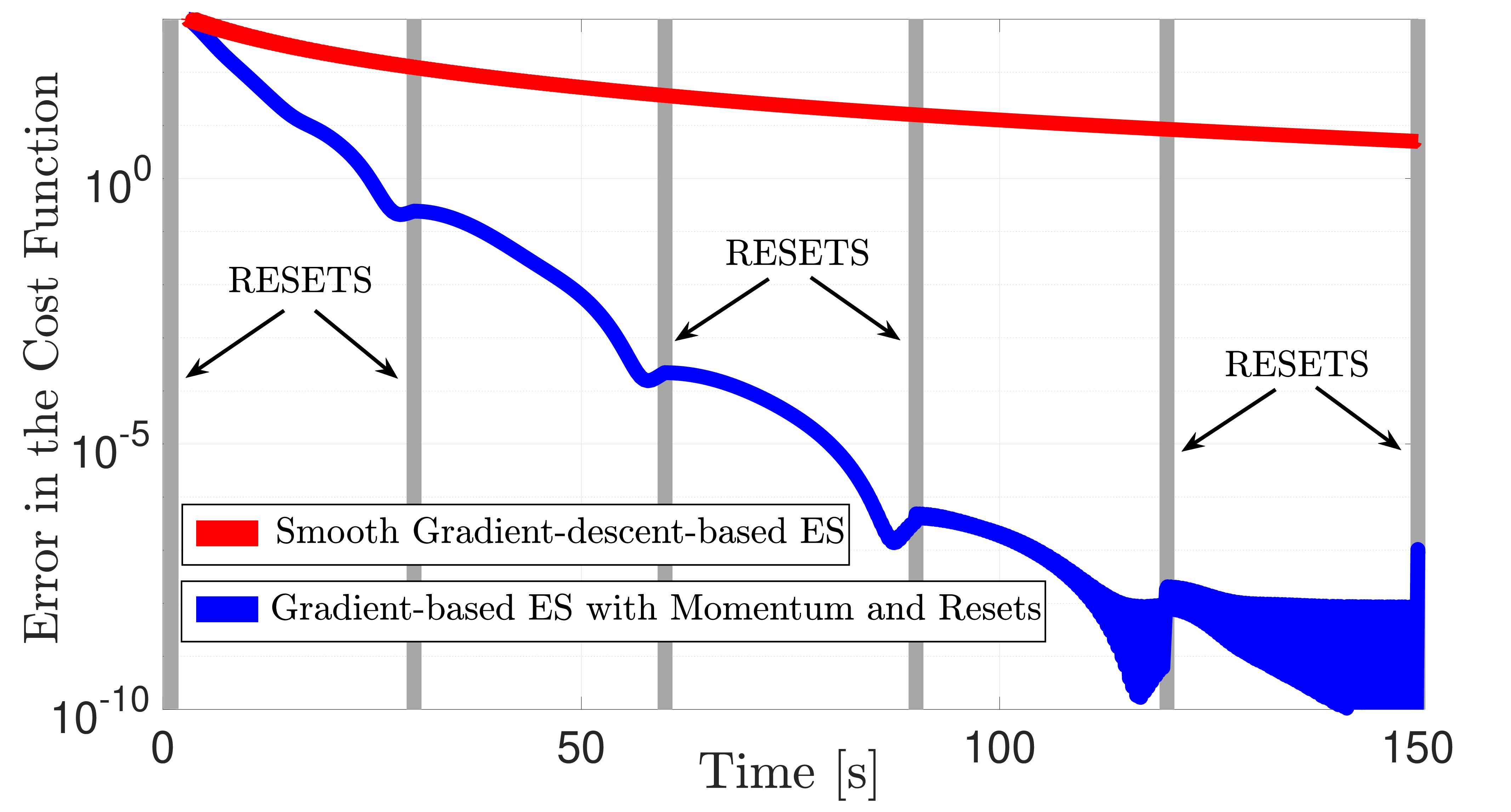}
\end{tcolorbox}
\caption{Trajectories $\hat{u}$ generated by hybrid set-seeking dynamics with momentum and resets, compared to those of \eqref{ESCODEVanilla}. \label{fig112035}}
\end{figure}
\begin{proposition}[Seeking with Momentum and Resets]\label{propositionmomentum}
Suppose that $J$ is strictly convex and its unique minimizer is given by $u^*$, $\nabla J$ is globally Lipschitz and $\rho=\frac{1}{2}$. If $\alpha=0$, then the HDS with flow and jump sets given by \eqref{flowsetmomentum} and \eqref{jumpsetmomentum}, and flow and jump maps given by \eqref{ODESSDSF} and \eqref{JUMPSSSDSF}, satisfies all the items of Assumption \ref{assumption2} with $\hat{u}=s_1$, $\hat{z}=\tau$, $r=1$, $\xi=s_2$, $\kappa=-1$, $C_{u,z}=\mathbb{R}^n\times[T_0,T]$, $D_{u,z}=\mathbb{R}^n\times\{T\}$, $\Lambda_{\xi,c}=\Lambda_{\xi,d}=\mathbb{R}^n$, $\mathcal{O}=\{u^*\}$, $\Psi_z=[T_0,T]$, and $\Psi_f=\{0\}$. If $\alpha=1$, and $J$ is also strongly convex, then Assumption \ref{assumption2} also holds for $T$ sufficiently large and in this case the set $\mathcal{A}$ is UGES.
\end{proposition}
Figure \ref{fig112035} shows the performance of the hybrid set-seeking constructed with \eqref{HDS0121} as the target decision-making system, compared to the standard smooth gradient descent-based ES algorithm \eqref{ESCODEVanilla}. As observed, the transient performance of the controller can be substantially improved under a suitable tuning of the resetting parameters $[T_0,T]$. For applications of similar ES algorithms in the context of traffic light systems, we refer the reader to \cite{ESC_Traffic}.

The hybrid set-seeking dynamics \eqref{HDS0121} can be extended to multi-agent systems, as well as game theoretic settings where the goal is to converge to a Nash equilibrium \cite{TAC21Momentum_Nash,krilavsevic2023learning}. However, in these applications, the resets need to be coordinated to emulate the centralized "accelerated" performance observed in Figure \ref{fig112035}, see \cite{TAC21Momentum_Nash}. Typically, such coordination mechanisms also need to be set-valued to guarantee global convergence properties \cite{SyncPovedaTeelAutomatica}. Adaptive resetting mechanisms that implement resets whenever state-dependent conditions are satisfied can also be considered for the design of hybrid set-seeking systems, see \cite{teel2019first}.   
\subsection{Seeking with Aggressive Gradient Exploration and Newton Fine-Tuning}
When the cost function \( J \) is smooth and its Hessian matrix is positive definite, the decision-making problem \eqref{optimization_problem} can be efficiently solved using Newton-like ES algorithms, which offer superior transient performance---allowing user-specified convergence rates---compared to traditional gradient-based ES methods \cite{MISONewton,GalarzaPovedaDallanese}. However, real-time estimation and inversion of the Hessian matrix in Newton-like ES algorithms can be highly sensitive to measurement noise and may lead to instabilities, particularly when the initial conditions are far from the minimizers of $J$. While increasing the frequency of the dither signals can help mitigate these challenges, it may require potentially unfeasible shorter sampling intervals to avoid aliasing during implementations. One potential solution to this issue is to use the more stable gradient-based ES algorithm when the system is far from the optimizer, and to switch to a Newton-like ES method for fine-tuning convergence near the optimum, a methodology that is not uncommon in standard optimization \cite[pp. 23]{martens2012training}. This idea can be formalized and implemented using hybrid set-seeking dynamics whenever a reasonable estimate of the optimal value of the cost $J$ is known. To simplify our presentation, we focus on the 2-dimensional setting, and we incorporate a Hessian estimator into the flow map of \eqref{flowmapcomplete}. Namely, instead of \eqref{flowmapcomplete}, we let $\xi=(\xi_1,\xi_2)$ and consider a hybrid set-seeking system with feedback law \eqref{main_input}, overall state $(x_{u,z},\xi,\mu)$, and flow map
\begin{figure}[t!]
 \begin{tcolorbox}[colback=ivoryA, colframe=ivoryA]
\centering
\includegraphics[width=1\linewidth]{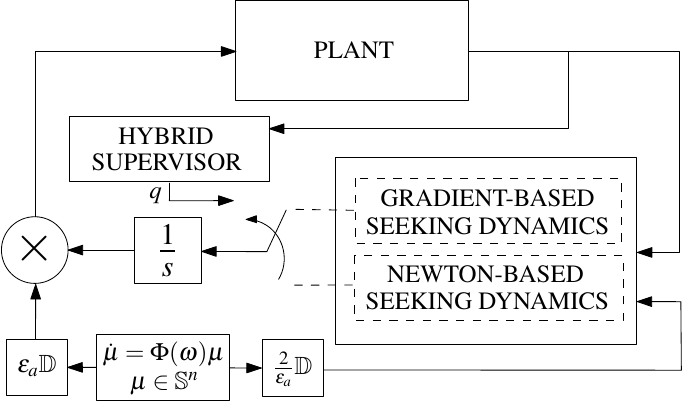}
\end{tcolorbox}
\caption{Scheme of hybrid set-seeking system that adaptively switches between gradient-based ES and Newton-based ES. \label{fig11202345678}.}
\end{figure}

\begin{figure*}[t!]
\begin{tcolorbox}[colback=iceblue!80, colframe=iceblue!80]
\centering
\includegraphics[width=0.33\linewidth]{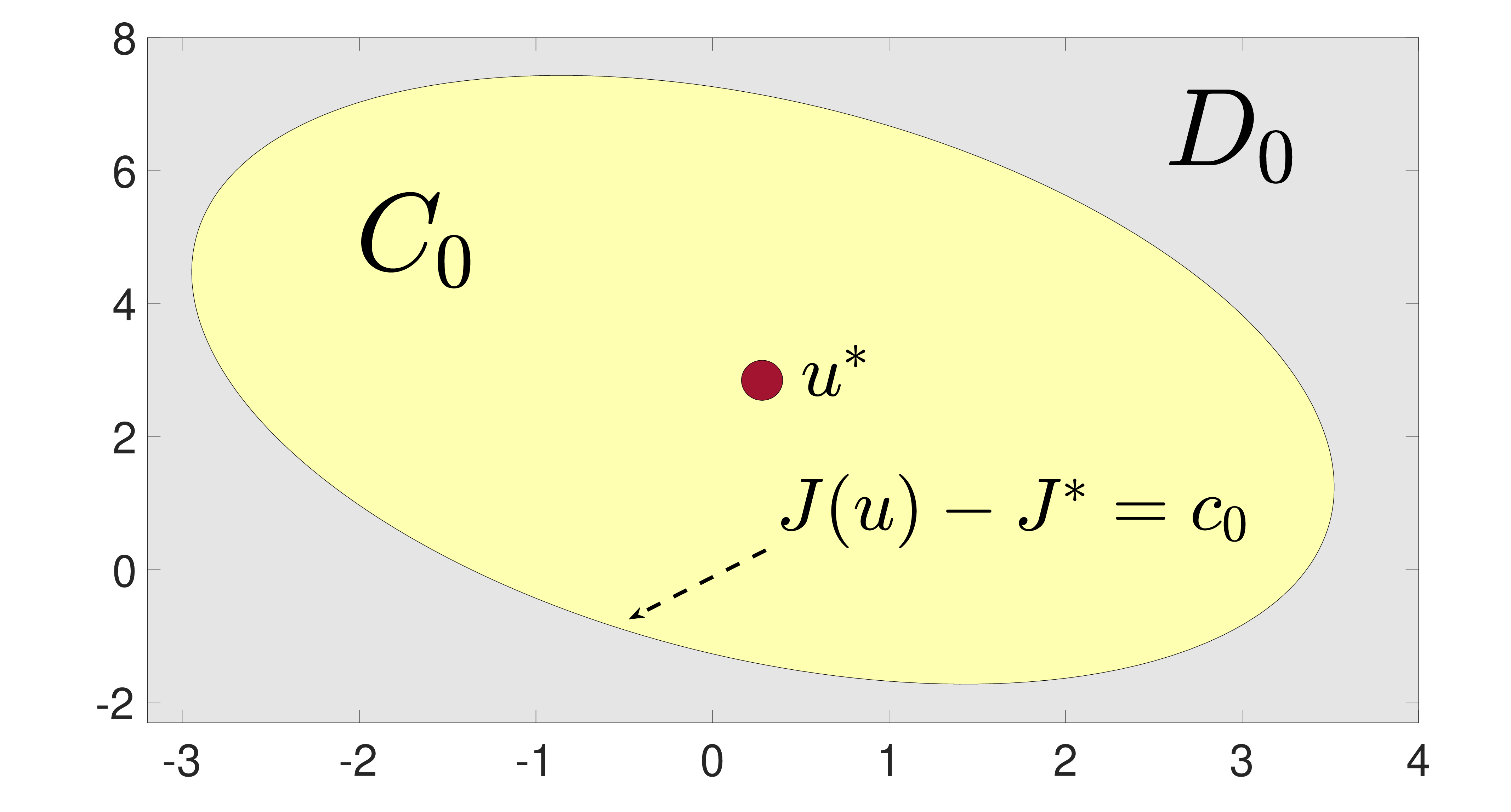}\includegraphics[width=0.33\linewidth]{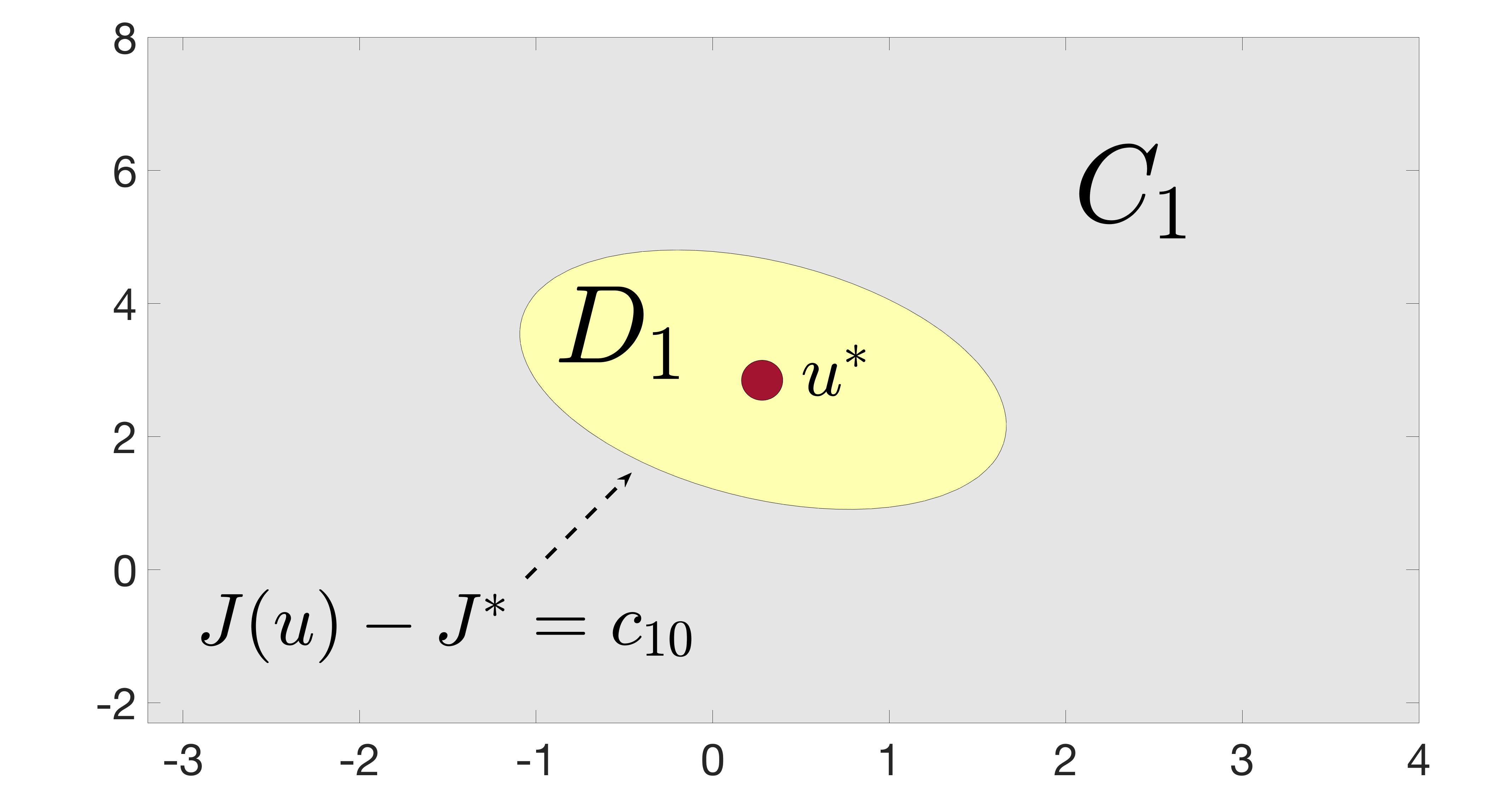}
\hspace{-0.15cm}\includegraphics[width=0.33\linewidth]{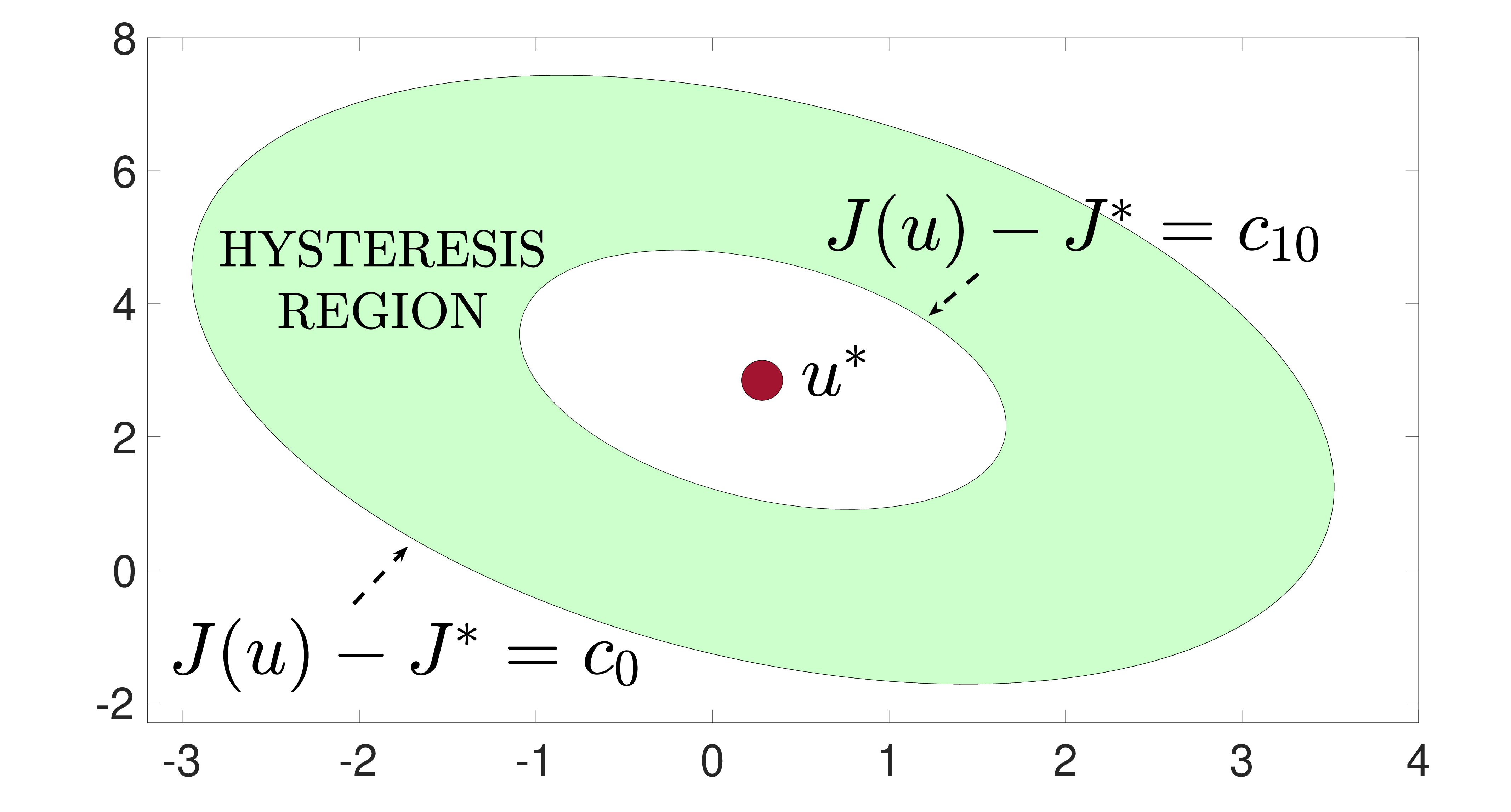}
\end{tcolorbox}
\caption{\tcr{Representation of sets defined in \eqref{setsflowqs} and \eqref{setsjumpqs}, as well as the region of the space where hysteresis is induced in the switched controller. \label{fig11202345678}}}
\end{figure*}

\begin{align}\label{flowmapng}
 \dot{x}\in F(x):&= \left(  \begin{array}{c}
\hat{F}(x_{u,z},\xi)\\
 -\frac{k_f}{\varepsilon_f}\left(\xi_1-\frac{2}{\varepsilon_a}y\mathbb{D}\mu\right)\\
-\frac{k_h}{\varepsilon_f}\left(y N(\mu)\xi_2-y\mathbb{D}\mu \right)\\
\Phi(\omega)\mu
\end{array}\right),
\end{align}
where $x_{u,z}=(\hat{u},\hat{z})\in\mathbb{R}^3$, and
\begin{equation}\label{mapaux11}
\hat{F}(x_{u,z},\xi)=\left(\begin{array}{c}
-\hat{z}\xi_1-(1-\hat{z})\xi_2\\
0
\end{array}\right).
\end{equation}
In this system, the auxiliary state $\hat{z}\in\mathcal{Q}:=\{0,1\}$ acts as a logic mode, and $\xi_2\in\mathbb{R}^{2}$ can be seen as an estimate of the vector  $\nabla^2 J^{-1}\nabla J$, where $\nabla^2 J$ is the Hessian of $J$ \cite{attouch2001second,labar2019newton}. The matrix-valued function  $N:\mathbb{R}^2\to\mathbb{R}^{2\times 2}$ can be designed as in \cite{MISONewton}, so that it has entries satisfying $N_{11}=\frac{16}{\varepsilon_{a}^2}(\mu_1^2 -\frac{1}{2})$, $N_{22}=\frac{16}{\varepsilon_{a}^2}(\mu_3^2 -\frac{1}{2})$, $N_{12}=\frac{4}{\varepsilon_{a}^2}\mu_1\mu_3$, and $N_{12}=N_{21}$. 
To define the sets $C$ and $D$, we follow a uniting control approach \cite{RSanfeliceBook}, where the logic state $\hat{z}$ indicates the current algorithm used by the controller. In particular, when the value of the cost $J$ is sufficiently close to the optimal value $J^*$, we use $\hat{z}=0$. On the other hand, when the value of $J$ is greater than $J^*$, we use $\hat{z}=1$. While this approach might require knowledge of $J^*$, for many ES applications where the cost to be minimized corresponds to some error function, it is known that $J^*=0$. On the other hand, if $J^*$ is unknown, a rough estimate can be used for the initial construction of the controller.

Building on these ideas, the switching zones for $\hat{z}$ are defined using
tunable thresholds specified by constants $c_0$ and $c_{10}$, with
$c_0 > c_{10} > 0$. To characterize these thresholds, we define the sets
\begin{subequations}\label{setsflowqs}
\begin{align}
C_0&:=\{\hat{u}\in\mathbb{R}^2:J(\hat{u})-J^*\le c_0\},\\
C_1&:=\overline{\mathbb{R}^2\backslash \{\hat{u}\in\mathbb{R}^2:J(\hat{u})-J^*< c_{10}\}},
\end{align}
\end{subequations}
as well as the sets
\begin{subequations}\label{setsjumpqs}
\begin{align}
D_0&:=\overline{\mathbb{R}^2\backslash \{\hat{u}\in\mathbb{R}^2:J(\hat{u})-J^*< c_0\}},\\
D_1&:=\left\{\hat{u}\in\mathbb{R}^2:J(\hat{u})-J^*\le c_{10}\right\}.
\end{align}
\end{subequations}
Using the constructions \eqref{setsflowqs} and \eqref{setsjumpqs}, the flow and jump sets for the state $x_{u,z}$ are given by
\begin{subequations}\label{setsswitchednewton}
\begin{align}
C_{u,z}&:=\left(C_0\times\{0\}\right)\cup \left(C_1\times\{1\}\right)\\
D_{u,z}&:=\left(D_0\times\{0\}\right)\cup \left(D_1\times\{1\}\right).
\end{align}
\end{subequations}
and the discrete-time dynamics of $x_{u,z}$ are given by
\begin{equation}\label{mapaux11223}
\hat{G}(x_{u,z})=\left(\begin{array}{c}
\hat{u}\\
1-\hat{z}
\end{array}\right).
\end{equation}
This scheme can be seen as using a "hybrid supervisor" to decide when to switch, while precluding Zeno behavior owing to the hysteresis mechanism induced by the design of the flow and jump sets. A cartoon of the switching zones in $\mathbb{R}^2$ is shown in the left plot of Figure \ref{fig11202345678}. 

\begin{figure}[b]
\begin{tcolorbox}[colback=iceblue!80, colframe=iceblue!80]
\centering
\includegraphics[width=\linewidth]{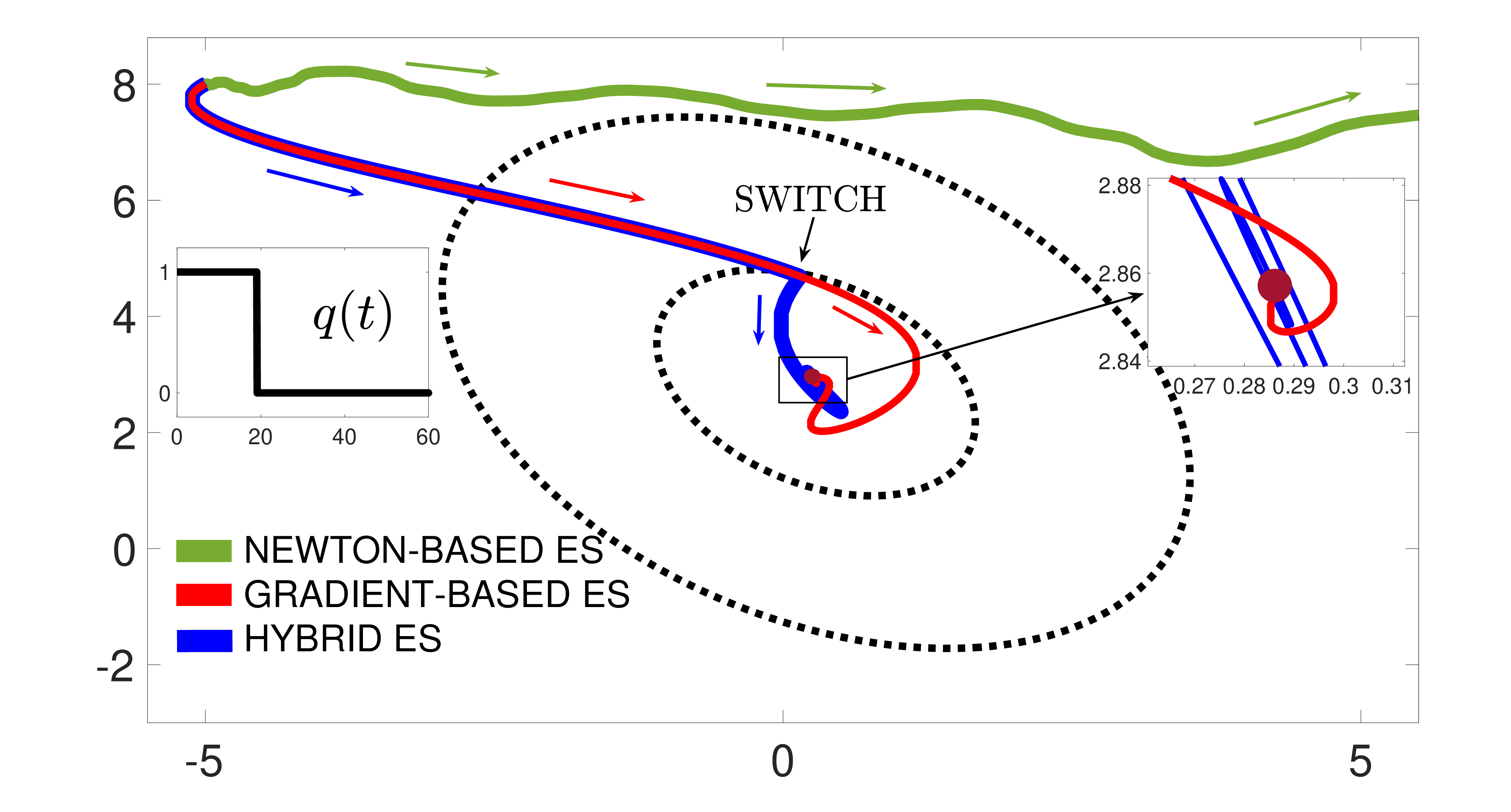}
\end{tcolorbox}
\caption{\tcr{Trajectory $\hat{u}$ of hybrid set-seeking system (shown in color blue) switching between Newton-based and gradient-based ES, compared to the trajectories obtained via standard gradient-based ES (shown in color red) and standard Newton-based ES (shown in color green).} \label{fig1120a}}
\end{figure}

\begin{figure}[t!]
\begin{tcolorbox}[colback=iceblue!80, colframe=iceblue!80]
\centering
\includegraphics[width=\linewidth]{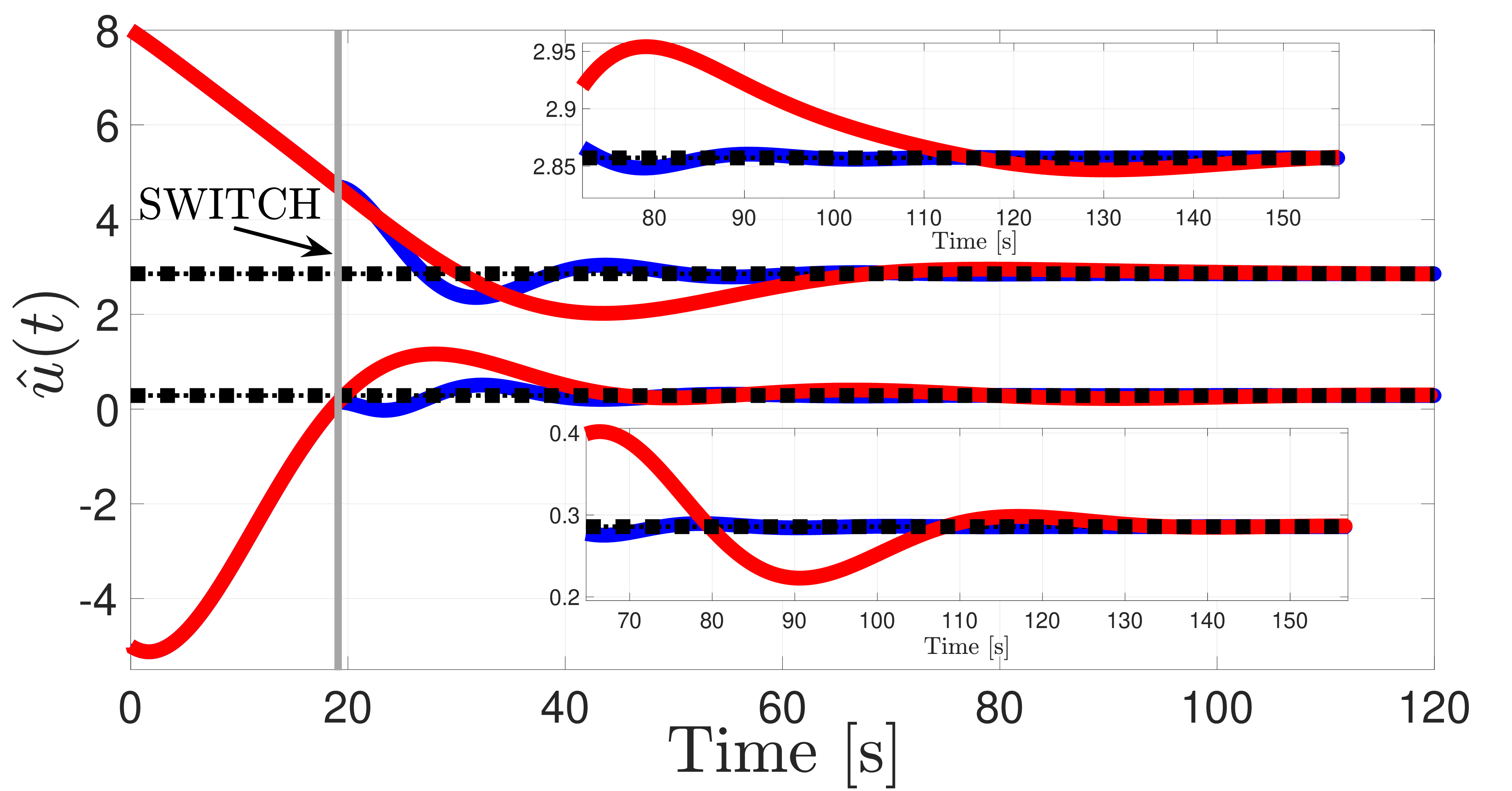}
\end{tcolorbox}
\caption{\tcr{Evolution in time of trajectories $\hat{u}$ for the hybrid set-seeking system switching between Newton-based and Gradient-based ES (shown in color blue). The red trajectory corresponds to the standard Gradient-based ES.}  \label{fig1120bced}}
\end{figure}

\begin{proposition}[Switched Newton-Gradient Dynamics]\label{proposition1}
Consider the hybrid set-seeking dynamics \eqref{originalsystemcomplete01}, with $y=J(u)$, flow map \eqref{flowmapng}, data \eqref{mapaux11}-\eqref{mapaux11223}, and feedback law \eqref{main_input}. Then, item (a) of Assumption 1 holds with $r=1$ and $\delta=0$. Moreover, if $J$ is smooth, radially unbounded and strongly convex, then item (b) also holds for the respective target hybrid decision-making dynamics with $\xi_1=\nabla J$, $\xi_2=\nabla^2J^{-1}\nabla J$, $\mathcal{O}=\{u^*\}$ and $\Psi=\{0\}$. \QEDB
\end{proposition}
We illustrate the Newton-Gradient Switched seeking dynamics using a simple two-dimensional numerical example for a quadratic cost function $J$, with the results shown in Figure \ref{fig1120a}. The figure depicts in color green a diverging trajectory that occurs when Newton-like ES dynamics are applied far from the optimal point using a moderate frequency $\omega$. 
As shown by the blue trajectory, incorporating the hybrid uniting mechanism introduced in this section not only eliminates the instability but also significantly improves the transient performance of the controller---outperforming the more robust yet slower gradient-based ES dynamics, represented by the red trajectory. Figure \ref{fig1120bced} highlights this improvement in transient behavior, reducing the convergence time by more than half.

\begin{pullquote}
Hybrid set-seeking systems can incorporate momentum with scheduled or adaptive restarting to accelerate convergence and reduce overshoots, while retaining robustness
\end{pullquote}
\begin{figure}[t!]
\begin{tcolorbox}[colback=ivoryA, colframe=ivoryA]
\centering
\includegraphics[width=0.95\linewidth]{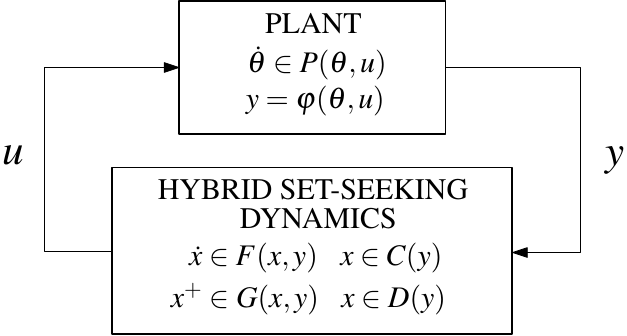}
\end{tcolorbox}
\caption{Dynamic plant interconnected with hybrid set-seeking dynamics for model-free feedback optimization.\label{fig112f}}
\end{figure}
% \begin{figure}[t!]
% \centerline{\includegraphics[width=0.99\linewidth]{instability.jpg}}
% \caption{Trajectories of the hybrid ES switching between Newton and gradient-based ES. \label{fig1120123}}
% \end{figure}
%\vspace{5cm}
%%%%%%%%%%%%%%%%%%%%%

% \begin{remark}[High-derivative Estimators]
% Dynamic Hessian estimators for the implementation of Newton-based algorithms could also be considered in (\ref{estimator_dynamics_1}). This is a relevant case since switching between gradient and Newton-based algorithms is a common approach to improve the convergence rate during the optimization process. 
% \end{remark}
% %
% \begin{remark}\label{remark_MAS2}
% \tcbb{In the multi-agent systems case (see also Remark \ref{Remark:MAS_1}), the output $y$ in (\ref{estimator_dynamics_1}) will correspond to a vector with entries given by the outputs of each of the individual subsystems, and the product $y\cdot \mathbb{D}\mu$ will correspond to the Hadamard product. In this case every agent will implement its own individual oscillator with matrix (\ref{individual_oscillator}).}
% \end{remark}

% \begin{remark}\label{remark_MAS3}
% \tcbb{For gradient-based Nash seeking schemes for MAS, Assumptions \ref{BigAssumption}-\ref{BigAssumption2} must hold replacing $\nabla J$ by the vector of individual partial derivatives $\left[\frac{\partial J_1(\hat{u})}{\partial \hat{u}_1},\ldots,\frac{\partial J_N(\hat{u})}{\partial \hat{u}_N}\right]^\top.$}
% \end{remark} 

\begin{sidebar}{Singularly Perturbed Hybrid Dynamical Systems}\label{singularSPHDS}
\sdbarinitial{S}tability properties of extremum seeking systems rely on multiple time-scale separations induced in the closed-loop system via appropriate tuning of the control parameters. Singular perturbation theory is a well-established field in the literature of differential equations that enables the study of systems for which certain states evolve in a faster time scale and converge to a quasi-steady state manifold, where the original system can be "reduced" \cite{Kokotovic_SP_Book,KhalilBook}. Here, we review some extensions of singular perturbations to hybrid systems \cite{SAA12,Wang:12_Automatica,abdelgalil2023multi} that fit the models considered in this paper. Further details can be found in \cite{SAA12}.

\vspace{-0.3cm}
\subsection{Hybrid Systems with "Fast" and "Slow" States}
Consider a hybrid system with states $(x,\theta)\in\mathbb{R}^{n_1}\times\mathbb{R}^{n_2}$ and the following dynamics: 
\begin{subequations}\label{SPSHDS1}
\begin{align}
&(x,\theta)\in C:= C_x\times \Theta,~~~~\left\{\begin{array}{l}\dot{x}\in F_x(x,\theta)\\
\varepsilon\dot{\theta}\in F_{\theta}(x,\theta)
\end{array}\right.
\label{flow_dynamics1}\\
&(x,\theta)\in D:= D_x\times \Theta,~~~~(x^+,\theta^+)\in G(x,\theta),\label{jump_dynamics1}
\end{align}
\end{subequations}
where the set-valued mappings $F_x:\mathbb{R}^{n_1}\times\mathbb{R}^{n_2}\rightrightarrows\mathbb{R}^{n_1}$ and  $F_{\theta}:\mathbb{R}^{n_1}\times\mathbb{R}^{n_2}\rightrightarrows\mathbb{R}^{n_2}$ characterize the continuous-time dynamics of $x$ and $\theta$, respectively. Similarly,  $G:\mathbb{R}^{n_1}\times\mathbb{R}^{n_2}\rightrightarrows\mathbb{R}^{n_1+n_2}$ describes the jump map of the system. The flow set $C$ is comprised by the Cartesian product of the sets $C_x\subset\mathbb{R}^{n_1}$ and $C_z\subset\mathbb{R}^{n_2}$, while the jump set $D$ is comprised of the sets $D_x\subset \mathbb{R}^{n_1}$ and $D_{\theta}\subset\mathbb{R}^{n_2}$. We make the assumption that system \eqref{SPSHDS1} satisfies the hybrid basic conditions and that $\Theta$ is a compact set. 

In system \eqref{SPSHDS1}, $\varepsilon\in\mathbb{R}_{>0}$ is a small parameter that induces a time-scale separation between the continuous-time dynamics of the state $x\in\mathbb{R}^{n_1}$ and the continuous-time dynamics of $\theta\in\mathbb{R}^{n_2}$. In particular, when $\varepsilon$ is small, $x\in\mathbb{R}^{n_1}$ behaves as a ``slow'' state in the system, while $\theta\in\mathbb{R}^{n_2}$ behaves as a ``fast'' state. %We let $z=(x,\theta)\in\mathbb{R}^n$ denote the overall state of the system, where $n_1+n_2=n$. 
When the fast dynamics in system \eqref{flow_dynamics1} define a suitable "quasi-steady" state, the behavior of the state $x$ can be approximately predicted by a \emph{reduced} system. For example, this is the case when $\theta$ converges to a continuously differentiable quasi-steady state manifold $m:\mathbb{R}^n\to\mathbb{R}^m$, as is usually the case in the literature of differential equations. However, in the context of hybrid systems, it is common to have fast dynamics for which the quasi-steady state is actually characterized by a set-valued mapping $M:\mathbb{R}^n\rightrightarrows\mathbb{R}^m$. In such cases, we assume that $M$ is outer-semicontinuous and locally bounded, and that for each $\rho>0$ the set
\begin{equation}
M_{\rho}=\{(x,\theta): \theta\in M(x),~x\in \rho\mathbb{B}\}
\end{equation}
is UGAS for the following \emph{boundary layer} dynamics:
\begin{equation}\label{boundarylayerdynamics}
(x,\theta)\in \left(C_x\cap \rho\mathbb{B}\right)\times\Theta,~~~~\dot{x}=0,~~~\dot{\theta}\in F_{\theta}(x,\theta).
\end{equation}
Note that in \eqref{boundarylayerdynamics}, the state $x$ is taken as constant from any initial condition on $C_x\cap \rho\mathbb{B}$. 

Using the mapping $M$, the \emph{reduced hybrid system}, with state $x_r\in\mathbb{R}^n$, can be defined using the following dynamics:
\begin{subequations}\label{reducedHDS}
\begin{align}
&x\in C_x,~~\dot{x}\in F_r(x):=\overline{\text{co}}\{f:f\in F_x(x,\theta),\theta\in M(x)\}\\
&x\in D_x,~~x^+\in G_r(x):=\{f_1:(f_1,f_2)\in G_x(x,\theta),~\theta\in\Theta\}.
\end{align}
\end{subequations}
By construction, and by the assumptions on the data of system \eqref{SPSHDS1}, the reduced hybrid system \eqref{reducedHDS} satisfies the Hybrid Basic Conditions. 

\vspace{-0.3cm}
\begin{theorem}[SGPAS in Singularly Perturbed HDS]\label{theorem3SP}
Suppose there exists a compact set $\mathcal{A}\subset\mathbb{R}^{n_1}$ that is UGAS for the reduced HDS \eqref{reducedHDS}. Then, the singularly perturbed HDS \eqref{SPSHDS1} renders the set $\mathcal{A}\times\Theta$ SGPAS as $\varepsilon\to0^+$.
\end{theorem}

\vspace{-0.3cm}
When the hybrid system \eqref{SPSHDS1} possesses additional structure, stronger results such as UGAS or UGES can be established. For instance, in the context of hybrid extremum seeking with hybrid filters operating on a faster time scale, the flow and jump sets in \eqref{SPSHDS1} take the form \( C = C_x \times C_{\theta} \) and \( D = D_x \times D_{\theta} \), where \( C_{\theta} \) and \( D_{\theta} \) are not necessarily compact. By imposing appropriate Lyapunov-based conditions on the reduced dynamics, boundary layer dynamics, and certain interconnection properties, UGAS or UGES can usually be established for these dynamics, provided that \( \varepsilon \) is sufficiently small \cite[Sec. 5]{abdelgalil2023multi}. This approach mirrors the classical composite Lyapunov method based on quadratic-like Lyapunov functions used to analyze singularly perturbed differential equations \cite{Saberi}.

\vspace{-0.3cm}
\subsection{Remark 1: Averaging using Singular Perturbations}
There are close connections between averaging and certain models of singularly perturbed systems. For example, when the fast state in \eqref{flow_dynamics1} does not settle to a point, but instead oscillates, it is still possible to define a well-posed reduced hybrid system. For example, when the continuous-time dynamics in \eqref{flow_dynamics1} have the form
\begin{align}
(x,\theta)\in C:= C_x\times \Theta,~~~~\left\{\begin{array}{l}\dot{x}=f_x(x,\theta)\\
\varepsilon\dot{\theta}=f_{\theta}(x,\theta)
\end{array}\right.
\label{flow_dynamics1}
\end{align}
and the jump map in \eqref{jump_dynamics1} is of the form $G(x,\theta)=G_x(x)\times\{\theta\}$, we can define a hybrid reduced system using the data
\begin{subequations}\label{timeinvariantaveraging}
\begin{align}
\bar{x}\in C_x,~~\dot{\bar{x}}=\bar{f}(\bar{x}),\\
\bar{x}\in D_x,~~\bar{x}^+\in G(\bar{x}),
\end{align}
\end{subequations}
where $\bar{f}$ is a continuous function that satisfies
\begin{equation}\label{average0assumption01}
\left|\frac{1}{T}\int_{0}^{T} \Big(f_x(x,\theta(s))-\bar{f}(x,p)\Big)\text{d}s\right|\leq \gamma_K(T), 
\end{equation}
for all $\mathcal{T}\in\mathbb{R}_{\geq0}$, all $x\in C_x\cap K$, all compact sets $K\subset\mathbb{R}^{n_1}$, and all solutions $\theta:[0,T]\to\mathbb{R}^{n_2}$ of the dynamics $\dot{\theta}(t)=f_{\theta}(x,\theta(t))$, where $\gamma_K(\cdot)$ is a continuous, non-increasing function that satisfies $\lim_{\mathcal{T}\to\infty}\gamma_{K}(\mathcal{T})=0$, and which is allowed to depend on $K$. Note that the definition of $\bar{f}$ in \eqref{average0assumption01} is similar to the average map defined in \eqref{average0assumption}. Indeed, one can study oscillatory hybrid systems with fast states generated by time-invariant systems via averaging using the model \eqref{timeinvariantaveraging}. This observation is important in order to exploit robustness properties in well-posed hybrid systems with compact attractors.

\end{sidebar}
\section{Part 3: Hybrid Set-Seeking Systems for Dynamic Plants}\label{ESCStatic}
In the previous section, we presented results and illustrative examples of
hybrid set-seeking systems applied to plants modeled as static maps with
output $y=J(u)$. In Part~III, we extend these results to settings in which
the plant exhibits dynamics. By leveraging singular perturbation results for hybrid inclusions (see Singularly Perturbed Hybrid Dynamical Systems), the standard methodology used to incorporate dynamic plants into smooth ES controllers can be naturally extended to the hybrid systems framework. This extension enables real-time set-seeking in dynamical systems.

\subsection{Model of the Plant}
\label{sec_model_inclusion}
Consider a plant with state $\theta\in\mathbb{R}^m$, input $u\in\mathbb{R}^n$ and output $y\in\mathbb{R}$, modeled by the constrained differential inclusion
\begin{equation}\label{dynamicstheta_1}
\dot{\theta}\in P(\theta,u),~~~~~(\theta,u)\in \Lambda_{\theta}\times\mathbb{U},~~~~~~y=\varphi(\theta,u),
\end{equation}
where $P:\mathbb{R}^m\times\mathbb{R}^n\rightrightarrows\mathbb{R}^m$ is a set-valued mapping, $\varphi:\mathbb{R}^m\times\mathbb{R}^n\rightarrow\mathbb{R}$ is an output function, $\theta$ evolves in the compact set $\Lambda_{\theta}:=\lambda_{\theta}\mathbb{B}\subset\mathbb{R}^m$, $\lambda_{\theta}\in\mathbb{R}_{>0}$ can be taken arbitrarily large to encompass any solution of interest, and $u$ evolves in the closed set $\mathbb{U}:=\hat{\mathbb{U}}+\mathbb{B}$, where $\hat{\mathbb{U}}\subset\mathbb{R}^n$. The unitary inflation on $\hat{\mathbb{U}}$ is motivated by the fact that the control signal $u$ will satisfy $u(t,j)\in\hat{\mathbb{U}}+\varepsilon_a\mathbb{B}\subset\mathbb{U}$ for all $(t,j)$ in the domain of the solutions, where $\varepsilon_a\in(0,1)$ is a tunable parameter.

To have enough regularity in the closed-loop system, we assume that $\hat{\mathbb{U}}$ is a closed set, that $P(\cdot,\cdot)$ is OSC, LB, and convex-valued relative to $\mathbb{R}^m\times\mathbb{U}$, and that $\varphi(\cdot,\cdot)$ is continuously differentiable. In this sense, the model \eqref{dynamicstheta_1} is quite general as it captures differential equations with a continuous right-hand side, as well as discontinuous plants modeled by their corresponding Krasovskii regularizations (see "Krasovskii Solutions of ODEs" in "Continuous-Time Set-Valued Dynamical Systems"). Common examples of plants with a discontinuous right-hand side include mechanical systems with Coulomb friction \cite[Chapter 12]{DiscontinuousSystemsOrlov}, systems arbitrarily switching between a finite number of continuous vector fields \cite[Example 2.14]{bookHDS}, and plants with uncertain models and internal discontinuous feedback controllers, such as those based on sliding mode control \cite{KhalilBook}.

\subsection{Model-Free Feedback Optimization}
For dynamic plants, the goal in ES is to regulate the input $u$ towards the set of points that optimizes the steady-state input-to-output map of the plant. To ensure a well-defined optimization problem, we assume that the plant \eqref{dynamicstheta_1} exhibits suitable stability properties for each fixed input \( u \). This can be achieved by first designing an inner controller that stabilizes \eqref{dynamicstheta_1} with respect to an external reference. Specifically, we assume that there exists a set-valued mapping $H:\mathbb{R}^{n}\rightrightarrows\mathbb{R}^{m}$ that is OSC and LB relative to $\mathbb{U}$, such that $H(\mathbb{U})\subset \Lambda_{\theta}$, and for each $\rho>0$ the compact set  
\begin{equation}
\mathbb{M}_\rho:=\{(\theta,u):~\theta\in H(u), u\in \mathbb{U}\cap \rho\mathbb{B}\} 
\end{equation}
is UGAS for the dynamics
\begin{equation}\label{restrictedplantdynamics}
(\theta,u)\in\Lambda_{\theta}\times(\mathbb{U})\cap\rho\mathbb{B},~~~\dot{\theta}\in P(\theta,u),~~~\dot{u}=0.
\end{equation}
In words, the above property establishes that for each fixed input $u$, the plant dynamics \eqref{dynamicstheta_1} have a well-defined (potentially set-valued) quasi-steady state, characterized by the set-valued mapping $H$. This property generalizes the classic assumptions made for ES in differential equations, e.g., \cite[Assumptions 1 and 2]{TanAndNesic2006Local}, for the case that the plant under control is given by a constrained set-valued dynamical system. %Note that we do not assume that $H(u)\subset\Lambda_{\theta}$ for all $u\in\mathbb{U}$, since for the closed-loop system we will always restrict $u$ to lie on a compact set that can be selected arbitrarily large to encompass all complete solutions of interest. Once $u$ has been restricted to a compact set, the constant $\lambda_{\theta}$ can be selected large enough to guarantee the containment $H(u)\subset\Lambda_{\theta}$.   
\begin{example}[Oscillator with Coulomb Friction]\label{Example1}
Consider a simple harmonic oscillator with Coulomb friction and controllable velocity offset, given in open loop by the discontinuous dynamics
\begin{equation}\label{Example1_Equation}
\dot{\theta}_{1}=\theta_2-u,~~~~\dot{\theta}_{2}=-\frac{B}{M}\text{sgn}(\theta_2-u)-\frac{K}{M}\theta_1,~~~~\dot{u}=0,
\end{equation}\noindent
with $\Lambda_{\theta}=\lambda_{\theta}\mathbb{B}\subset\mathbb{R}^2$, $\mathbb{U}=\mathbb{R}$, $(B,K,M)\in\mathbb{R}_{>0}^{3}$, and $\lambda_{\theta}>\frac{B}{K}>0$ is selected sufficiently large to characterize all the complete solutions of interest. The Krasovskii regularization of this system affects only the dynamics of $\theta_2$, and is given by
\begin{equation}
\dot{\theta}_{2}\in P_2(\theta,u)=\left\{\begin{array}{l}
-\frac{B}{M}-\frac{K}{M}\theta_1~~~~~~~~~\text{if}~~\theta_2>u~\\
\left[-\frac{B}{M},\frac{B}{M}\right]-\frac{K}{M}\theta_1~~\text{if}~\theta_2=u\\
\frac{B}{M}-\frac{K}{M}\theta_1~~~~~~~~~~~~\text{if}~~\theta_2<u.
\end{array}\right.
\end{equation}\noindent
In this case, for each $\rho>0$, system \eqref{restrictedplantdynamics} renders the set $\mathbb{M}_{\rho}:=\{(\theta,u):\theta \in [-\frac{B}{K},\frac{B}{K}]\times\{u\}, u\in\mathbb{R}\cap\rho\mathbb{B}\}$ UGAS \cite[Section 3]{Paden}, thus satisfying the assumptions made on \eqref{dynamicstheta_1}. \QEDB
\end{example}
\begin{example}[Fast Switching Plant]\label{Example2}
Consider an open-loop switched linear system given by $\dot{\theta}=p_q(\theta,u):=A_q\theta+B_qu$, $\dot{u}=0$, $\mathbb{U}:=\mathbb{R}$, $(\epsilon,\lambda_{\theta})\in\mathbb{R}^2_{>0}$, with matrices
 \begin{align}\label{DynamicsExample2}
 A_q=\left[\begin{array}{cc} -1 & \frac{3}{2}-\frac{5}{4}q \\ -\frac{9}{4}+\frac{5}{4}q & -1 \end{array}\right],~ B_q=\left[\begin{array}{c} 1 \\\frac{9}{4}-\frac{5}{4}q  \end{array}\right],
\end{align}\noindent
where $q\in\{1,2\}$. Under arbitrarily fast switching of $q$, this system is conveniently modeled by the differential inclusion 
\begin{equation}\label{exampleplantswitching}
\dot{\theta}\in\{\alpha p_1(\theta,u)+(1-\alpha)p_2(\theta,u), \alpha\in[0,1]\},~\dot{u}=0, 
\end{equation}
which is OSC, LB and convex-valued. %Moreover, using the Lyapunov function $V=(\theta_1-u)^2+\theta_2^2$,
For each $\rho>0$ system \eqref{exampleplantswitching} with flow set $\Lambda_{\theta}\times\left(\mathbb{U}\cap\rho\mathbb{B}\right)$
renders UGAS the set  $\mathbb{M}_{\rho}:=\{(\theta,u):\theta \in \{u\}\times\{0\}, u\in\mathbb{R}\cap\rho\mathbb{B}\}$, thus satisfying the required assumption. \QEDB
\end{example}

To ensure a well-defined single-valued steady-state model-free optimization problem, we also assume that for each \( u \in \mathbb{U} \) and any pair \( \theta, \theta' \in H(u) \), the condition \( \varphi(\theta, u) = \varphi(\theta', u) \) holds. In Example \ref{Example1}, any continuous output function \( \varphi(\cdot, \cdot) \) that depends only on \( \theta_2 \) in its first argument will satisfy this assumption. On the other hand, in Example \ref{Example2},  any continuous output function satisfies this assumption. If $H$ is singled-valued, as in the traditional literature of ES, then this assumption is also satisfied.

We can now define the \emph{response map} of the plant \eqref{dynamicstheta_1}: 
\begin{equation}\label{response_map}
J(u):=\{\varphi(\theta,u):~ \theta\in H(u)\},
\end{equation}
which is assumed to be continuously differentiable. As in the static case, the set-seeking problem is given by
\begin{equation}\label{optimization_problem2} \text{optimize}~~J(u),~~~\text{s.t.}~~u\in \hat{\mathbb{U}}. \end{equation}
For consistency, we denote again the set of solutions to \eqref{optimization_problem2} as $\mathcal{O}\subset\mathbb{R}^n$, which is assumed to be non-empty and compact.
\subsection{Closed-Loop System and Stability Properties}
To interconnect the plant \eqref{dynamicstheta_1} with the hybrid-set seeking dynamics \eqref{originalsystemcomplete01} via the feedback law \eqref{main_input}, a time-scale separation usually needs to be induced between the flows of both systems. To achieve this separation, the flow map \eqref{flowmapcomplete} is multiplied by a small gain $k>0$, leading to the following closed-loop hybrid system:
\begin{subequations}\label{originalsystemcomplete}
 \begin{align}
 C:&=C_{u,z}\times \Lambda_{\xi}\times\mathbb{S}^n\times\Lambda_{\theta},\label{flowsetcompleteFinal}\\
 \dot{x}\in F(x):&= \left(\begin{array}{c}
 k\hat{F}_{\delta}(x_{u,z},\xi)\\
 -k\frac{k_f}{\varepsilon_{f}}~ \left(\xi-\frac{2}{\varepsilon_a} \varphi(\theta,\hat{u}+\varepsilon_a\mathbb{D}\mu)\cdot\mathbb{D}\mu\right)\\
k\Phi(\omega)\mu\\
P(\theta,~\hat{u}+\varepsilon_a\mathbb{D}\mu)\\
\end{array}\right),\label{flowmapcomplete3}\\
D:&=D_{u,z}\times \Lambda_{\xi}\times\mathbb{S}^n\times \Lambda_{\theta},\label{jumpsetcomplete3Final}\\
x^+\in G(x):&=\left(  \begin{array}{c} 
\hat{G}_{\delta}(x_{u,z})\\
\xi\\
\mu\\
\theta
  \end{array}\right),\label{jumpmapcompleteFinal}
\end{align}
\end{subequations}\noindent
where $x:=(x_{u,z},\xi,\mu,\theta)\in\mathbb{R}^{r+4n+m}$, with tunable parameters $(k,\varepsilon_a,\varepsilon_f,\varepsilon_{\omega})$. Figure \ref{fig112f} shows a high-level block diagram of the interconnection between the plant and the controller. The dependence of the hybrid dynamics on the output \( y \) highlights that the specific structure of the components is application-dependent.

To study the stability properties of the closed-loop system \eqref{originalsystemcomplete} via singular perturbation theory, we introduce a new continuous-time variable given by $s=kt$. Since $\frac{d}{dt}x=k\frac{d}{ds}x$, the continuous-time dynamics \eqref{flowmapcomplete3} in the new $s$-time scale become
\begin{equation}
\frac{d}{ds}x\in F(x):= \left(\begin{array}{c}
 \hat{F}_{\delta}(x_{u,z},\xi)\\
 -\frac{k_f}{\varepsilon_f}~ \left(\xi-\frac{2}{\varepsilon_a} \varphi(\theta,\hat{u}+\varepsilon_a\mathbb{D}\mu)\cdot\mathbb{D}~\mu\right)\\
\Phi(\omega)\mu\\
\frac{1}{k}P(\theta,~\hat{u}+\varepsilon_a\mathbb{D}~\mu)\\
\end{array}\right),\label{flowmapcomplete3ave}
\end{equation}
which, combined with \eqref{flowsetcompleteFinal}, \eqref{jumpsetcomplete3Final}, and \eqref{jumpmapcompleteFinal} generates a singularly perturbed hybrid system of the form \eqref{SPSHDS1} with $\varepsilon=k$ (see "Singularly Perturbed Hybrid Dynamical Systems"). We denote this system as $\mathcal{H}_{s}$ to emphasize that the continuous-time dynamics evolve in the $s$-time scale. For this system, it follows that the reduced hybrid dynamics \eqref{reducedHDS} are precisely given by the hybrid-set seeking systems \eqref{originalsystemcomplete01} with $y=J(u)$, whose stability properties were studied in Theorems \ref{theorem1}-\ref{theorem4}. Since the boundary layer dynamics of $\mathcal{H}_{s}$ are precisely given by the plant dynamics \eqref{restrictedplantdynamics}, we can use Theorem \ref{theorem3SP} (see "Singularly Perturbed Hybrid Systems") to obtain the following stability result for hybrid set-seeking systems with plants in the loop \cite{PoTe17Auto}: 
\begin{figure*}[t!]
\begin{tcolorbox}[colback=ivoryA, colframe=ivoryA]
\centering
\includegraphics[width=0.95\linewidth]{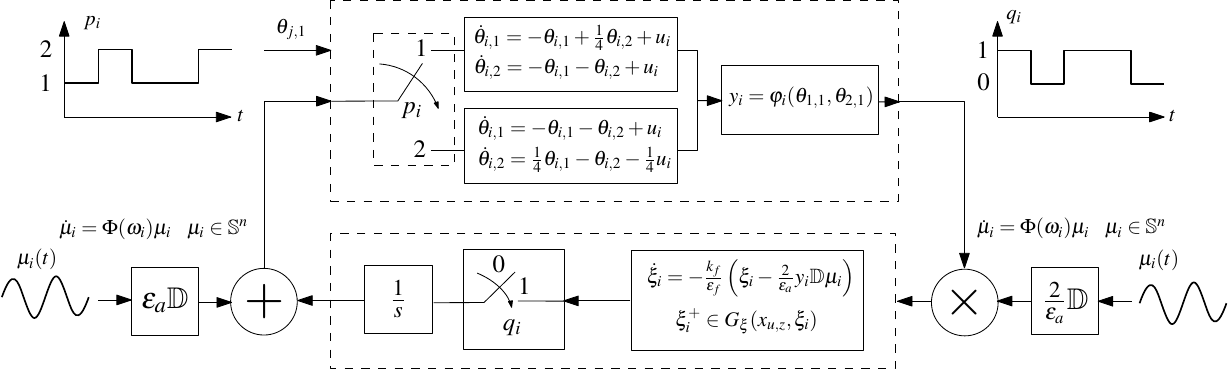}
\end{tcolorbox}
\caption{\tcr{Block diagram representation of Nash-seeking dynamics under intermittent updates for players characterized by switching dynamical systems. The logic mode $q_i$ indicates if the dynamics $\hat{u}_i$ are inactive $(\dot{\hat{u}}_i=0)$ or active $(\dot{\hat{u}}_i=\xi)$}.
\label{NashScheme}}
\end{figure*}

\begin{theorem}[Stability Properties in Dynamic Maps]\label{dynamictheorem}
Consider the HDS $\mathcal{H}_{s}$, and suppose that Assumption \ref{keyassumptions} holds for the hybrid reduced set-seeking dynamics \eqref{originalsystemcomplete01}. Then, the set $\mathcal{A}\times\Lambda_{\xi}\times\mathbb{S}^n\times\Lambda_{\theta}$ is SGPAS as $(\delta,\varepsilon_f,\varepsilon_{a},\varepsilon_{\omega},k)\to 0^+$.
\end{theorem}

Based on Theorem \ref{dynamictheorem} and the material presented in Parts II-III of this article, several comments are in order:  
\begin{itemize}
\item The semi-global practical stability properties of Theorem \ref{dynamictheorem} are established with respect to compact sets, where the main state $x_{u,z}$ is steered towards a neighborhood of $\mathcal{A}=\mathcal{O}\times \Psi$, i.e., $\hat{u}$ converges to a neighborhood of the optimal solution of problem \eqref{optimization_problem2}. Since the actual input $u$ is defined via the feedback law \eqref{main_input}, it will also converge to a neighborhood of order $\varepsilon_a$ of the set $\mathcal{O}$, which is consistent with the classic results in ES.
\item It is possible to refine the convergence statement for the state $\theta$, such that $\theta(t,j)\to H(\mathcal{O}) $ as $(t,j)\to\infty$. Without additional structure or regularity conditions, this property holds only if jumps do not dominate the hybrid system's dynamics---that is, the system must not exhibit Zeno behavior or purely discrete solutions. Such pathological behaviors can be ruled out by applying any of the temporal regularization strategies discussed in Part 2 using dynamic timers, which impose dwell-time or average dwell-time constraints on the jump times.
\item Similarly to the static map case, the result of Theorem \ref{dynamictheorem} asserts suitable $\mathcal{KL}$ bounds of the form \eqref{KLbound2} for all solutions of the system. In this sense, we do not insist on the uniqueness of solutions for any of the proposed controllers. Instead, we ask that every solution satisfies the desired property. When multiple solutions exist, the practitioner can select the most meaningful solution for the application of interest based on application-dependent considerations.
\item As discussed in Part II, the result of Theorem \ref{dynamictheorem} can be extended to cover set-seeking systems that use hybrid filters and hybrid dither generators as part of the feedback loops. Similarly, more complex plants can also be incorporated by using averaging tools for hybrid systems that exhibit hybrid boundary layer dynamics. Due to space constraints, we do not elaborate on these extensions here, and instead we refer readers interested in these subjects to the references \cite{SP_HDS_hybrid_boundary_layer}, \cite{abdelgalil2023multi}, \cite{Kutadinata_Moase}, \cite{Kutadinata:14_Traffic}, and \cite{PovedaCDC18}.
\end{itemize}

\tcr{We finish this section with a numerical example in the context of Nash equilibrium seeking problems, and by highlighting the fact that the stability results presented in  Theorems \ref{theorem1}, 4 and \ref{dynamictheorem} are amenable to standard discretization techniques, such as (consistent) Euler methods or Runge-Kutta methods (c.f., Corollary 2). This follows by the fact that the proposed hybrid controllers are robust to small perturbations, including those that emerge when the flow map is discretized. For further details on this subject, we refer the reader to \cite{Simulator_HDS} and \cite[Sec.4]{PovedaNaliAuto20}.}

\begin{figure*}
\begin{tcolorbox}[colback=iceblue!80, colframe=iceblue!80]
\centering
\includegraphics[width=0.7\linewidth]{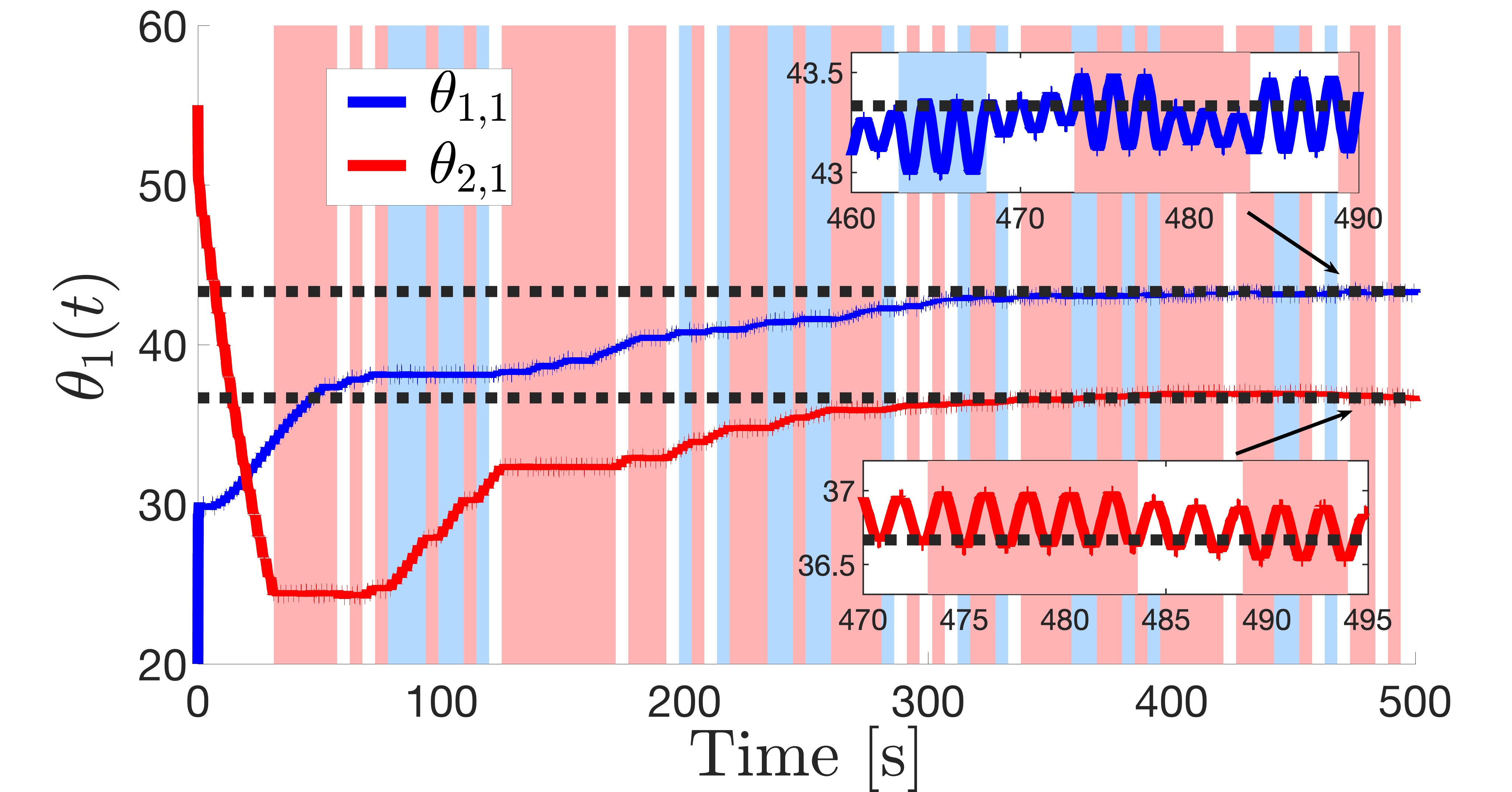}
\end{tcolorbox}
\caption{Evolution in time of the states $\theta_{i,1}$, for each player $i\in\{1,2\}$, under intermittent updates of the Nash seeking dynamics. The colors (light blue and light red) of the shaded regions indicate times when $q_i=0$ for each player $i$.\label{figNashSimulation}}
\end{figure*}

\begin{example}[Nash-Seeking with Intermittent Updates]\label{nashseekingballexample}
\tcr{Consider a two-player noncooperative game, where each player $i\in\{1,2\}$ is characterized by a switching plant with state $\theta_i=(\theta_{i,1},\theta_{i,2})$, and dynamics of the form \eqref{DynamicsExample2}, that is} 
\begin{align}
\dot{\theta}_{i,1}&=-\theta_{i,1}+\left(\frac{3}{2}-\frac{5}{4}p_i\right)\theta_{i,2}+u_i,\\ \dot{\theta}_{i,2}&=-\theta_{i,2}+\left(-\frac{9}{4}+\frac{5}{4}p_i\right)\theta_{i,1}+\left(\frac{9}{4}-\frac{5}{4}p_i\right)u_i, 
\end{align}
\tcr{where $p_i:\text{dom}(q)\to\{1,2\}$ is the switching signal, which, for simplicity, is assumed to be the same for both players. Each player has a quadratic output function $\varphi_i$ that depends on the first state of each player, i.e.,   $y_i=\varphi_i(\tilde{\theta}_1)=\tilde{\theta}_1^\top Q_i\tilde{\theta}_1+\tilde{\theta}_1^\top b_i+c_i$, where $\tilde{\theta}_1=(\theta_{1,1},\theta_{2,1})$, $Q_i\in\mathbb{R}^{2\times2},b_i\in\mathbb{R}^{2}$ and $c_i\in\mathbb{R}$ are selected as in \cite[Sec.6.3]{PoTe17Auto}. The goal of the players is to individually maximize their own steady-state input-to-output cost function \eqref{response_map} using only individual output measurements $y_i$. However, unlike traditional Nash-seeking problems where players continuously update their actions \cite{Frihauf12a,poveda2022fixed}, we study the scenario where the inputs to the player dynamics are updated only sporadically, which is common under computational or communication constraints. To model this behavior, we assign to each player a logic state $q_i\in\{0,1\}$ that multiplies the dynamics of the player's states $\hat{u}_i$, $i\in\{1,2\}$. In particular, we consider the simple Nash-seeking rule $\dot{\hat{u}}_i=q_i\xi_i$, where $\xi_i$ is generated by an individual low-pass filter implemented by each player. Figure \ref{NashScheme} shows a block-diagram representation of the closed-loop system of each player $i\in\{1,2\}$ under the seeking dynamics with switching logic states $q_i$, and switching plant with logic states $p_i$. To study the stability properties of the overall switching system, we let $\ell\in\{1,2,3,4\}$ be a logic state that indexes the four possible combinations of the pair $(q_1,q_2)$, where $\ell=1$ corresponds to $q_1=q_2=1$. Using $\ell\in Q=\{1,2,3,4\}$ as a logic state, and $Q=Q_s\cup Q_u$ where $Q_s=\{1\}$ and $Q_u=\{2,3,4\}$, we impose an average dwell-time bound of the form \eqref{averagedwellinequality} on the jump times of $\ell^+$ using the hybrid automaton \eqref{dwell_time_dynamics}, and an average activation bound \eqref{time_ratio_constraint} using the hybrid monitor \eqref{time_ratio_constraint}. Additionally, since the switching plant admits a common Lyapunov function, we leverage proposition \eqref{proposition_arbitrary} and consider the differential inclusion \eqref{arbitary_switching} obtained from the two modes of the plant. The resulting system has the form \eqref{originalsystemcomplete}, with main states $x_{u,z}=(\hat{u},\ell,\tau_{1},\tau_{2})$. Figure \ref{figNashSimulation} shows the trajectories of the state variables \( \theta_1 \) for the players under the Nash-seeking dynamics with intermittent updates. The red shaded regions indicate time intervals during which player 2 does not update its action, while the blue shaded region corresponds to periods when player 1 does not update \( \hat{u}_1 \). As observed, for a moderate value of \( \eta_2 \), the trajectories still converge to a neighborhood of the Nash equilibrium, indicated by the dotted lines. However, the convergence time increases as \( \eta_2 \to 1 \). For related results in scenarios where \( q_i \) is allowed to take negative values—causing the nominal Nash-seeking dynamics to become unstable for bounded periods of time (i.e., in the presence of adversarial players)—we refer the reader to \cite[Sec.6.3]{PoTe17Auto}.
} 

\end{example}

% \begin{figure}[t!]
% \centerline{\includegraphics[width=0.99\linewidth]{CommunicationNash.jpg}}
% \caption{Nash seeking\label{fig112f}}
% \end{figure}

\section{CONCLUDING REMARKS}
Hybrid set-seeking systems extend traditional continuous-time extremum-seeking algorithms by enabling feedback schemes, algorithms, and plants that incorporate both continuous-time and discrete-time dynamics, and whose stability properties are studied with respect to general sets as opposed to equilibrium points. These systems naturally arise in applications that require switching between multiple controllers to meet stringent transient or steady-state specifications that cannot be achieved with smooth feedback alone. They are also relevant for "intelligent" controllers that use "if-then" rules. Hybrid set-seeking systems may also emerge when the seeking dynamics are smooth, but the plant under control exhibits hybrid behavior, such as jumps or impacts. Other potential applications include developing general, real-time, model-free optimization algorithms for cyber-physical systems, where digital components are integrated into physical plants, such as in sampled-data systems and event-triggered control. The analytical development of hybrid set-seeking systems builds on the traditional tools used for extremum-seeking control: averaging and singular perturbation theory. By extending these tools to hybrid systems, new families of model-free hybrid set-seeking systems can be constructed to emulate model-based hybrid decision-making algorithm. In this sense, working with hybrid time domains and dynamics that satisfy appropriate regularity conditions allows for the property of ``closeness of solutions", which is crucial in the study of smooth ES systems, to be established for hybrid set-seeking systems as well. By leveraging both this property and the uniform stability of the target dynamics, semi-global stability can be demonstrated for general hybrid set-seeking dynamics.

The examples of hybrid set-seeking systems presented in this paper only begin to reveal the broad potential of hybrid extremum seeking and hybrid set-seeking in model-free control and optimization. These approaches open the door to novel controllers and algorithms capable of incorporating "if-then" logic through hybrid systems tools. Furthermore, such algorithms are expected to benefit from techniques developed in related fields, including computer science, signal processing, and communications.

The approaches studied in this paper rely on conventional averaging techniques, often referred to as "first-order" averaging. However, over the past decade, a rich body of work in extremum seeking has emerged based on Lie-bracket and second-order averaging techniques, see \cite{DurrLieBracket,TrackingES,labar2022extremum,labar2019newton,michalowsky2014multidimensional,labar2022iss,Grushkovskaya18,abdelgalil2024initialization}. These methods are often better suited for applications involving geometric constraints, such as systems evolving on manifolds or non-holonomic systems. For a monograph on Lie-bracket-based extremum seeking with applications to model-free stabilization, we refer the reader to \cite{NonSmoothESC_Krstic}. For hybrid tools similar to those discussed in this paper, which are applicable to Lie-bracket averaging-based systems, we refer the readers to \cite{abdelgalil2023lie,abdelgalil2024hybrid}.

Finally, in addition to the development of these novel algorithms, several open problems remain in the context of hybrid seeking systems. These include the development of a comprehensive stability framework for seeking systems that incorporate stochastic behaviors during flows and jumps, as well as developing tools for interconnected networked seeking systems, i.e., ``open extremum-seeking systems''---an under-explored area with numerous potential applications in complex engineering, biological, and socio-technical systems where multiple individual agents implement different adaptation and optimization algorithms for real-time decision making.

\section{Acknowledgments}
Many of the developments reported in this article are a
result of research supported by the Air Force Office of
Scientific Research under grant number FA9550-22-1-0211, and the National Science Foundation under grant numbers CNS 1947613, ECCS CAREER 2305756, and CMMI 2228791.

\section{Author Information}

\begin{IEEEbiography}{{J}orge I. Poveda}{\,}(poveda@ucsd.edu) received B.S. degrees in Electronics Engineering and Mechanical Engineering, both from the University of Los Andes, Bogot\'a, Colombia, in 2012. He received M.Sc. and Ph.D. degrees in Electrical and Computer Engineering from the University of California, Santa Barbara, in 2016 and 2018, respectively. Subsequently, he was a Postdoctoral Fellow at Harvard University in 2018 and Assistant Professor at the University of Colorado (Boulder). Since 2022, he has been with the Department of Electrical and Computer Engineering at the University of California, San Diego, where he is currently Associate Professor.  He is the recipient of the CRII  and CAREER awards from NSF, the Young Investigator awards from AFOSR and SHPE, the 2023 Donald P. Eckman Award from AACC, the 2023 Best Paper Award from IEEE Transactions on Control of Network Systems, and the 2013 CCDC Outstanding Scholar Fellowship  and 2020 Best Ph.D. Dissertation award from UCSB. Furthermore, he is an advisor for students selected as finalists for the Best Student Paper award at the 2024 American Control Conference and as winners of the Young Author Award at the 2024 IFAC Conference on Analysis and Design of Hybrid Systems. He was also a finalist for the Best Student Paper Award at the IEEE Conference on Decision and Control in 2017 (as a student) and 2021 (as a co-author). He has served as Associate Editor for Automatica, Nonlinear Analysis: Hybrid Systems, and IEEE Control Systems Letters.
\end{IEEEbiography}

\begin{IEEEbiography}{Andrew R. Teel}{\,}
(teel@ucsb.edu) received his A.B. degree in Engineering Sciences from Dartmouth College in Hanover, New Hampshire, in 1987 and his M.S. and Ph.D. degrees in Electrical Engineering from the University of California, Berkeley, in 1989 and 1992, respectively. After receiving his Ph.D., he was a postdoctoral fellow at the Ecole des Mines de Paris in Fontainebleau, France. In 1992 he joined the faculty of the Electrical Engineering Department at the University of Minnesota, where he was an assistant professor until 1997. Subsequently, he joined the faculty of the Electrical and Computer Engineering Department at the University of California, Santa Barbara, where he is currently a Distinguished Professor and director of the Center for Control, Dynamical systems, and Computation.  His research interests are in nonlinear and hybrid dynamical systems, with a focus on stability analysis and control design. He has received NSF Research Initiation and CAREER Awards, the 1998 IEEE Leon K. Kirchmayer Prize Paper Award, the 1998 George S. Axelby Outstanding Paper Award, and was the recipient of the first SIAM Control and Systems Theory Prize in 1998. He was the recipient of the 1999 Donald P. Eckman Award and the 2001 O. Hugo Schuck Best Paper Award, both given by the American Automatic Control Council, and also received the 2010 IEEE Control Systems Magazine Outstanding Paper Award.  In 2016, he received the Certificate of Excellent Achievements from the IFAC Technical Committee on Nonlinear Control Systems.  In 2020, he and his co-authors received the HSCC Test-of-Time Award. He is Editor-in-Chief for Automatica, and a Fellow of the IEEE and of IFAC. 
\end{IEEEbiography}

\bibliographystyle{IEEEtran}
\bibliography{references2}

%\section{``About this Issue'' Summary}
%
%This article explores extremum-seeking systems with hybrid dynamics in the loop, combining continuous- and discrete-time behaviors. Such systems can arise in high-performance and cyber-physical applications with digital components or logic-based rules. To study hybrid extremum-seeking systems, we first review the control-theoretic foundations of perturbation theory for hybrid inclusions, focusing on averaging and singular perturbations. These tools support the design and analysis of hybrid extremum-seeking algorithms for robust, model-free feedback optimization of static or dynamic plants, incorporating switching mechanisms, resets, intermittent updates, and slowly varying parameters. The article offers a tutorial-style introduction accessible to readers who are new to hybrid systems.
%
\endarticle

\end{document}